\documentclass[aap,MSNbibl,dvips]{arximspdf}
\usepackage{mathbh}
\usepackage{graphicx}
\doi{10.1214/12-AAP890} %kopijuoti is PTS
\volume{23}
\issue{5}
\pubyear{2013}
\firstpage{1913}
\lastpage{1966}

\makeatletter
\newcommand{\angler}{\rangle}
\newcommand{\anglel}{\langle}
\newcommand{\rrVert}{\Vert}
\newcommand{\rrvert}{\vert}
\newcommand{\llVert}{\Vert}
\newcommand{\llvert}{\vert}
\newcommand{\eqref}[1]{(\ref{#1})}
\newtheorem{theorem}{Theorem}[section]
\newproclaim{assumption}[theorem]{Assumption}
\newtheorem{lemma}[theorem]{Lemma}
\newtheorem{conj}[theorem]{Conjecture}
\newproclaim{remark}[theorem]{Remark}
\newtheorem{prop}[theorem]{Proposition}
\newtheorem{cor}[theorem]{Corollary}
\newtheorem{lemmas}{Lemma}[section]

\newcommand{\N}{\mathbb{N}}
\newcommand{\R}{\mathbb{R}}
\newcommand{\E}{\mathbb{E}}
\renewcommand{\P}{\mathbb{P}}
\newcommand{\eps}{\varepsilon}
\newcommand{\1}{\mathbh{1}}
\newcommand{\ld}{\operatorname{ld}}
\newcommand{\sigmab}{\bar{\sigma}}

\newcommand{\dl}{\delta}
\makeatother

\begin{document}
\begin{frontmatter}

\title{Divergence of the multilevel Monte Carlo Euler method for
nonlinear stochastic differential equations\thanksref{T1}}
\thankstext{T1}{Supported in part by the research project ``Numerical solutions of
stochastic differential equations with nonglobally Lipschitz continuous coefficients''
and by the Collaborative Research Centre $701$ ``Spectral Structures and Topological
Methods in Mathematics,'' both funded by the German Research Foundation.}

\runtitle{Divergence of the MLMC method}

\begin{aug}
\author[A]{\fnms{Martin}~\snm{Hutzenthaler}\ead[label=e1]{hutzenthaler@bio.lmu.de}},
\author[B]{\fnms{Arnulf}~\snm{Jentzen}\ead[label=e2]{ajentzen@math.princeton.edu}}\break
\and
\author[C]{\fnms{Peter~E.}~\snm{Kloeden}\corref{}\ead[label=e3]{kloeden@math.uni-frankfurt.de}}
\runauthor{M. Hutzenthaler, A. Jentzen and P. E. Kloeden}
\affiliation{University of Munich (LMU), Princeton University
and\break
Goethe University Frankfurt am Main}
\address[A]{M. Hutzenthaler\\
LMU BioCenter\\
Department Biologie II\\
University of Munich (LMU)\\
D-82152~Planegg-Martinsried\\
Germany\\
\printead{e1}} %adresu isvedimo komanda gale!
\address[B]{A. Jentzen\\
Program in Applied and\\
\quad Computational Mathematics\\
Princeton University\\
Princeton, New Jersey 08544-1000\\
USA\\
\printead{e2}}
\address[C]{P.~E. Kloeden\\
Institute for Mathematics\\
Goethe University Frankfurt am Main\\
D-60054 Frankfurt am Main\\
Germany\\
\printead{e3}}
\end{aug}

% HISTORY:
\received{\smonth{10} \syear{2011}}
\revised{\smonth{6} \syear{2012}}

% ABSTRACT
%
\begin{abstract}
The Euler--Maruyama scheme
is known to diverge strongly
and numerically weakly
when applied to nonlinear
stochastic differential
equations (SDEs) with
superlinearly growing and
globally one-sided Lipschitz
continuous drift coefficients.
Classical Monte Carlo simulations
do, however, not suffer from this
divergence behavior of Euler's
method because this divergence
behavior happens on
\textit{rare events}.
Indeed,
for such nonlinear SDEs
the classical
Monte Carlo Euler method has been
shown to converge by exploiting that the Euler approximations
diverge only on events whose probabilities
decay to zero very rapidly.
Significantly more efficient
than the classical Monte Carlo Euler
method is the recently introduced
multilevel Monte Carlo Euler method.
%in [Giles, M.~B., Multilevel Monte Carlo
%Path Simulation, \textit{Operations Research}
The main observation of this article is
that
this multilevel Monte Carlo Euler
method does---in contrast to classical Monte
Carlo methods---not converge in general
in the case of such nonlinear SDEs.
More precisely, we establish divergence of the
multilevel Monte Carlo Euler method
for a family of SDEs with superlinearly
growing and globally one-sided Lipschitz
continuous drift coefficients.
In particular, the multilevel Monte
Carlo Euler method diverges
for these nonlinear SDEs
on an event that is not at all
rare but has
\textit{probability one}.
As a consequence for applications, we recommend
not to use the multilevel Monte Carlo Euler method
for SDEs with superlinearly
growing nonlinearities. Instead we
propose to combine
the multilevel Monte Carlo method with a
slightly modified Euler method.
More precisely,
we show that
the multilevel Monte
Carlo method
combined with
a
tamed Euler method
converges
for nonlinear SDEs with globally one-sided Lipschitz continuous
drift coefficients
and preserves
its strikingly higher order convergence
rate from the Lipschitz case.
\end{abstract}

% KEYWORDS
\begin{keyword}[class=AMS]
\kwd{60H35}
\end{keyword}
\begin{keyword}
\kwd{Rare events}
\kwd{nonlinear stochastic differential equations}
\kwd{nonglobally Lipschitz continuous}
\end{keyword}

\end{frontmatter}

%s1 #&#
\section{Introduction}
\label{secintro}
We consider the following setting
in this introductory
section.
Let $ T \in(0,\infty) $,
$
d, m \in\mathbb{N} :=  \{ 1, 2, \ldots \}
$,
let
$  ( \Omega, \mathcal{F},
\mathbb{P}  ) $
be a probability space
with a normal filtration
$ ( \mathcal{F}_t )_{ t \in[0,T] } $
and let $ \xi\dvtx\Omega
\rightarrow\mathbb{R}^d $ be an
$ \mathcal{F}_0
/ \mathcal{B}(\mathbb{R}^d) $-measurable
mapping with
$
\mathbb{E} [
\| \xi\|^p_{ \mathbb{R}^d }
] < \infty
$
for all
$ p \in[1,\infty) $.
Moreover, let
$ \mu\dvtx\mathbb{R}^d \rightarrow
\mathbb{R}^d $ be a smooth
globally one-sided Lipschitz continuous
function with at most polynomially
growing derivatives, and let
$ \sigma\dvtx\mathbb{R}^d
\rightarrow\mathbb{R}^{ d \times m }
$ be a smooth globally Lipschitz
continuous function with at most
polynomially growing derivatives.
In particular, we assume that
there exists a real number
$ c \in(0,\infty) $ such that
$ \anglel  x-y, \mu( x ) - \mu( y )
\angler_{ \mathbb{R}^d }
\leq c \| x - y \|^2_{ \mathbb{R}^d }
$
and
$ \| \sigma( x ) - \sigma( y )
\|_{ \mathbb{R}^{ d \times m } }
\leq c \| x - y \|_{ \mathbb{R}^d } $
for all $ x, y \in\mathbb{R}^d $.
These assumptions ensure the existence
of an up to indistinguishability unique
adapted stochastic process
$
X \dvtx[0,T] \times
\Omega\rightarrow
\mathbb{R}^d
$
with continuous sample paths
solving the stochastic
differential equation (SDE)
%
%e1 #&#
\begin{equation}
\label{eqSDEintro} d X_t = \mu( X_t ) \,dt + \sigma(
X_t ) \,dW_t,\qquad X_0 = \xi
\end{equation}
for $ t \in[0,T] $; see, for example,
Alyushina~\cite{Alyushina1987},
Theorem~1 in
Krylov~\cite{Krylov1990} or
Theorem 2.4.1
in Mao~\cite{m97}.
The function $ \mu$ is the
drift coefficient,
and the function
$ \sigma$ is the
diffusion coefficient
of the SDE~\eqref{eqSDEintro}.
%%The drift coefficient
%$ \mu$
%is the infinitesimal mean
%% of the process $ X $
%and
%%the diffusion coefficient
%$ \sigma\cdot\sigma^{ * } $
%is the infinitesimal covariance
%matrix of the process $ X $.
Our goal in this
introductory section is then to efficiently compute
the deterministic real number
%
%e2 #&#
\begin{equation}
\label{eqgoal} \mathbb{E} \bigl[ f( X_T ) \bigr],
\end{equation}
where $ f \dvtx\mathbb{R}^d
\rightarrow\mathbb{R} $ is a
smooth function with at most
polynomially growing derivatives.
Note that this question is not treated in
the standard literature in computational
stochastics
(see, e.g.,
Kloeden and Platen~\cite{kp92} and
Milstein~\cite{m95})
%Milstein and Tretjakov~\cite{mt04}
which concentrates on SDEs
with globally Lipschitz continuous
coefficients rather than the
SDE~\eqref{eqSDEintro}.
%since the drift coefficient $ \mu$
%of the SDE~\eqref{eqSDEintro},
%in general, fails to be globally Lipschitz
%continuous.
% which concentrates
% on globally Lipschitz continuous coefficients.
The computation of statistical
quantities of the form
\eqref{eqgoal} for
SDEs with nonglobally Lipschitz
continuous coefficients is an
%crucial issue
%major concern
%fundamental problem
important aspect in financial
engineering, in particular, in
option pricing.
%It seems that SDEs with nonglobally
%Lipschitz nonlinearities better
%models
%In particular, recently developed
%demonstrates
%For instance, in applications from
%option pricing,
%the solution process of the SDE
%models stock prizes, interest rates
%or other type of underlyings.
For details the reader is refereed
%to the introductory section
%in Higham, Mao, Pan and Szpruch~\cite{hmps10}
%and
to the
monographs
Lewis~\cite{l00},
Glasserman~\cite{g04},
Higham~\cite{h04}
and
Szpruch~\cite{s10}.
%Another application involving
%the computation of statistical
%quantities of the solution
%of a SDE is the study of
%groundwater
%flow. For details the reader is refereed
%to Section~2 in~\cite{cgss00}.

In order to simulate
the quantity~\eqref{eqgoal}
on a computer, one has
to discretize both the solution
process
$
X \dvtx[0,T] \times\Omega
\rightarrow\mathbb{R}^d
$
of the
SDE~\eqref{eqSDEintro}
as well as the underlying
probability space
$
( \Omega, \mathcal{F}, \mathbb{P} )
$.
The simplest method
for discretizing the
SDE~\eqref{eqSDEintro}
is the Euler method
(a.k.a.~Euler--Maruyama
method).
More formally, the Euler
approximations
$
Y^N_n \dvtx\Omega
\rightarrow\mathbb{R}^d
$,
$
n \in \{ 0, 1, \ldots, N  \}
$,
$
N \in\mathbb{N}
$,
for the SDE~\eqref{eqSDEintro} are defined
recursively through
$ Y^N_0 :=\xi$
and
%
%e3 #&#
\begin{equation}
\label{eqdefeuler} Y^N_{ n + 1 } := Y^N_n
+ \mu \bigl( Y^N_n \bigr) \cdot\frac{T}{N} +
\sigma \bigl( Y^N_n \bigr) ( W_{{ (n+1) T }/{ N }} -
W_{{ n T }/{ N }
} )
\end{equation}
for all $ n \in \{ 0, 1, \ldots, N  \} $
and all $ N \in\mathbb{N} $.
Convergence of Euler's method both
in the strong as well as in the numerically
weak sense is well known
in the case of globally Lipschitz continuous
coefficients $ \mu$ and $ \sigma$
of the SDE; see, for example, Section 14.1
in Kloeden and Platen~\cite{kp92}
and Section 12
in Milstein~\cite{m95}.
The case of superlinearly growing
and hence nonglobally
Lipschitz continuous coefficients
of the SDE
is more subtle.
Indeed, Theorem~2.1 in
the recent article~\cite{hjk11}
shows in the presence of noise
that Euler's method diverges
to infinity
both in the strong and numerically
weak sense
if the coefficients of the SDE
grow superlinearly; see Theorem~\ref{thmeEdivergence} below
for a generalization hereof.
In this situation,
Theorem~2.1
in~\cite{hjk11} also proves
the existence of events
$
\Omega_N \in
\mathcal{F}
$,
$
N \in\mathbb{N}
$,
and of
real numbers
$ \theta, c \in(1,\infty) $
such that
$
\mathbb{P}[
\Omega_N ]
\geq
\theta^{
( - N^{ \theta} )
}
$
and
$
| Y^N_N( \omega) |
\geq
c^{ ( c^N ) }
$
for all
$
\omega\in
\Omega_N
$,
$
N \in\mathbb{N}
$.
Clearly, this implies the divergence
of absolute moments of the
Euler approximation, that is,
$
\lim_{ N \rightarrow\infty}
\mathbb{E} [
| Y^N_N |^p
] = \infty
$
for all
$ p \in(0,\infty) $.

The classical method for discretizing expectations is the Monte Carlo
Euler method.
Let
$
Y^{ N, k }_n \dvtx\Omega
\rightarrow\mathbb{R}^d
$,
$ n \in \{ 0, 1, \ldots, N  \} $,
$ N \in\mathbb{N} $,
for $ k \in\mathbb{N} $
be suitable independent copies of the
Euler
approximations~\eqref{eqdefeuler}; see Section~\ref
{secdivergenceMCE}
for the
precise definition.
The Monte Carlo Euler approximation
of \eqref{eqgoal} with
$ N \in\mathbb{N} $ time steps
and $ N^2 $
Monte Carlo runs
(see Duffie and Glynn~\cite{dg95b}
for more details on this choice)
is then the random
real number
%
%e4 #&#
\begin{equation}
\label{eqdefmceuler} \frac{1}{N^2} \Biggl( \sum
_{ k = 1 }^{ N^2 } f \bigl( Y^{ N, k }_N
\bigr) \Biggr).
\end{equation}
Convergence of the Monte Carlo Euler approximations~\eqref{eqdefmceuler}
is well known in the case of globally Lipschitz continuous
coefficients $\mu$ and $\sigma$; see, for example, Section 14.1
in Kloeden and Platen~\cite{kp92}
and Section 12
in Milstein~\cite{m95}.
Recently, convergence
of the Monte Carlo
Euler
approximations~\eqref{eqdefmceuler}
has also been established for
the SDE~\eqref{eqSDEintro}.
More formally,
Corollary~3.23
in
\cite{HutzenthalerJentzen2012PhiArxiv}
(which generalizes
Theorem~2.1 in~\cite{hj11})
implies
%
%e5 #&#
\begin{equation}
\label{eqmcconvergence} \lim_{ N \rightarrow\infty} \Biggl\llvert \mathbb{E} \bigl[ f(
X_T ) \bigr] - \frac{ 1 }{ N^2 } \Biggl( \sum
_{ k = 1 }^{ N^2 } f \bigl( Y^{ N, k }_N
\bigr) \Biggr) \Biggr\rrvert = 0
\end{equation}
$\mathbb{P}$-almost surely; see also
Theorem~\ref{thmMCEconvergence} below.
The Monte Carlo Euler method
is thus \textit{strongly consistent}
(see, e.g.,
Nikulin~\cite{n01},
Cram{\'e}r~\cite{c99} or
Appendix~A.1 in
Glasserman~\cite{g04})
for the SDE~\eqref{eqSDEintro}.
The reason why
convergence~\eqref{eqmcconvergence}
of the Monte Carlo Euler
method does hold
although the Euler approximations
diverge is
as follows.
%that the
%events $ \Omega_N $,
%$ N \in\mathbb{N} $,
%on which Euler's
%method diverges
%(see the previous paragraph)
%are \textit{rare
%events} whose probabilities
%decay to zero very rapidly
%(see Lemma~4.6 in~\cite{hj11}
%for details).
%In other words:
The events $ \Omega_N $,
$ N \in\mathbb{N} $,
on which Euler's method diverges
(see Theorem~\ref{thmeEdivergence} below)
are \textit{rare events}
and their probabilities decay
to zero faster than any polynomial
in $ N $
as $ N \to\infty$;
see Lemma 2.6 in~\cite{hjk10b}
for details.
Therefore,
for large
$ N \in\mathbb{N} $
the event
$ \Omega_N $
is too unlikely
to occur in any
of $ N^2 $ Monte Carlo
simulations
in \eqref{eqdefmceuler}.
%See also Section~3 in
%simulations
%illustrating~\eqref{eqmcconvergence}.

% MLMC
Considerably
more efficient than the Monte Carlo
Euler method are the so-called
multilevel Monte Carlo Euler methods
% which has
%, in the context of SDEs,
% been introduced recently by
in Giles~\cite{g08b};
see also
Creutzig et al.~\cite{cdmr09},
Dereich~\cite{dereich2011},
Giles~\cite{g08a},
Giles, Higham and
Mao~\cite{ghm09},
Heinrich~\mbox{\cite{h98,heinrich01}},
Heinrich and Sindambiwe~\cite{hs99}
and
Kebaier~\cite{k05}
for related results.
In this method, time is discretized through
the Euler method
and expectations are approximated
by the multilevel Monte Carlo method.
More formally, let
$
Y^{ N, l, k }_n \dvtx\Omega
\rightarrow\mathbb{R}^d
$,
$
n \in \{ 0, 1, \ldots, N  \}
$,
$ N \in\mathbb{N} $,
for $ l \in\N_0:= \{ 0, 1, 2, \ldots \} $
and $ k \in\mathbb{N} $
be suitable independent copies
of the Euler
approximations~\eqref{eqdefeuler}; see Section~\ref{sectamedconvergence}
for the precise definition.
Then the
multilevel Monte Carlo
Euler approximations
for the SDE~\eqref{eqSDEintro} which we investigate in this article
are defined as
%
%e6 #&#
\begin{eqnarray}
\label{eqmlmc} \frac{1}{N} \sum_{ k = 1 }^{ N }
f \bigl( Y_1^{ 1, 0, k } \bigr) + \sum
_{ l = 1 }^{ \log_2( N ) } \frac{ 2^l }{ N } \Biggl( \sum
_{ k = 1 }^{ { N }/{ 2^l } } f \bigl( Y_{ 2^l }^{ 2^l, l, k }
\bigr) - f \bigl( Y_{ 2^{ (l-1) } }^{
2^{ (l-1) }, l, k
} \bigr) \Biggr)
\end{eqnarray}
for $ N \in\{ 2^1, 2^2, 2^3, \ldots\} $.
Clearly, there are also other multilevel Monte Carlo methods than~(\ref{eqmlmc}); see, for example, Giles~\cite{g08b}
for more details. For simplicity, we refer to~(\ref{eqmlmc}) as the multilevel Monte Carlo Euler method throughout this article.
%
%To be accurate, every Monte Carlo method
%with several levels is a multilevel Monte Carlo method.
%We concentrate on~\eqref{eqmlmc} as it is the simpliest
%multilevel Monte Carlo Euler method having the optimal
%order of convergence.
In the case of globally Lipschitz
continuous coefficients of the
SDE~\eqref{eqSDEintro},
this method has been shown to
converge
significantly faster
to the target quantity~\eqref{eqgoal}
than the Monte Carlo Euler
method~\eqref{eqdefmceuler}.
More precisely, in the case
of globally Lipschitz continuous coefficients
$ \mu$ and $\sigma$,
the multilevel Monte
Carlo Euler method~\eqref{eqmlmc}
converges with
order $ \frac{1}{2}-$ while the Monte Carlo Euler
method converges
with order
$ \frac{1}{3}-$ with respect to
the computational effort; see Section~1 in
Giles~\cite{g08b}
or
Creutzig et al.~\cite{cdmr09}
for details.
In the general setting of the
SDE~\eqref{eqSDEintro} where
$ \mu$ does not need to be
globally Lipschitz continuous,
convergence of the multilevel
Monte Carlo Euler
method~\eqref{eqmlmc} remained
an open question.

%%%%%%%%%%%%%%%%%%%
%
%f1 #&#
\begin{figure}

\includegraphics{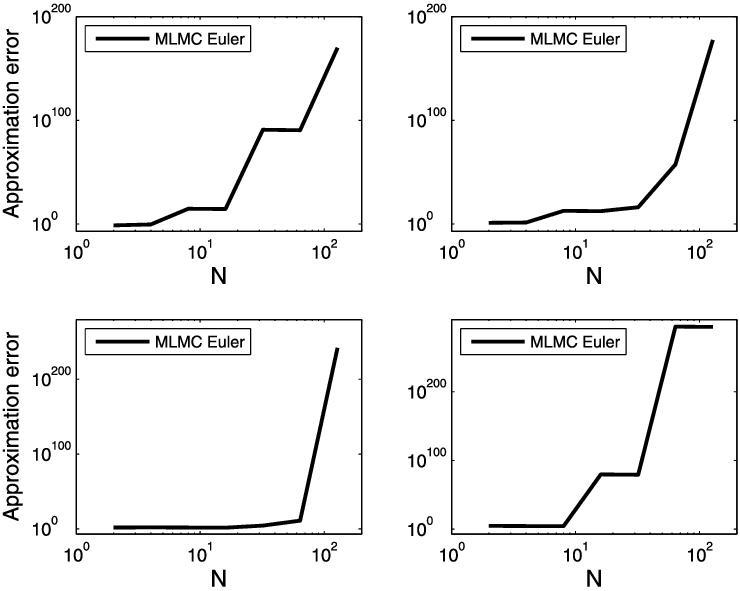}

\caption{Four sample paths of the
approximation error of the multilevel
Monte Carlo Euler
approximation~\protect\eqref{eqmlmc}
for the SDE~\protect\eqref{eqSDEx5}
for $N\in\{2^1,2^2,\ldots,2^{7}\}$
with $ T = 1 $.}
\label{ffivesigma1}
\end{figure}
%
%%%%%%%%%%%%%%%%%%%

The
convergence~\eqref{eqmcconvergence} of the
Monte Carlo Euler method,
and the fact that Euler's method
diverges on very rare events only,
shaped our first guess that the
multilevel Monte Carlo Euler method
should converge too.
However, convergence
of the multilevel Monte Carlo Euler method
fails to hold in the general setting of
the SDE~\eqref{eqSDEintro}.
To prove this, it suffices to establish
nonconvergence for one counterexample
which we choose to be as follows.
Let
$ d = m = 1 $,
let $ \mu( x ) = - x^5 $,
$ \sigma(x) = 0 $,
$ f(x) = x^2 $
for all $ x \in\mathbb{R} $
and let $ \xi\dvtx \Omega
\rightarrow\mathbb{R} $
be standard normally distributed.
Clearly, this choice satisfies
the assumptions of the
SDE~\eqref{eqSDEintro} and
the SDE~\eqref{eqSDEintro}
thus reduces
to the random ordinary differential equation
%
%e7 #&#
\begin{equation}
\label{eqSDEx5} d X_t = - X_t^5 \,dt,\qquad
X_0 = \xi
\end{equation}
for $ t \in[0,T] $.
The main observation of this article is that the
approximation error of the
multilevel Monte Carlo Euler method for the SDE~\eqref{eqSDEx5}
diverges to infinity.
More formally,
Theorem~\ref{thmconj} below
implies
%
%e8 #&#
\begin{eqnarray}
\label{eqdivergencex5} &&\mathop{\lim_{N \rightarrow\infty}}_{
\log_2(N) \in\mathbb{N}
} \Biggl
\llvert \mathbb{E} \bigl[ ( X_T )^2 \bigr] -
\frac{1}{N} \sum_{ k = 1 }^{ N } \bigl(
Y_1^{ 1, 0, k } \bigr)^{ 2 }
\nonumber
\\[-9pt]
\\[-9pt]
\nonumber
&&\hspace*{19pt}\qquad{}- \sum
_{ l = 1 }^{ \log_2( N ) } \frac{ 2^l }{ N } \Biggl( \sum
_{ k = 1 }^{{ N }/{ 2^l } } \bigl( Y_{ 2^l }^{ 2^l, l, k }
\bigr)^{ 2 } - \bigl( Y_{ 2^{ (l-1) } }^{
2^{ (l-1) }, l, k
}
\bigr)^{ 2 } \Biggr) \Biggr\rrvert = \infty
\end{eqnarray}
$\mathbb{P}$-almost surely.
Note that the multilevel Monte
Carlo Euler method diverges
%in \eqref{eqdivergencex5}
on an event that is not rare
but has
\textit{probability one}.
Thus---in contrast to classical
Monte Carlo simulations---the multilevel Monte Carlo Euler
method
is
very sensitive
to the rare events on which Euler's
method diverges in the sense
of Theorem~\ref{thmeEdivergence}
below.
To visualize the
divergence~\eqref{eqdivergencex5},
Figure~\ref{ffivesigma1} depicts
four random sample paths
of the
approximation error of the
multilevel Monte Carlo Euler method~\eqref{eqmlmc}
for the SDE~\eqref{eqSDEx5}
with $ T=1 $ and shows explosion
even for small values
of $N\in\{2^1,2^2,2^3,\ldots\}$.
We emphasize that we are only able
to establish the\vadjust{\goodbreak}
divergence~\eqref{eqdivergencex5}
for the simple SDE~\eqref{eqSDEx5}.
Even in this simple case,
the proof of
the divergence~\eqref{eqdivergencex5} is
rather involved and requires
precise estimates on the speed
of~divergence of Euler's method
for the random
ordinary differential
equation~\eqref{eqSDEx5}
on an appropriate event of
instability; see below
for an outline.

Comparing the convergence result~\eqref{eqmcconvergence}
for the Monte Carlo Euler method
and the divergence result~\eqref{eqdivergencex5}
for the multilevel Monte Carlo
Euler method reveals a remarkable
difference between the classical
Monte Carlo Euler method
and the new multilevel Monte Carlo
Euler method.
The classical Monte Carlo Euler
method applies both
to SDEs with globally Lipschitz continuous coefficients
and to SDEs with
possibly superlinearly growing
coefficients
such as our SDE~\eqref{eqSDEintro}.
The multilevel Monte Carlo Euler method, however,
produces often completely
wrong values
%cannot be used for a
%number of
%interesting
%applications involving
in the case of
SDEs with
superlinearly growing nonlinearities.
This is particularly unfortunate
as SDEs with superlinearly growing
nonlinearities
are
very \mbox{important} in applications; see, for example,~\cite{l00,hmps11,s10}
for applications in financial engineering.
We recommend
not to use the
multilevel Monte Carlo Euler
method for
applications
with such nonlinear SDEs.

Nonetheless, the multilevel Monte Carlo
method can be used for SDEs with
nonglobally Lipschitz continuous
coefficients
when being
combined with a strongly convergent
numerical approximation method.
For example, in~\cite{hjk10b}
the following slight modification
of the Euler method~\eqref{eqdefeuler}
is proposed.
Let
$ Z^N_n
\dvtx\Omega\rightarrow\mathbb{R} $,
$ n \in \{ 0, 1, \ldots, N  \} $,
$ N \in\N$,
be defined recursively
through $ Z^N_0 := \xi$
and
%
%e9 #&#
\begin{eqnarray}\quad
\label{eqtamed} Z^N_{ n + 1 } := Z^N_n
+ \frac{
\mu( Z^N_n )
\cdot
{T}/{N}
}{
1 +
{T}/{N}
\cdot
\|
\mu( Z^N_n )
\|_{ \mathbb{R}^d }
} + \sigma \bigl( Z^N_n \bigr) (
W_{{ (n+1) T }/{ N }
} - W_{{ n T }/{ N }
} )
\end{eqnarray}
for all
$ n \in \{ 0, 1, \ldots, N-1  \} $
and all
$ N \in\N$.
Following~\cite{hjk10b}
we refer to this numerical approximation
as a tamed Euler
method.
%~\eqref{eqtamed}.
Additionally, let
$ Z^{N,l,k}_n \dvtx
\Omega\rightarrow\mathbb{R} $,
$ n \in \{ 0, 1, \ldots, N  \} $,
$ N \in\N$,
for
$ l \in\N_0 $
and
$ k \in\N$
be independent copies
of the tamed Euler
approximations~\eqref{eqtamed}.
In Theorem~\ref{thmtamedconvergence}
below we then prove convergence of the
multilevel Monte Carlo tamed Euler method
for all locally Lipschitz continuous test
functions on the path space whose
local Lipschitz constants grow at most polynomially.
%More formally,
%the pathwise convergence for the functional~\eqref{eqgoal}
%is as follows.
In particular,
Theorem~\ref{thmtamedconvergence}
below
implies
the existence
of finite random variables
$ C_{ \varepsilon} \dvtx
\Omega\rightarrow[0,\infty) $,
$ \varepsilon\in(0,\frac{1}{2}) $,
such that
%
%e10 #&#
\begin{eqnarray}
\label{eqmlmctamedrate} &&\Biggl\llvert \mathbb{E} \bigl[ f( X_T )
\bigr] - \frac{1}{N} \sum_{ k = 1 }^{ N } f
\bigl( Z_1^{ 1, 0, k } \bigr) \nonumber\\
&&\quad{}- \sum_{ l = 1 }^{ \log_2( N ) }
\frac{ 2^l }{ N } \Biggl( \sum_{ k = 1 }^{{ N }/{ 2^l } } f
\bigl( Z_{ 2^l }^{ 2^l, l, k } \bigr) - f \bigl( Z_{ 2^{ (l-1) } }^{
2^{ (l-1) }, l, k
}
\bigr) \Biggr) \Biggr\rrvert\\
&&\qquad \leq \frac{
C_{ \varepsilon}
}{
N^{( {1}/{2} - \varepsilon
)
}
}\nonumber
\end{eqnarray}
for all $ N \in
\{ 2^1, 2^2, 2^3, \ldots\} $
and all
$ \varepsilon\in(0,\frac{1}{2}) $
$ \mathbb{P} $-almost surely.
To sum it up, the classical Monte
Carlo Euler method converges
[see \eqref{eqmcconvergence}],
the new multilevel Monte Carlo Euler
method,
in general, fails to converge
[see \eqref{eqdivergencex5}]
and the new multilevel Monte Carlo
tamed Euler method converges
and preserves its striking higher
convergence order from the Lipschitz
case [see \eqref{eqmlmctamedrate}].
Thus, concerning applications,
the message of
this article is
that the multilevel Monte Carlo Euler
method~\eqref{eqmlmc} \textit{needs to be modified}
appropriately when being applied to SDEs
with superlinearly growing nonlinearities.
This is a crucial difference
to the classical
Monte Euler method which has been shown
to converge for such SDEs and
which
does
not need to be modified.
However, when
modified appropriately
[see, e.g., \eqref{eqtamed}],
the multilevel Monte
Carlo method preserves its strikingly higher
convergence order from the global Lipschitz case
and is significantly more efficient
than the classical Monte Carlo
Euler method, even for such nonlinear SDEs.
Thereby, this article motivates
future research in the construction
and the analysis of
``appropriately modified''
numerical approximation
methods.\looseness=1
%and
%``appropriately modified''
%multilevel Monte Carlo
%methods which converge
%even in the case of SDEs
%with superlinearly growing
%nonlinearities.

For the interested reader,
we
now outline the central ideas in the proof
of~\eqref{eqdivergencex5}.
For this we use
the random variables
$
\xi^{ l, k } \dvtx\Omega
\rightarrow\mathbb{R}
$,
$ l \in\mathbb{N}_0 $,
$ k \in\mathbb{N} $,
defined by
$
\xi^{ l, k } :=
Y^{ M, l, k }_0
$
for all $ M \in\mathbb{N} $,
$ l \in\mathbb{N}_0 $,
$ k \in\mathbb{N} $.
Then we note for every
$ M, k \in\mathbb{N} $,
$ l \in\mathbb{N}_0 $
and every
$ \omega\in\Omega$
that
$ | Y^{ M, l, k }_n( \omega) | $
is strictly increasing in
$ n \in\{ 0, 1, \ldots, M \} $
if and only if
$
| \xi^{ l, k }( \omega) | =
| Y^{ M, l, k }_0(\omega) |
>
( 2 M
)^{{1}/{4} }
T^{ - {1}/{4} }
$.
It turns out that
$ | Y^{ M, l, k }_n(\omega) | $
increases in $ n \in\{ 0, 1, \ldots, M \} $
double exponentially
fast for all
$ \omega\in\{
| \xi^{ l, k } | >
( 2 M
)^{{1}/{4} }
T^{ - {1}/{4} }
\} $,
$ l \in\mathbb{N}_0 $
and all $ k, M \in\mathbb{N} $; see Lemma~\ref{lupperesti} and
Corollary~\ref{cordoublegrowth} below
for details.
A central observation in our proof
of the
divergence~\eqref{eqdivergencex5}
is then that
the behavior of the multilevel Monte
Carlo Euler method is dominated by
the highest level that produces
such
double exponentially fast
increasing trajectories. More precisely,
a key step in our proof of
\eqref{eqdivergencex5} is to
introduce the random
variables
$
L_N \dvtx
\Omega\to
\{1,2,\ldots,\ld(N)\}
$,
$N\in\{2^1,2^2,2^3, \ldots\}
$,
by
%
%e11 #&#
\begin{eqnarray}
L_N&:=&\max \biggl( \{ 1 \} \cup \biggl\{ l \in \bigl\{ 1, 2, \ldots,
\ld(N) \bigr\} \dvtx
\nonumber
\\[-8pt]
\\[-8pt]
\nonumber % \Big(
&&\hspace*{57pt}\exists k \in \biggl\{ 1, 2, \ldots,
\frac{N}{2^l} \biggr\} \dvtx \bigl| \xi^{l,k}\bigr | > 2^{{ l }/{ 4 } }
T^{ - { 1 }/{ 4 } } % \Big)
\biggr\} \biggr)
\end{eqnarray}
for all $N\in\{2^1,2^2,2^3,\ldots\}$.
Using the random variables
$ L_N $,
$ N \in\{ 2^1, 2^2,\break 2^3, \ldots\} $,
we now rewrite the
multilevel Monte
Carlo Euler method in
\eqref{eqdivergencex5}
as
%
%e12 #&#
%e13 #&#
\begin{eqnarray}
\label{eqrewriteMLintro} && \frac{1}{N} \sum_{ k = 1 }^{ N }
\bigl( Y_1^{ 1, 0, k } \bigr)^{ 2 } + \sum
_{ l = 1 }^{ \log_2( N ) } \frac{ 2^l }{ N } \Biggl( \sum
_{ k = 1 }^{{ N }/{ 2^l } } \bigl( Y_{ 2^l }^{ 2^l, l, k }
\bigr)^{ 2 } - \bigl( Y_{ 2^{ (l-1) } }^{
2^{ (l-1) }, l, k
}
\bigr)^{ 2 } \Biggr)
\\
&&\qquad = \mathop{\sum_{
l \in\{ 0, 1, \ldots, \log_2(N) \}}}_{
l \neq L_N - 1, L_N
}
\frac{ 2^l }{ N } \sum_{ k = 1 }^{{ N }/{ 2^l } } \bigl(
Y_{ 2^l }^{ 2^l, l, k } \bigr)^{ 2 }
\nonumber
\\[-8pt]
\\[-8pt]
\nonumber
&&\qquad\quad{} - \mathop{\sum
_{
l \in\{ 1, 2, \ldots, \log_2(N) \}
}}_{
l \neq L_N
} \frac{ 2^{l} }{ N } \sum
_{ k = 1 }^{{ N }/{ 2^{l} } } \bigl( Y_{ 2^{ (l-1) } }^{
2^{ (l-1) }, l, k
}
\bigr)^{ 2 }
\\
 \label{eqcentralterms}&&\quad\qquad{} + \frac{ 2^{ (L_N-1) } }{ N } \sum_{ k = 1 }^{ { N }/{ (2^{ (L_N-1) }) } }
\bigl( Y_{ 2^{ (L_N-1) } }^{
2^{ (L_N-1) }, L_N-1, k
} \bigr)^{ 2 }
\\\label{eq15}
&&\qquad\quad{}+
\frac{ 2^{ L_N } }{ N } \sum_{ k = 1 }^{{ N }/{ (2^{ L_N }) } } \bigl(
Y_{ 2^{ L_N } }^{ 2^{ L_N }, L_N, k } \bigr)^{ 2 } %\\ &
\\\label{eq16}
&&\qquad\quad{}-
\frac{ 2^{ L_N } }{ N } \sum_{ k = 1 }^{{ N }/{( 2^{ L_N }) } } \bigl(
Y_{ 2^{ ( L_N -1) } }^{
2^{ ( L_N -1) }, L_N , k
} \bigr)^{ 2 }
\end{eqnarray}
for all $N\in\{2^1,2^2,2^3,\ldots\}$.
Due to the definition of $ L_N $,
$ N \in\{ 2^1, 2^2, 2^3, \ldots\} $,
it turns out that the asymptotic
behavior of the multilevel Monte Carlo
Euler method~\eqref{eqrewriteMLintro}
is essentially
determined by the three
sums
in \eqref{eqcentralterms}--(\ref{eq16}); see inequality~\eqref{eqresultcase1},
estimate~\eqref{eqresultcase2}
and
inequalities~\eqref{eqresultcase3}, \eqref{eqresultcase4}
in the proof of
Theorem~\ref{thmconj} for details.
In order to investigate these three
summands, we---roughly speaking---quantify the value of the largest
summand
in each of the three sums
in \eqref{eqcentralterms}--(\ref{eq16}).
For this we introduce the
random variables
$
\eta_N
\dvtx
\Omega\to[0,\infty)
$
and
$
\theta_N
\dvtx
\Omega\to[0,\infty)
$
for
$ N \in\{ 2^1, 2^2, 2^3, \ldots\} $
by
%
%e14 #&#
\begin{equation}
\eta_N:=\max \biggl\{ \bigl| \xi^{L_N,k} \bigr| \in\mathbb{R} \dvtx k
\in \biggl\{ 1, 2, \ldots, \frac{N}{2^{L_N}} \biggr\} \biggr\}
\end{equation}
and
%
%e15 #&#
\begin{equation}
\theta_N := \max \biggl\{ \bigl| \xi^{ L_N - 1, k} \bigr| \in\mathbb{R} \dvtx
k \in \biggl\{ 1, 2, \ldots, \frac{N}{2^{(L_N-1)}} \biggr\} \biggr\}
\end{equation}
for all
$N\in\{2^1,2^2,2^3,\ldots\}$.
Using the random variables
$ \eta_N \dvtx
\Omega\rightarrow
[0,\infty) $
and
$ \theta_N \dvtx
\Omega\rightarrow
[0,\infty) $
for $ N \in\{ 2^1, 2^2, 2^3,
\ldots\} $ we then distinguish
between three different cases
[see inequality~\eqref{eqresultcase1},
inequality~\eqref{eqresultcase2}
and
inequalities~\eqref{eqresultcase3},
\eqref{eqresultcase4}
below].
First, on the events
$
\{
\eta_N
>
2^{{ (L_N+1) }/{ 4 } }
T^{ - {1}/{4} }
\}
\in
\mathcal{F}
$,
$ N \in\{ 2^1, 2^2, 2^3, \ldots\} $,
the middle sum in
\eqref{eq15}
will be positive with
large absolute value
and will essentially
determine the
behavior of
the multilevel Monte
Carlo Euler
approximations~\eqref{eqrewriteMLintro}; see
estimate~\eqref{eqresultcase1}
for details.
Second, on the events
$
\{
\eta_N
\leq
2^{{ (L_N+1) }/{ 4 } }
T^{ - {1}/{4} }
\}
\cap
\{
\eta_N < \theta_N
\}
\in
\mathcal{F}
$,
$ N \in\{ 2^1, 2^2, 2^3, \ldots\} $,
the sum in
\eqref{eqcentralterms}
will be positive with large absolute
value and will essentially
determine the
behavior of
the multilevel Monte
Carlo Euler
approximations~\eqref{eqrewriteMLintro}; see
inequality~\eqref{eqresultcase2}
for details.
Finally, on the events
$
\{
\eta_N
\leq
2^{{ (L_N+1) }/{ 4 } }
T^{ - {1}/{4} }
\}
\cap
\{
\eta_N > \theta_N
\}
\in
\mathcal{F}
$,
$ N \in\{ 2^1, 2^2, 2^3, \ldots\} $,
the sum in~\eqref{eq16}
will be negative with
large absolute value
and will essentially
determine the
behavior of
the multilevel Monte
Carlo Euler
approximations~\eqref{eqrewriteMLintro}; see
inequalities~\eqref{eqresultcase3}
and \eqref{eqresultcase4}
for details.
This very rough outline of
the case-by-case analysis in
our proof
of \eqref{eqdivergencex5} also
illustrates that the multilevel
Monte Carlo Euler
approximations~\eqref{eqrewriteMLintro}
assume both positive
(first and second case) as well as
negative values
(third case)
with large absolute values.
We add that this
case-by-case analysis argument
in our proof of
\eqref{eqdivergencex5} requires
that the probability that the
random variables
$ \eta_N $ and $ \theta_N $
are close to each other
in some sense
must decay rapidly to zero
as $ N \in\{ 2^1, 2^2, 2^3, \ldots\} $
goes to infinity; see inequality~\eqref{eqAN4}
below.
We verify the above decaying
of the probabilities in
Lemma~\ref{lAN4} below which is
a crucial step in our proof
of \eqref{eqdivergencex5}.
Additionally, we add that the level
$L_N$ is approximately of
order $ \log( \log( N ) ) $
as $ N $ goes to infinity; see Lemma~\ref{lAN1} for the precise
assertion.
In view of the above case-by-case
analysis of the multilevel Monte Carlo
Euler method,
we find it quite remarkable to
observe that the essential behaviour
of the multilevel Monte Carlo
Euler method
in \eqref{eqdivergencex5}
is determined by
the levels around the order
$ \log( \log( N ) ) $
as $ N $ goes to infinity.

The remainder
of this article
is organized as follows.
Theorem~\ref{thmeEdivergence}
in
Section~\ref{secdivergenceeE}
slightly
generalizes the result on strong
and weak divergence of the Euler
method of Hutzenthaler, Jentzen
and Kloeden~\cite{hjk11}.
Convergence
of the Monte Carlo Euler method
is reviewed in
Section~\ref{secdivergenceMCE}.
The main result of this article, that is,
divergence of the multilevel Monte Carlo Euler method
for the SDE~\eqref{eqSDEx5},
is presented and
proved in Section~\ref{secmultilevelpathwisedivergence}.
We believe that the multilevel Monte Carlo Euler method
diverges more generally and formulate this
as Conjecture~\ref{conj}
in Section~\ref{secdivergencemlMCE}.
Section~\ref{sectamedconvergence}
contains our proof of almost
sure and strong convergence
of the multilevel Monte Carlo tamed Euler method
for all locally Lipschitz continuous
test functions on the path space
whose
local Lipschitz constants grow
at most polynomially.
%

%
% %%%%%%%%%%%%%%%%%%%%%%%%%%%%%%%%%%%%%%%%%%%%%%%%%%%%%%%%%%%%%%%%%%
% Section
% %%%%%%%%%%%%%%%%%%%%%%%%%%%%%%%%%%%%%%%%%%%%%%%%%%%%%%%%%%%%%%%%%%
%
%s2 #&#
\section{Divergence of
the Euler method}
\label{secdivergenceeE}
Throughout this section assume that the following
setting is fulfilled.
Let $ T \in(0, \infty) $,
let $ ( \Omega, \mathcal{F}, \mathbb{P} ) $
be a probability space with a
filtration $ ( \mathcal{F}_t )_{ t \in[0,T] } $
and let $ W \dvtx[0,T] \times\Omega
\rightarrow\mathbb{R} $ be a one-dimensional
standard
$ ( \mathcal{F}_t )_{ t \in[0,T] } $-Brownian
motion.
Additionally, let $ \xi\dvtx\Omega\rightarrow\mathbb{R} $ be an $
\mathcal{F}_0 / \mathcal{B}( \mathbb{R} ) $-measurable mapping and let
$ \mu, \sigma\dvtx\mathbb{R} \rightarrow\mathbb{R} $ be two $
\mathcal{B}( \mathbb{R} ) / \mathcal{B}( \mathbb{R} ) $-measurable
mappings. We then define the Euler approximations
$
Y_n^N \dvtx
\Omega\rightarrow\mathbb{R}
$,
$
n \in\{ 0, 1, \ldots, N \} $, $ N \in\N
$,
recursively
by $ Y_0^N := \xi$ and
%
%e16 #&#
\begin{eqnarray}
\label{eqexplicitEuler} Y_{ n + 1 }^N := Y_n^N
+ \mu \bigl( Y_n^N \bigr) \cdot \frac{ T }{ N } +
\sigma \bigl( Y_n^N \bigr) \cdot ( W_{ { ( n + 1 ) T }/{ N } } -
W_{ { n T }/{ N } } )
\end{eqnarray}
for all $ n \in\{ 0, 1, \ldots, N - 1 \} $ and all $ N \in\N$.
The following theorem generalizes
Theorem~2.1 in Hutzenthaler,
Jentzen and
Kloeden~\cite{hjk11}.
%
% %%%%%%%%%%%%%%%%%%%%%%%%%%%%%%%%%%%
% Theorem
% %%%%%%%%%%%%%%%%%%%%%%%%%%%%%%%%%%%
%
%th2.1 #&#
\begin{theorem}[(Strong and
weak divergence of the Euler method)]
\label{thmeEdivergence}
Assume that the above setting
is fulfilled,
and let
$ \alpha, c \in( 1, \infty) $
be real numbers such that
$
\llvert  \mu( x ) \rrvert
+
\llvert  \sigma( x ) \rrvert
\geq
\frac{ \llvert  x \rrvert^{ \alpha}
}{ c }
$
for all $ x \in\mathbb{R} $
with $ | x | \geq c $.
Moreover, assume that
$
\mathbb{P} [
\sigma( \xi) \neq0
] > 0
$ or that there
exists a real number
$ \beta\in(1,\infty)$
such that
$
\P [ | \xi| \geq x  ]
\geq
\beta^{
( -x^\beta)
}
$
for all $x\in[1,\infty)$.
Then there exists a real
number $ \theta\in(1,\infty) $
and a sequence of nonempty
events $ \Omega_N \in\mathcal{F} $,
$ N \in\N$, such that
$
\mathbb{P} [
\Omega_N
]
\geq
\theta^{  ( - N^\theta ) }
$ and
$
\llvert
Y^N_N( \omega)
\rrvert
\geq
c^{
(
(
{ (\alpha+ 1) }/{ 2 }
)^{
N
}
)
}
$
for all $ \omega\in\Omega_N $
and all $ N \in\N$.
In particular, the
Euler
approximations~\eqref{eqexplicitEuler} satisfy
$ \lim_{ N \rightarrow\infty} \mathbb{E} [
| Y_N^N |^p
] = \infty$
for all $ p \in(0,\infty) $.
\end{theorem}

Theorem~\ref{thmeEdivergence}
immediately follows from
Lemmas~\ref{lemmatails} and~\ref{lemmaeEdivergence} below.
More results on Euler's method for SDEs
with possibly superlinearly growing
nonlinearities can, for example, be found
in~\cite{gk96b,g98b,mt04,mt05}
and in the references therein.
%
% %%%%%%%%%%%%%%%%%%%%%%%%%%%%%%%%%%%
% Lemmata
% %%%%%%%%%%%%%%%%%%%%%%%%%%%%%%%%%%%
%
%le2.2 #&#
\begin{lemma}[(Tails
of $ Y^N_1 $,
$ N \in\mathbb{N} $)]
\label{lemmatails}
Assume that the above setting
is fulfilled and let
$
\mathbb{P} [
\sigma( \xi) \neq0
] > 0
$.
Then there
exists a real number
$ \beta\in(1,\infty)$
such that
$
\P [ | Y^N_1 | \geq x  ]
\geq
\beta^{
(
- ( N x )^\beta
)
}
$
for all $ x \in[1,\infty) $
and all $ N \in\mathbb{N} $.
\end{lemma}
\begin{pf}
By assumption we have
$ \mathbb{P} [ | \sigma( \xi)
| > 0  ] > 0 $. Therefore,
there exists a
real number $ K \in( 1, \infty) $
such that
%
%e17 #&#
\begin{eqnarray}
\label{eqdefvartheta} \vartheta := \mathbb{P} \biggl[ \bigl| \sigma( \xi) \bigr| \geq
\frac{ 1 }{ K }, | \xi| + T \bigl| \mu( \xi) \bigr| \leq K \biggr] \in ( 0, \infty) .
\end{eqnarray}
Moreover, we have
\begin{eqnarray*}
&&\mathbb{P} \bigl[ \bigl| Y^N_1 \bigr| \geq x \bigr]\\
&&\qquad = \mathbb{P}
\biggl[ \biggl| \xi+ \mu( \xi) \frac{T}{N} + \sigma( \xi) W_{ {T}/{N} } \biggr|
\geq x \biggr] \\
&&\qquad\geq \mathbb{P} \bigl[\bigl | \sigma( \xi) W_{ {T}/{N} }\bigr | - | \xi|
- T \bigl| \mu( \xi) \bigr| \geq x \bigr]
\\
&&\qquad \geq \mathbb{P} \biggl[\bigl | \sigma( \xi) \bigr| \geq \frac{ 1 }{ K }, | \xi| + T \bigl|
\mu( \xi) \bigr| \leq K , \bigl| \sigma( \xi) W_{ {T}/{N} }\bigr | - | \xi| - T \bigl| \mu( \xi)
\bigr| \geq x \biggr]
\\
&&\qquad \geq \mathbb{P} \biggl[ \bigl| \sigma( \xi)\bigr | \geq \frac{ 1 }{ K }, | \xi| + T \bigl|
\mu( \xi) \bigr| \leq K , \frac{1}{K} | W_{ {T}/{N} } | - K \geq x \biggr]
\end{eqnarray*}
for all $ x \in[1,\infty) $
and all $ N \in\mathbb{N} $.
Definition~\eqref{eqdefvartheta}
and Lemma~4.1
in~\cite{hjk11}
therefore show
\begin{eqnarray*}
\mathbb{P} \bigl[ \bigl| Y^N_1 \bigr| \geq x \bigr] & \geq&
\mathbb{P} \biggl[ \bigl| \sigma( \xi)\bigr | \geq \frac{ 1 }{ K }, | \xi| + T \bigl| \mu(
\xi) \bigr| \leq K \biggr] \cdot \mathbb{P} \biggl[ \frac{1}{K} |
W_{ {T}/{N} } | - K \geq x \biggr] %\\ &
\\
&=& \vartheta \cdot \mathbb{P}
\bigl[ T^{ - {1}/{2} } | W_{ T } | \geq T^{ - {1}/{2} }
N^{ {1}/{2} } K (x + K) \bigr]
\\
& \geq&\hspace*{-1pt} \frac{ \vartheta}{ 4 \sqrt{T} } \cdot \exp \bigl( - T^{-1} N K^2
( x + K )^2 \bigr)\hspace*{-0.5pt}\geq\hspace*{-0.5pt} \frac{ \vartheta}{ 4 \sqrt{T} } \cdot \exp \bigl( - 4
T^{-1} K^4 ( N x )^2 \bigr)
\\
& = &\frac{ \vartheta}{ 4 \sqrt{T} } \cdot \bigl( e^{ 4 T^{-1} K^4 } \bigr)^{
(
-  ( N x  )^2
)
} \geq
\biggl( e^{ 4 T^{-1} K^4 } + \frac{
4 \sqrt{T}
}{
\vartheta
} \biggr)^{
(
- 2  ( N x  )^2
)
}
\end{eqnarray*}
for all $ x \in[1,\infty) $
and all $ N \in\mathbb{N} $.
This completes the proof
of Lemma~\ref{lemmatails}.
\end{pf}
%
%le2.3 #&#
\begin{lemma}
\label{lemmaeEdivergence}
Assume that the above setting\vspace*{1pt}
is fulfilled
and let
$ \alpha, c \in( 1, \infty) $
be real numbers such that
$
\llvert  \mu( x ) \rrvert
+
\llvert  \sigma( x ) \rrvert
\geq
\frac{ \llvert  x \rrvert^{ \alpha}
}{ c }
$
for all $ x \in\mathbb{R} $
with $ | x | \geq c $.
Moreover, assume that there exist
real numbers
$ N_0 \in\{ 0, 1, 2, \ldots\} $,
$\beta\in(1,\infty)$
such that
$
\P [ |Y_{N_0}^N|
\geq x  ] \geq
\beta^{ ( -(N x)^\beta )}
$
for all
$ x\in[1,\infty) $
and all
$
N \in
\mathbb{N} \cap
\{ N_0, N_0 + 1, \ldots\}
$.
Then there exists a real
number $ \theta\in(1,\infty) $
and a sequence of nonempty
events $ \Omega_N \in\mathcal{F} $,
$
N \in
\mathbb{N} \cap
\{ N_0, N_0 + 1, \ldots\}
$,
such that
$
\mathbb{P} [
\Omega_N
]
\geq
\theta^{  ( - N^\theta ) }
$
and
$
\llvert
Y^N_N( \omega)
\rrvert
\geq
c^{
(
(
{( \alpha+ 1) }/{ 2 }
)^{
N
}
)
}
$
for all $ \omega\in\Omega_N $
and all
$
N \in
\mathbb{N} \cap
\{ N_0, N_0 + 1, \ldots\}
$.
In particular, the
Euler
approximations~\eqref{eqexplicitEuler} satisfy
$ \lim_{ N \rightarrow\infty} \mathbb{E} [
| Y_N^N |^p
] = \infty$
for all $ p \in(0,\infty) $.
\end{lemma}
\begin{pf}%[Proof of Lemma~\ref{lemmaeEdivergence}]
Define
real numbers
$ r_N \in[0, \infty) $,
$ N \in\N$, by
%
%e18 #&#
\begin{equation}
\label{eqdefrN} r_N := \max \biggl( c, \biggl( \frac{ 2 N c }{ T }
\biggr)^{{ 2 }/{ ( \alpha- 1 ) }
} \biggr)
\end{equation}
for all $ N \in\N$.
We also use the function
$ \operatorname{sgn}
\dvtx\mathbb{R} \rightarrow\mathbb{R} $ defined by
$ \operatorname{sgn}( x ) := 1 $
for all $ x \in[0, \infty) $
and by $ \operatorname{sgn}( x ) := -1 $
for all $ x \in( -\infty, 0 ) $.
Furthermore, we
define events
$ \Omega_N \in\mathcal{F} $,
$
N \in
\mathbb{N} \cap
\{ N_0, N_0 + 1, \ldots\}
$,
by
%
%e19 #&#
\begin{eqnarray}
\label{eqOmegaNdef} \Omega_N &:=& \Biggl( \bigcap
_{ n = N_0 }^{ N - 1 } \biggl\{ \omega\in\Omega\dvtx
\operatorname{sgn} \Bigl( \mu \bigl( Y_n^N( \omega)
\bigr)  \sigma \bigl( Y_n^N( \omega) \bigr) \Bigr)
\nonumber\\
&&\hspace*{32pt}{}\times \bigl( W_{ { ( n + 1 ) T }/{ N }
}( \omega) - W_{ { n T }/{ N }
}( \omega) \bigr) \geq
\frac{ T }{ N } \biggr\} \Biggr)
\\
&&{}\cap \bigl\{ \omega\in\Omega\dvtx\bigl | Y_{ N_0 }^N( \omega)\bigr |
\geq ( r_N )^{
(
(
{( \alpha+ 1)
}/{ 2 }
)^{ N_0 }
)
} \bigr\}\nonumber
\end{eqnarray}
for all
$
N \in
\mathbb{N} \cap
\{ N_0, N_0 + 1, \ldots\}
$.
In particular, the definition of\break
$
( \Omega_N )_{
N \in
\mathbb{N} \cap
\{ N_0, N_0 + 1, \ldots\}
}
$
implies
%
%e20 #&#
\begin{eqnarray}
\label{eqOmegaNfollow} %\lefteqn{
&&\biggl\llvert \mu \bigl(
Y_n^N( \omega) \bigr) \cdot \frac{ T }{ N } + \sigma
\bigl( Y_n^N( \omega) \bigr) \cdot \bigl(
W_{ { ( n + 1 ) T }/{ N } }( \omega) - W_{ { n T }/{ N } }( \omega) \bigr) \biggr\rrvert
%}
\nonumber
\\[-8pt]
\\[-8pt]
\nonumber
&&\qquad= \frac{ T }{ N } \bigl\llvert \mu \bigl(
Y_n^N( \omega) \bigr) \bigr\rrvert + \bigl\llvert \sigma
\bigl( Y_n^N( \omega) \bigr) \bigr\rrvert \cdot \bigl
\llvert W_{ { ( n + 1 ) T }/{ N } }( \omega) - W_{ { n T }/{ N } }( \omega) \bigr\rrvert
\end{eqnarray}
for all
$ n \in \{ N_0, N_0+1, \ldots,
N-1  \} $,
$ \omega\in\Omega_N $ and
all
$
N \in
\mathbb{N} \cap
\{ N_0, N_0 + 1, \ldots\}
$.\vadjust{\goodbreak}

In the next step let
$
N \in
\mathbb{N} \cap
\{ N_0, N_0 + 1, \ldots\}
$
and
$ \omega\in\Omega_N $
be arbitrary.
We then claim
%
%e21 #&#
\begin{equation}
\label{eqclaim1} \bigl\llvert Y_n^N( \omega) \bigr
\rrvert \geq ( r_N )^{
(  (
{( \alpha+ 1) }/{ 2 }
)^{ n }
)
}
\end{equation}
for all
$ n \in\{ N_0, N_0 + 1, \ldots, N \} $.
We now show \eqref{eqclaim1}
by induction on
$ n \in\{ N_0, N_0 + 1, \ldots, N \} $.
The base case $ n = N_0 $
follows from
definition~\eqref{eqOmegaNdef}
of $\Omega_N$.
For the induction step
assume that~\eqref{eqclaim1}
holds for one
$ n \in\{ N_0, N_0 + 1, \ldots, N - 1
\} $.
In particular, this implies
%
%e22 #&#
\begin{eqnarray}
\label{eqclaim2} \bigl\llvert Y_n^N( \omega) \bigr
\rrvert \geq ( r_N )^{
(
(
{ (\alpha+ 1) }/{ 2 }
)^{ n }
)
} \geq r_N \geq c > 1 .
\end{eqnarray}
Moreover,
definition~\eqref{eqexplicitEuler},
the triangle inequality and
equation~\eqref{eqOmegaNfollow}
yield
\begin{eqnarray*}
 \bigl\llvert Y_{ n + 1 }^N( \omega) \bigr\rrvert &\geq& \biggl
\llvert \mu \bigl( Y_n^N( \omega) \bigr) \cdot
\frac{ T }{ N } + \sigma \bigl( Y_n^N( \omega) \bigr)
\cdot \bigl( W_{ { ( n + 1 ) T }/{ N } }( \omega) - W_{ { n T }/{ N } }( \omega) \bigr) \biggr
\rrvert \\
&&{}- \bigl\llvert Y_n^N( \omega) \bigr\rrvert
\\
&=& \frac{ T }{ N } \bigl\llvert \mu \bigl( Y_n^N(
\omega) \bigr) \bigr\rrvert + \bigl\llvert \sigma \bigl( Y_n^N(
\omega) \bigr) \bigr\rrvert \cdot \bigl\llvert W_{ { ( n + 1 ) T }/{ N } }( \omega) -
W_{ { n T }/{ N } }( \omega) \bigr\rrvert \\
&&{}- \bigl\llvert Y_n^N(
\omega) \bigr\rrvert \\
&\geq& \frac{ T }{ N } \bigl( \bigl\llvert \mu \bigl(
Y_n^N( \omega) \bigr) \bigr\rrvert + \bigl\llvert \sigma
\bigl( Y_n^N( \omega) \bigr) \bigr\rrvert \bigr) - \bigl
\llvert Y_n^N( \omega) \bigr\rrvert ,
\end{eqnarray*}
and
the estimate
$ \llvert  \mu(x) \rrvert
+ \llvert  \sigma(x) \rrvert
\geq\frac{
\llvert  x \rrvert^{ \alpha} }{ c }
$
for all $ x \in\mathbb{R} $
with $ \llvert  x \rrvert  \geq c $,
inequality~\eqref{eqclaim2}
and
definition~\eqref{eqdefrN}
therefore show
\begin{eqnarray*}
 \bigl\llvert Y_{ n + 1 }^N( \omega) \bigr\rrvert &\geq&
\frac{ T }{ N c } \bigl\llvert Y_n^N( \omega) \bigr
\rrvert^{ \alpha} - \bigl\llvert Y_n^N( \omega)
\bigr\rrvert \geq \frac{ T }{ N c } \bigl\llvert Y_n^N(
\omega) \bigr\rrvert^{ \alpha} - \bigl\llvert Y_n^N(
\omega) \bigr\rrvert^{{ ( \alpha+ 1 ) }/{ 2 }}
\\
&=& \bigl\llvert Y_n^N( \omega) \bigr
\rrvert^{{( \alpha+ 1 )}/{ 2 }
} \biggl( \frac{ T }{ N c } \bigl\llvert
Y_n^N( \omega) \bigr\rrvert^{ { ( \alpha- 1 ) }/{ 2 } } - 1 \biggr)
\\
&\geq &\bigl\llvert Y_n^N( \omega) \bigr
\rrvert^{{ ( \alpha+ 1 )}/{ 2 }
} \biggl( \frac{ T }{ N c } ( r_N
)^{ { ( \alpha- 1 ) }/{ 2 } } - 1 \biggr) \geq \bigl\llvert Y_n^N(
\omega) \bigr\rrvert^{{ ( \alpha+ 1 ) }/{ 2 }
} .
\end{eqnarray*}
The induction hypothesis
hence yields
\[
\bigl\llvert Y_{ n + 1 }^N( \omega) \bigr\rrvert \geq \bigl
\llvert Y_n^N( \omega) \bigr\rrvert^{{ ( \alpha+ 1 ) }/{ 2 }
} \geq
\bigl( ( r_N )^{
(
(
{ (\alpha+ 1) }/{ 2 }
)^{ n }
)
} \bigr)^{{ ( \alpha+ 1 ) }/{ 2 }
} = (
r_N )^{
(
(
{ (\alpha+ 1 )}/{ 2 }
)^{ (n+1) }
)
} .
\]

Inequality~\eqref{eqclaim1} thus
holds for all
$ n \in\{ N_0, N_0 + 1, \ldots, N \} $,
$ \omega\in\Omega_N $
and all
$
N \in
\mathbb{N} \cap
\{ N_0, N_0 + 1, \ldots\}
$.
In particular, we obtain
%
%e23 #&#
\begin{equation}
\label{eqinductionend} \bigl\llvert Y_N^N( \omega) \bigr
\rrvert \geq ( r_N )^{
(
(
{ (\alpha+ 1) }/{ 2 }
)^{ N }
)
} \geq c^{
(
(
{( \alpha+ 1) }/{ 2 }
)^{ N }
)
}
\end{equation}
for all $ \omega\in\Omega_N $ and all
$
N \in
\mathbb{N} \cap
\{ N_0, N_0 + 1, \ldots\}
$.
Additionally,
Lemma~4.1 in~\cite{hjk11}
yields
%
%e24 #&#
\begin{eqnarray}
\label{eqincrementestimate} && \mathbb{E} \bigl[ 1_{
\{
\operatorname{sgn} (
\mu ( Y_n^N  )
\cdot
\sigma ( Y_n^N  )
)
\cdot
(
W_{ { ( n + 1 ) T }/{ N }
}
-
W_{ { n T }/{ N }
}
)
\geq
\frac{ T }{ N }
\}
} | \mathcal{ F
}_{ { n T }/{ N } } \bigr]
\nonumber\\
&&\qquad = \mathbb{P} \biggl[ ( W_{
{ ( n + 1 ) T
}/{ N }
} - W_{
{ n T
}/{ N }
} ) \geq
\frac{ T }{ N } \biggr] = \mathbb{P} \biggl[ W_{ { T }/{ N } } \geq
\frac{ T }{ N } \biggr] \\
&&\qquad= \mathbb{P} \biggl[ T^{ - {1}/{2} } W_T
\geq \sqrt{ \frac{ T }{ N } } \biggr] \geq \frac{
e^{ -{ T }/{ N } }
\sqrt{ T }
}{
8 \sqrt{ N }
}
\nonumber
\end{eqnarray}
$\P$-almost surely for
all $ n \in\{ 0, 1, \ldots, N - 1 \} $
and all $ N \in\N$.
Therefore, we obtain
\begin{eqnarray*}
\mathbb{P} [ \Omega_N ] & = &\mathbb{P} \bigl[ \bigl| Y_{ N_0 }^N
\bigr| \geq ( r_N )^{
(  (
{
(\alpha+ 1)
}/{
2
}
)^{ N_0 }  )
} \bigr] \cdot \biggl( \mathbb{P} \biggl[
T^{ - {1}/{2} } W_T \geq \sqrt{ \frac{ T }{ N } } \biggr]
\biggr)^{
(N - N_0)
}
\\
& \geq& \beta^{
(
-  ( N r_N  )^{
(
(
{ (\alpha+ 1) }/{ 2 }
)^{ N_0 }
\beta
)
}
)
} \cdot \biggl( \mathbb{P} \biggl[ T^{ - {1}/{2} }
W_T \geq \sqrt{ \frac{ T }{ N } } \biggr] \biggr)^{
N
}
\\
&\geq& \beta^{
(
-  ( N r_N  )^{
(
(
{ (\alpha+ 1) }/{ 2 }
)^{ N_0 }
\beta
)
}
)
} \cdot \biggl( \frac{
e^{ -{ T }/{ N } }
\sqrt{ T }
}{
8 \sqrt{ N }
} \biggr)^{
N
}
\\
& \geq& e^{ - T } \cdot \beta^{
(
-  ( N r_N  )^{
(
(
{( \alpha+ 1) }/{ 2 }
)^{ N_0 }
\beta
)
}
)
} \cdot \biggl(
\frac{
\sqrt{ T }
}{
8 \sqrt{ N }
} \biggr)^{
N
}
\end{eqnarray*}
for all
$
N \in
\mathbb{N} \cap
\{ N_0, N_0 + 1, \ldots\}
$.
This shows the existence
of a real number
$ \theta\in(1,\infty) $ such that
%
%e25 #&#
\begin{equation}
\label{eqOmegaprobab} \mathbb{P} [ \Omega_N ] \geq
\theta^{
(
- N^{ \theta}
)
}
\end{equation}
for all
$
N \in
\mathbb{N} \cap
\{ N_0, N_0 + 1, \ldots\}
$.
Combining~\eqref{eqinductionend} and~\eqref{eqOmegaprobab} finally
gives
\begin{eqnarray*}
\lim_{ N \rightarrow\infty} \mathbb{E} \bigl[ \bigl\llvert Y_N^N
\bigr\rrvert^p \bigr] &\geq& \lim_{ N \rightarrow\infty} \mathbb{E} \bigl[
1_{ \Omega_N } \bigl\llvert Y_N^N \bigr
\rrvert^p \bigr] %\nonumber
\geq \lim_{ N \rightarrow\infty}
\bigl( \mathbb{P} [ \Omega_N ] \cdot c^{
(
p \cdot
(
{ (\alpha+ 1 )}/{ 2 }
)^{
N
}
)
} \bigr)\\
& \geq&
\lim_{ N \rightarrow\infty} \bigl( \theta^{  ( - N^{ \theta}  )
} \cdot c^{
(
p \cdot
(
{( \alpha+ 1) }/{ 2 }
)^{
N
}
)
} \bigr) =
\infty
\end{eqnarray*}
for all $ p \in(0,\infty) $.
This, \eqref{eqinductionend} and~\eqref{eqOmegaprobab}
then complete the proof of Theorem~\ref{thmeEdivergence}.
\end{pf}

%
% %%%%%%%%%%%%%%%%%%%%%%%%%%%%%%%%%%%%%%%%%%%%%%%%%%%%%%%%%%%%%%%%%%
% Section
% %%%%%%%%%%%%%%%%%%%%%%%%%%%%%%%%%%%%%%%%%%%%%%%%%%%%%%%%%%%%%%%%%%
%
%s3 #&#
\section{Convergence of
the Monte Carlo Euler method}
\label{secdivergenceMCE}
The Monte Carlo Euler method has been
shown to converge
with probability one
%in the pathwise sense
for one-dimensional SDEs
with
superlinearly growing
and globally
one-sided Lipschitz continuous
drift coefficients
and with
globally Lipschitz continuous
diffusion
coefficients;
see~\cite{hj11}.
The Monte Carlo Euler method
is thus \textit{strongly consistent}
(see, e.g.,
Nikulin~\cite{n01},
Cram{\'e}r~\cite{c99} or
Appendix~A.1 in
Glasserman~\cite{g04})
for such SDEs.
After having reviewed
this convergence
result of the Monte Carlo Euler
method,
we complement
in this section
this convergence result
with the behavior of moments of
the Monte Carlo Euler
approximations
for such SDEs.
More precisely,
an immediate consequence
of Theorem~\ref{thmeEdivergence} is
Corollary~\ref{thmMCEdivergence}
below which shows
for such SDEs
that the Monte Carlo Euler
approximations
diverge in the strong $L^p$-sense
for every $p \in[1,\infty)$.
We emphasize that this strong
divergence result does not reflect
the behavior of the Monte Carlo
Euler method in a simulation
and it is presented for completeness
only.
Indeed, the events
on which the
Euler approximations diverge
(see
Theorem~\ref{thmeEdivergence})
are \textit{rare events}, and their
probabilities decay to zero
very rapidly; see, for example, Lemma~4.5 in~\cite{hj11}
for details.
This is the reason why the Monte
Carlo Euler method is strongly consistent
and thus
does converge
%with probability one
according to~\cite{hj11}; see also Theorem~\ref{thmMCEconvergence}
below
and Corollary~3.23
in
\cite{HutzenthalerJentzen2012PhiArxiv}.

Throughout this section assume
that the following setting is fulfilled.
Let $ T \in( 0, \infty) $,
let $ ( \Omega, \mathcal{F},
\mathbb{P} ) $ be a probability
space with a normal
filtration $ ( \mathcal{F}_t
)_{ t \in[0,T] } $,
let $ W^{ k } \dvtx[0,T] \times
\Omega\rightarrow
\mathbb{R} $,
$ k \in\N$, be a family of
independent one-dimensional
standard $ ( \mathcal{F}_t
)_{ t \in[0,T] } $-Brownian motions and let $ \xi^{ k } \dvtx\Omega
\rightarrow\mathbb{R} $, $ k \in\N$, be a family of independent
identically distributed
$ \mathcal{F}_0
/ \mathcal{B}( \mathbb{R} ) $-measurable
mappings with $ \mathbb{E} [ | \xi^1 |^p  ] < \infty$ for
all $
p \in[1,\infty) $. Moreover, let $ \mu, \sigma\dvtx\mathbb{R}
\rightarrow\mathbb{R} $ be two $ \mathcal{B}( \mathbb{R} ) /
\mathcal
{B}( \mathbb{R} ) $-measurable mappings such that
there exists a predictable stochastic process
$ X \dvtx[0,T] \times\Omega\rightarrow\mathbb{R} $
which satisfies
$
\int_0^T
\llvert  \mu( X_s ) \rrvert
+
\llvert  \sigma( X_s ) \rrvert^2 \,ds
<
\infty
$
$ \mathbb{P} $-almost surely and
%
%e26 #&#
\begin{equation}
\label{eqSDEmceuler} X_t = \xi^1 + \int
_0^t \mu( X_s ) \,ds + \int
_0^t \sigma( X_s )
\,dW_s^1
\end{equation}
$ \mathbb{P} $-almost surely for
all $ t \in[0,T] $.
The drift coefficient $ \mu$ is the infinitesimal mean of the process
$ X $ and the diffusion coefficient $ \sigma$ is the infinitesimal
standard deviation of the process $ X $.
We then define a family $ Y_n^{ N, k } \dvtx\Omega\rightarrow\mathbb
{R} $, $ n \in
\{ 0, 1, \ldots, N \} $, $ N, k \in\N$,
of Euler approximations by
$ Y_0^{ N, k } := \xi^k $ and
%
%e27 #&#
\begin{equation}
Y_{ n + 1 }^{ N, k } := Y_n^{ N, k } + \mu \bigl(
Y_n^{ N, k } \bigr) \cdot \frac{ T }{ N } + \sigma \bigl(
Y_n^{ N, k } \bigr) \cdot \bigl( W^k_{ { ( n + 1 ) T }/{ N } }
- W^k_{ { n T }/{ N } } \bigr)
\end{equation}
for all $ n \in\{ 0, 1, \ldots, N - 1 \} $
and all $ N, k \in\N$.
For clarity of exposition we recall
the following
convergence theorem from~\cite{hj11}.
Its proof can be found in~\cite{hj11}.
%
%th3.1 #&#
\begin{theorem}[(Strong consistency
and convergence with probability
one of the Monte Carlo
Euler method)]
\label{thmMCEconvergence}
Assume that the above setting
is fulfilled, let
$ \mu, \sigma, f \dvtx\mathbb{R}
\rightarrow\mathbb{R} $
be four times continuously
differentiable and
let $ c \in[0,\infty) $ be
a real number such that
$
(
x - y
) \cdot
(
\mu(x) -
\mu(y)
)
\leq
c
\llvert  x - y
\rrvert^2
$,
$
\llvert
\sigma(x) -
\sigma(y)
\rrvert
\leq
c
\llvert  x - y
\rrvert
$
and
$
\llvert  \mu^{(4)}(x)
\rrvert
+
\llvert  \sigma^{(4)}(x)
\rrvert
+
\llvert  f^{(4)}(x)
\rrvert
\leq c (1 + \llvert  x \rrvert^c)
$
for all $ x \in\mathbb{R} $.
Then there exist finite
$ \mathcal{F} /
\mathcal{B}( [0,\infty) )
$-measurable
mappings
$ C_{ \varepsilon}
\dvtx\Omega\rightarrow
[0,\infty) $,
$ \varepsilon\in(0,1) $,
such that
%
%e28 #&#
\begin{equation}
\Biggl\llvert \mathbb{E} \bigl[ f( X_T ) \bigr] - \frac{ 1 }{ N^2 }
\Biggl( \sum_{ k = 1 }^{ N^2 } f \bigl(
Y_N^{ N, k } \bigr) \Biggr) \Biggr\rrvert \leq
\frac{
C_{ \varepsilon}
}{
N^{
(
1 - \varepsilon
)
}
}
\end{equation}
for all $ N \in\N$
and all $ \varepsilon\in(0,1) $
$\mathbb{P}$-almost surely.
\end{theorem}
In contrast to pathwise convergence
of the Monte Carlo Euler
method for SDEs
with globally one-sided Lipschitz
continuous drift
and globally Lipschitz
continuous diffusion
coefficients
(see
Theorem~\ref{thmMCEconvergence}
above
for details),
strong convergence of the
Monte Carlo Euler method, in general,
fails to hold for such SDEs which
is established in the following
corollary of Theorem~\ref{thmeEdivergence},
that is, in
Corollary~\ref{thmMCEdivergence}.
As mentioned above we emphasize
that Corollary~\ref{thmMCEdivergence}
does not reflect the behavior
of the Monte Carlo Euler method
in a practical simulation
because the events
on which the Euler approximations
diverge
(see Theorem~\ref{thmeEdivergence})
are rare events, and their
probabilities decay to zero very
rapidly; see Lemma~4.5 in~\cite{hj11}
for details.

%
%co3.2 #&#
\begin{cor}[(Strong divergence of
the Monte Carlo Euler method)]
\label{thmMCEdivergence}
Assume that the above setting
is fulfilled
and let
$ \alpha, c \in( 1, \infty) $
be real numbers such that
$
\llvert  \mu( x ) \rrvert
+
\llvert  \sigma( x ) \rrvert
\geq
\frac{ \llvert  x \rrvert^{ \alpha}
}{ c }
$
for all $ x \in\mathbb{R} $
with $ | x | \geq c $.
Moreover, assume that
$
\mathbb{P} [
\sigma( \xi^1 ) \neq0
] > 0
$
or that
there exists a real
number $\beta\in(1,\infty)$ such that
$
\P [
| \xi^1 | \geq x
]
\geq
\beta^{ (-x^\beta )}
$
for all $x\in[1,\infty)$.
Moreover, let
$ f \dvtx\mathbb{R} \rightarrow
\mathbb{R} $ be
$ \mathcal{B}( \mathbb{R} )
/ \mathcal{B}( \mathbb{R} )
$-measurable with
$ f( x ) \geq\frac{ 1 }{ c }
| x |^{ { 1 }/{ c } } - c $
for all $ x \in\mathbb{R} $.
Then
%
%e29 #&#
\begin{equation}
\label{eqsec3toshow} \lim_{ N \rightarrow\infty} \mathbb{E} \Biggl[ \Biggl\llvert
\mathbb{E} \bigl[ f( X_T ) \bigr] - \frac{ 1 }{ N^2 } \Biggl( \sum
_{ k = 1 }^{ N^2 } f \bigl( Y_N^{ N, k }
\bigr) \Biggr) \Biggr\rrvert^p \Biggr] = \infty
\end{equation}
for all $ p \in[1, \infty) $.
\end{cor}
\begin{pf}%[Proof of Corollary~\ref{thmMCEdivergence}]
The triangle inequality,
Jensen's inequality and
the estimate
$
f( x ) \geq
\frac{ 1 }{ c }
| x |^{ { 1 }/{ c } } - c
$
for all $ x \in\mathbb{R} $ give
%
%e30 #&#
\begin{eqnarray}
\label{eqsec3use} %
&& \Biggl\llVert \mathbb{E} \bigl[ f(
X_T ) \bigr] - \frac{ 1 }{ N^2 } \Biggl( \sum
_{ k = 1 }^{ N^2 } f \bigl( Y_N^{ N, k }
\bigr) \Biggr) \Biggr\rrVert_{ L^p( \Omega; \mathbb{R} ) } \nonumber\\
&&\qquad\geq \frac{ 1 }{ N^2 } \Biggl\llVert
\sum_{ k = 1 }^{ N^2 } f \bigl(
Y_N^{ N, k } \bigr) \Biggr\rrVert_{ L^p( \Omega; \mathbb{R} ) } - \mathbb{E}
\bigl[ \bigl| f( X_T ) \bigr| \bigr]
\nonumber
\\[-8pt]
\\[-8pt]
\nonumber
&&\qquad \geq \frac{ 1 }{ N^2 } \mathbb{E} \Biggl[ \sum_{ k = 1 }^{ N^2 }
f \bigl( Y_N^{ N, k } \bigr) \Biggr] - \mathbb{E} \bigl[ \bigl| f(
X_T ) \bigr| \bigr] = \mathbb{E} \bigl[ f \bigl( Y_N^{ N, 1 }
\bigr) \bigr] - \mathbb{E} \bigl[\bigl | f( X_T ) \bigr| \bigr]\\
&&\qquad \geq
\frac{ 1 }{ c } \cdot \mathbb{E} \bigl[ \bigl\llvert Y_N^{ N, 1 }
\bigr\rrvert^{ { 1 }/{ c } } \bigr] - c - \mathbb{E} \bigl[ \bigl| f( X_T
) \bigr| \bigr] \nonumber
\end{eqnarray}
for all $ N \in\N$
and all $ p \in[1,\infty) $.
Combining \eqref{eqsec3use}
and
Theorem~\ref{thmeEdivergence}
then shows~\eqref{eqsec3toshow}
in the case
$ \mathbb{E} [ | f( X_T ) |  ]
< \infty$.
In the case
$ \mathbb{E} [ | f( X_T ) |  ]
= \infty$,
the estimate $ f(x) \geq-c $
for all $ x \in\mathbb{R} $ shows
$ \mathbb{E} [ f( X_T )  ]
= \infty$, and this implies~\eqref{eqsec3toshow}
in the case
$ \mathbb{E} [ | f( X_T ) |  ]
= \infty$.
The proof of Corollary~\ref{thmMCEdivergence}
is thus completed.
\end{pf}

%%%%%%%%%%%%%%%%%%%%%%%%%%%%%%%%%
%% Section
%%%%%%%%%%%%%%%%%%%%%%%%%%%%%%%%%

%s4 #&#
\section{Counterexamples to
convergence of the multilevel
Monte Carlo Euler method}
\label{secmultilevelpathwisedivergence}
Theorem~\ref{thmconj} below establishes divergence with probability one
of the multilevel Monte Carlo Euler method~\eqref{eqmlmc}
for the SDE~\eqref{eqSDEx5}.
This, in particular, proves that
the multilevel Monte Carlo
Euler method is
in contrast to the classical Monte
Carlo Euler method
not \textit{consistent}
(see, e.g.,
Nikulin~\cite{n01},
Cram{\'e}r~\cite{c99} or
Appendix~A.1 in
Glasserman~\cite{g04})
for the SDE~\eqref{eqSDEx5}.

Throughout this section assume that
the following setting is fulfilled.
Let $ T, \sigmab\in( 0, \infty) $,
let $ ( \Omega, \mathcal{F},
\mathbb{P} ) $ be a probability space
and let
$ \xi^{ l, k } \dvtx\Omega
\rightarrow\mathbb{R} $,
$ l \in\N_0 $,
$ k \in\N$, be a
family
of independent normally distributed
$ \mathcal{F} /
\mathcal{B}( \mathbb{R} )
$-measurable mappings with mean
zero
and standard deviation $\sigmab$.
Moreover, let $ X \dvtx[0,T] \times\Omega\rightarrow\mathbb{R} $
be the unique stochastic process with continuous sample paths which
fulfills the SDE
%
%e31 #&#
\begin{equation}
\label{eqSDEE} dX_t = -X_t^5 \,dt,\qquad
X_0 = \xi
\end{equation}
for $ t \in[0,T] $.
We then define a family of Euler approximations
$ Y_n^{ N, l, k } \dvtx\Omega\rightarrow\mathbb{R} $, $ n \in\{ 0,
1, \ldots, N \} $, $ N \in\N$,
$ l \in\N_0 $,
$ k \in\N$, by
$ Y_0^{ N, l, k } := \xi^{ l, k } $ and
%
%e32 #&#
\begin{equation}
\label{eqcounterexample} Y_{ n + 1 }^{ N, l, k } :=
Y_{ n }^{ N, l, k } - \bigl( Y_{ n }^{ N, l, k }
\bigr)^5 \cdot \frac{ T }{ N }
\end{equation}
for all $ n \in\{ 0, 1, \ldots, N - 1 \} $,
$ N \in\N$,
$ l \in\N_0 $
and all $ k \in\N$.
%
%%%%%%%%%%%%%%%%%%%%%%%%%%%%%%%%%
%
%th4.1 #&#
\begin{theorem}[{[Main result of
this article: Divergence
with probability one
of the multilevel Monte
Carlo Euler
method for the
SDE~\eqref{eqSDEE}]}]
\label{thmconj}
Assume that the above setting
is fulfilled.
Then
%
%e33 #&#
\begin{eqnarray}
\label{eqconj} \quad \mathop{\lim_{N \rightarrow\infty}}_{ \ld( N ) \in\N
} \Biggl| \frac{1}{N}
\sum_{ k = 1 }^{ N }\bigl | Y_1^{ 1, 0, k }
\bigr|^p + \sum_{ l = 1 }^{ \ld( N ) }
\frac{ 2^l }{ N } \sum_{ k = 1 }^{{ N }/{ 2^l } } \bigl( \bigl|
Y_{ 2^l }^{ 2^l, l, k } \bigr|^p - \bigl| Y_{ 2^{ (l-1) } }^{
2^{ (l-1) }, l, k
}
\bigr|^p \bigr) % -
% \mathbb{E}\Big[
% \big| X_T \big|^p
% \Big]
\Biggr| =\infty
\end{eqnarray}
$\P$-almost surely for all $p\in(0,\infty)$.
\end{theorem}
%

%f2 #&#
\begin{figure}

\includegraphics{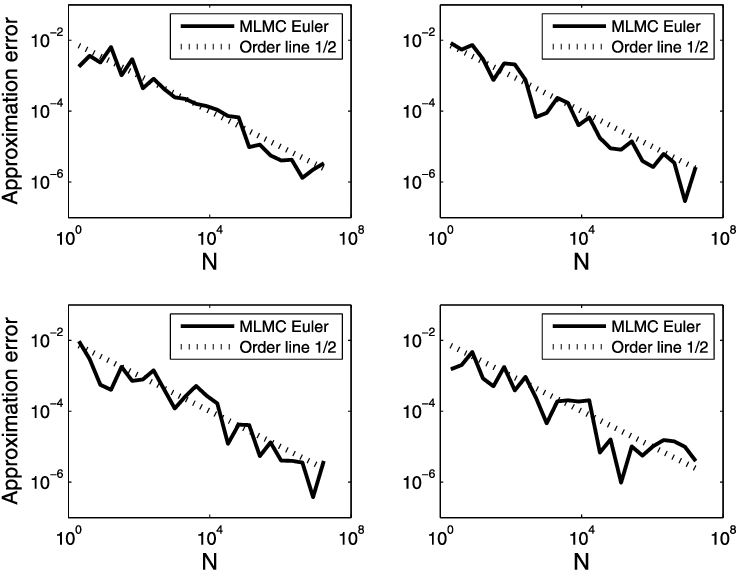}

\caption{Four sample paths of the approximation error of the multilevel
Monte Carlo Euler approximation
in~\protect\eqref{eqconj}
%for the SDE~\eqref{eqSDEx5}
where $T=1$, $\sigmab=0.1$, $p=2$ and $N\in\{2^1,2^2,\ldots,2^{22}\}$.}
\label{ffivesigma01}
\end{figure}

The proof of
Theorem~\ref{thmconj}
is postponed to
Section~\ref{secproofconj}
below.

%%%%%%%%%%%%%%%%%%%%%%%%%%%%%%%%%
%% Subsection
%%%%%%%%%%%%%%%%%%%%%%%%%%%%%%%%%

%s4.1 #&#
\subsection{Simulations}
\label{secSimulations}

\begin{figure}

\includegraphics{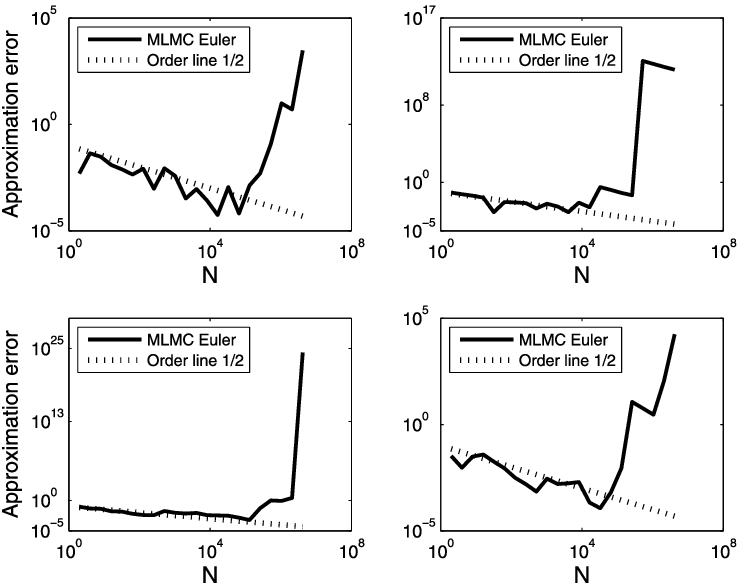}

\caption{Four sample paths of the approximation error of the multilevel
Monte Carlo Euler approximation
in~\protect\eqref{eqconj}
%for the SDE~\eqref{eqSDEx5}
where $T=1$, $\sigmab=\frac{1}{3}$, $p=2$ and $N\in\{2^1,2^2,\ldots
,2^{22}\}$.}
\label{ffivesigma033}
\end{figure}

We illustrate
Theorem~\ref{thmconj}
with numerical simulations.
To this end we observe that
the exact solution of the random
ordinary differential
equation~\eqref{eqSDE}
satisfies
%
%e34 #&#
\begin{equation}
\label{eqexactrep} X_t = \frac{ \xi}{
( 1 + 4 t \xi^4
)^{ { 1 }/{ 4 } }
}
\end{equation}
for all $t\in[0,1]$.
The real number
$ \E [ (X_1)^2  ]
$
can then
be computed approximatively
by numerical integration or
by the Monte Carlo method.
Figure~\ref{ffivesigma1}
depicts four random sample paths
of the approximation error of the
multilevel Monte Carlo Euler
approximations
in the case $T=1$ and $\sigmab=1$
in \eqref{eqSDE}
where
$
\E [ (X_1)^2  ]
\approx0.28801
$
(calculated with the \texttt{integrate}-function of~\texttt{R}).
The sample paths clearly diverge even for small
$
N \in
\{ 2^1, 2^2, 2^3, \ldots\}
$.
For some other SDEs, however, pathwise divergence does not emerge
for small
$
N \in\{ 2^1, 2^2, 2^3, \ldots\}
$.
For example, let us choose a standard deviation as small as
$
\sigmab= 0.1
$
in~\eqref{eqSDEE}
where
$
T=1
$.
Here the exact value
satisfies
$
\E [ (X_1)^2  ]
\approx0.009971
$
(calculated with the \texttt{integrate}-function of \texttt{R}).
Then sample paths of
the multilevel Monte Carlo Euler approximation seem to converge even
for reasonably large
$N\in\{2^1,2^2,2^3,\ldots\}$; see Figure~\ref{ffivesigma01} for four
sample paths.
%%%%%%%%%%%%%%%%%%%
%
%%%%%%%%%%%%%%%%%%%
So the sample paths of
the multilevel Monte Carlo Euler method
for some SDEs
first seem to converge, but diverge as
$N\in\{2^1,2^2,2^3,\ldots\}$
becomes sufficiently large.
To see this in a plot, we tried different values of $\sigmab$
and found sample paths in case of $\sigmab=\frac{1}{3}$ and $T=1$
which first seem to convergence
to the exact value
$
\E [ (X_1)^2  ]
\approx0.09248
$
(calculated with the \texttt{integrate}-function of \texttt{R})
but diverge for larger values of $N\in\{2^1,2^2,2^3,\ldots\}$; see
Figure~\ref{ffivesigma033} for four sample paths.
%%%%%%%%%%%%%%%%%%%%
%

%f3 #&#

%
%%%%%%%%%%%%%%%%%%%%

%%%%%%%%%%%%%%%%%%%%%%%%%%%%%%%%%
%% Subsection
%%%%%%%%%%%%%%%%%%%%%%%%%%%%%%%%%

%s4.2 #&#
\subsection{\texorpdfstring{Proof of Theorem~\protect\ref{thmconj}}{Proof of Theorem 4.1}}
\label{secproofconj}
First of all, we introduce more notation
in order to prove
Theorem~\ref{thmconj}.
Let $y_n^{N,x}\in\R$,
$n\in\{0,1,\ldots,N\}$,
$ N \in\mathbb{N} $,
$ x \in\mathbb{R} $,
be defined recursively
through $y_0^{N,x}:=x$ and
%
%e35 #&#
\begin{equation}
\label{eqconjY} y_{n+1}^{N,x} := y_n^{N,x}
- \bigl(y_n^{N,x} \bigr)^5 \cdot
\frac{T}{N} = y_n^{N,x} \biggl( 1 -
\bigl(y_n^{N,x} \bigr)^4 \cdot
\frac{T}{N} \biggr)
\end{equation}
for all $n\in\{0,1,\ldots,N-1\}$,
$N\in\N$ and all $x\in\R$
and let
$ p \in(0,\infty) $
be fixed for the rest of this
section.
This notation enables us
to rewrite the multilevel Monte
Carlo Euler
approximation
in \eqref{eqconj}
as
%
%e36 #&#
\begin{eqnarray}
\label{eqrewrite2} %
&&\frac{1}{N} \sum
_{ k = 1 }^{ N } \bigl| Y_1^{ 1, 0, k }
\bigr|^p + \sum_{ l = 1 }^{ \ld( N ) }
\frac{ 2^l }{ N } \sum_{ k = 1 }^{{ N }/{ 2^l } } \bigl( \bigl|
Y_{ 2^l }^{ 2^l, l, k } \bigr|^p -\bigl | Y_{ 2^{ (l-1) } }^{
2^{ (l-1) }, l, k
}
\bigr|^p \bigr)
\nonumber\\
&&\qquad = \sum_{ l = 0 }^{ \ld( N ) } \frac{ 2^l }{ N } \sum
_{ k = 1 }^{{ N }/{ 2^l } } \bigl\llvert
Y_{ 2^l }^{ 2^l, l, k } \bigr\rrvert^p - \sum
_{ l = 1 }^{ \ld( N ) } \frac{ 2^l }{ N } \sum
_{ k = 1 }^{{ N }/{ 2^l } } \bigl\llvert Y_{ 2^{ (l-1) } }^{
2^{ (l-1) }, l, k
}
\bigr\rrvert^p \\
&&\qquad= \sum_{ l = 0 }^{ \ld( N ) }
\frac{ 2^l }{ N } \sum_{ k = 1 }^{{ N }/{ 2^l } } \bigl
\llvert y_{ 2^l }^{ 2^l, \xi^{l, k} } \bigr\rrvert^p - \sum
_{ l = 1 }^{ \ld( N ) } \frac{ 2^l }{ N } \sum
_{ k = 1 }^{{ N }/{ 2^l } } \bigl\llvert y_{ 2^{ (l-1) } }^{
2^{ (l-1) }, \xi^{l, k}
}
\bigr\rrvert^p \nonumber
\end{eqnarray}
for all
$ N \in\{ 2^1, 2^2, 2^3, \ldots\} $.
Additionally, let
$
L_N \dvtx
\Omega
\to
\{ 1, 2, \ldots, \ld(N) \}
$ %,
%$
% N \in\{ 2^1, 2^2, 2^3, \ldots
% \}
%$,
be defined as
%
%e37 #&#
\begin{eqnarray}
\label{eqdefLN} L_N&:=&\max \biggl( \{ 1 \} \cup \biggl\{ l \in \bigl
\{ 1, 2, \ldots, \ld(N) \bigr\} \dvtx % \Big(
\nonumber
\\[-8pt]
\\[-8pt]
\nonumber
&&\hspace*{58pt}\exists k \in \biggl\{ 1, 2,
\ldots, \frac{N}{2^l} \biggr\} \dvtx \bigl| \xi^{l,k} \bigr| >
2^{{ l }/{ 4 } } T^{ - {1}/{4} } % \Big)
\biggr\} \biggr)
\end{eqnarray}
for every
$ N \in\{ 2^1, 2^2, 2^3, \ldots\}$.
Furthermore, define
$
\eta_N \dvtx\Omega\to
[0,\infty)
$
and
$
\theta_N \dvtx\break\Omega\to[0,\infty)
$
%for
%$
% N \in\{ 2^1, 2^2, 2^3, \ldots\}
%$
by
%
%e38 #&#
\begin{equation}
\label{eqdefetaN} \eta_N := \max \biggl\{\bigl | \xi^{L_N,k} \bigr| \in
\mathbb{R} \dvtx k \in \biggl\{ 1, 2, \ldots, \frac{ N }{ 2^{ L_N } } \biggr\} \biggr
\}
\end{equation}
and
%
%e39 #&#
\begin{equation}
\label{eqdefthetaN} \theta_N := \max \biggl\{\bigl | \xi^{(L_N - 1), k }
\bigr| \in\mathbb{R} \dvtx k \in \biggl\{ 1, 2, \ldots, \frac{N}{2^{(L_N-1)}} \biggr\}
\biggr\}
\end{equation}
for every
$N\in\{ 2^1, 2^2, 2^3, \ldots\}$.
Moreover, we define the mappings
$ \lceil\cdot\rceil,
\lfloor\cdot\rfloor\dvtx\mathbb{R}
\rightarrow\mathbb{Z} $ by
$\lceil x \rceil:=
\min\{z\in\mathbb{Z} \dvtx z\geq x\}$
and by
$\lfloor x\rfloor:=\max\{z\in\mathbb{Z} \dvtx z\leq x\}$
for all $x\in\mathbb{R} $.
Additionally, we fix a real number
$ \delta\in(0,\frac{1}{2}) $
for the rest of this section.
In the next step the following
events are used in our analysis
of the multilevel Monte Carlo
Euler method.
Let
$ A_N^{(1)} $,
$ A_N^{(2)} $,
$ A_N^{(3)} $,
$ A_N^{(4)} \in\mathcal{F} $,
$ N \in\{ 2^1, 2^2, 2^3, \ldots\} $,
be defined by
%
%e40 #&#
%e41 #&#
%e42 #&#
%e43 #&#
\begin{eqnarray}
A_N^{(1)} &:=& \bigl\{ L_N < \bigl\lfloor 2
\operatorname{ld} \bigl( \sigmab^2 T^{{1}/{2}
}
\operatorname{ln}(N) \bigr) \bigr\rfloor \bigr\} ,\label{eqA1}
\\
\label{eqA2}A_N^{(2)} & :=& \biggl\{ \exists l \in \bigl\{ 0, 1, 2,
\ldots, \ld(N) \bigr\} \dvtx
\nonumber
\\[-8pt]
\\[-8pt]
\nonumber
&&\hspace*{7pt}
\biggl( \exists k \in \biggl\{ 1, 2, \ldots,
\frac{N}{2^l} \biggr\} \dvtx \bigl| \xi^{l,k} \bigr| \geq 2^{
{ (l-1) }/{ 4 }
}
T^{ - {1}/{4} } N \biggr) \biggr\},
\\
\label{eqA3}A_N^{(3)} & := &\bigl\{ \exists l \in\N, \bigl\lfloor 2
\operatorname{ld} \bigl( \sigmab^2 T^{
{1}/{2}
}
\operatorname{ln}(N) \bigr) \bigr\rfloor \leq l\leq \operatorname{ld}(N) + 1 \dvtx
%
 % \left(
% \frac{ 2^{ l } }{ T }
% \right)^{{1}/{4} }
\nonumber
\\[-8pt]
\\[-8pt]
\nonumber
&&\hspace*{12pt}2^{ {l}/{4} } T^{ - {1}/{4} }\leq\eta_N < 2^{ {l}/{4} }
T^{ - {1}/{4} } % \left( \frac{ 2^{l} }{ T }
% \right)^{{1}/{4} }
%
\bigl( 1 + 5^{ (- \dl\cdot2^{ (l-1) } ) }
\bigr) \bigr\} ,
\\
A_N^{(4)}&:=& \bigl\{ \llvert \eta_N-
\theta_N\rrvert \leq4^{ (-2^{(L_N-1)} )}\eta_N \bigr\}
\label{eqA4}
\end{eqnarray}
for all $N\hspace*{-0.5pt}\in\hspace*{-0.5pt}\{2^1,2^2,2^3,\ldots\}$.
Additionally, define
$ N_0 \in\{ 2^1, 2^2, 2^3, \ldots\}
$
and
$
N_1 \dvtx\break\Omega\hspace*{-0.5pt}\rightarrow
\{ 2^1, 2^2, 2^3, \ldots\}
\cup\{ \infty\}
$
by
$
N_0 :=
2^{
\lceil\exp (
4 \bar{ \sigma}^{ - 2 }
T^{ - {1}/{2} }
)
+ \bar{\sigma}^{ 8 } T^2
\rceil
}
$
and by
%
%e44 #&#
\begin{eqnarray}
\qquad N_1(\omega) &:=& \min \bigl( \{ \infty\} \cup \bigl\{ n
\in \bigl\{ N_0, 2^1N_0,
2^2N_0, \ldots \bigr\} \dvtx % \\
% &
\forall m \in \bigl\{ n, 2^1n, 2^2n, \ldots \bigr\} \dvtx
\\
&&\hspace*{133pt}\qquad\omega\notin A_m^{ (1) } \cup A_m^{ (2) }
\cup A_m^{ (3) } \cup A_m^{ (4) } \bigr
\} \bigr) \nonumber%
\end{eqnarray}
for all $ \omega\in\Omega$.
Next we prove
a few lemmas that we use in our
proof of Theorem~\ref{thmconj}.
%
%%%%%%%%%%%%%%%%%%%%%%%%%%%%%%%%%
%
%le4.2 #&#
\begin{lemma}[(Dynamics for small initial values)] \label{lstability}
Assume that the above setting
is fulfilled.
Then we have
$
\llvert y_{n}^{N,x}\rrvert
\leq
|x|
\leq
(\frac{2N}{T} )^{{1}/{4}}
$
for all
$
n \in\{ 0, 1, \ldots, N \}
$,
$
|x| \geq
( \frac{ 2 N }{ T }
)^{ { 1 }/{ 4 }
}
$
and
all $ N \in\N$.
\end{lemma}
\begin{pf}%[Proof of Lemma~\ref{lstability}]
Fix $N\in\N$ and $|x|\leq (\frac{2N}{T} )^{{1}/{4}}$.
We prove
$\llvert y_{n}^{N,x}\rrvert \leq\llvert x\rrvert $
by induction
on $n\in\{0,1,\ldots,N\}$.
The base case $n=0$ is trivial.
For the induction step $n\to n+1$, note that
the induction hypothesis implies
%
%e45 #&#
\begin{equation}
\bigl\llvert y_{n+1}^{N,x} \bigr\rrvert = \bigl
\llvert y_n^{N,x} \bigr\rrvert \cdot \biggl\llvert 1 -
\frac{T}{N} \bigl(y_n^{N,x} \bigr)^4
\biggr\rrvert \leq \bigl\llvert y_n^{N,x} \bigr\rrvert \leq
\llvert x\rrvert %
\end{equation}
for all $n\in\{0,1,\ldots,N-1\}$.
This completes the proof of Lemma~\ref{lstability}.
\end{pf}
%
%%%%%%%%%%%%%%%%%%%%%%%%%%%%%%%%%
%
%le4.3 #&#
\begin{lemma}[(Dynamics for large initial values)] \label{linstability}
Assume that the above setting
is fulfilled.
Then we have
$
\llvert  y_n^{N,x}
\rrvert \geq\llvert x\rrvert
\geq ( \frac{2N}{T}  )^{{1}/{4}}
$
for all $n\in\{0,1,\ldots,N\}$,
$
|x| \geq
( \frac{ 2 N }{ T }
)^{ { 1 }/{ 4 }
}
$
and
all $N\in\N$.
In particular, we have
%
%e46 #&#
\begin{equation}
\label{eqabsoluteequation} \bigl\llvert y_{n+1}^{N,x} \bigr
\rrvert = \bigl\llvert y_n^{N,x} \bigr\rrvert \biggl(
\frac{T}{N} \bigl(y_n^{N,x} \bigr)^4-1
\biggr)
\end{equation}
for all $n\in\{0,1,\ldots,N-1\}$,
$|x|\geq (\frac{2N}{T} )^{{1}/{4}}$
and
all $N\in\N$.
\end{lemma}
\begin{pf}%[Proof of Lemma~\ref{linstability}]
Fix $N\in\N$ and $|x|\geq (\frac{2N}{T} )^{{1}/{4}}$.
We prove
$\llvert y_{n}^{N,x}\rrvert \geq\llvert x\rrvert $
by induction
on $n\in\{0,1,\ldots,N\}$.
The base case $n=0$ is trivial.
For the induction step $n\to n+1$, note that
the induction hypothesis implies
%
%e47 #&#
\begin{eqnarray}
\bigl\llvert y_{n+1}^{N,x} \bigr\rrvert& = &\bigl\llvert
y_n^{N,x} \bigr\rrvert \cdot \biggl\llvert
\frac{T}{N} \bigl(y_n^{N,x} \bigr)^4 - 1
\biggr\rrvert = \bigl\llvert y_n^{N,x} \bigr\rrvert \biggl(
\frac{T}{N} \bigl(y_n^{N,x} \bigr)^4-1
\biggr)
\nonumber
\\[-8pt]
\\[-8pt]
\nonumber
& \geq &\bigl\llvert y_n^{N,x} \bigr\rrvert \geq
\llvert x\rrvert
\end{eqnarray}
for all $n\in\{0,1,\ldots,N-1\}$.
This completes the induction.
Assertion~\eqref{eqabsoluteequation}
then im\-me\-dia\-te\-ly follows
by taking absolute values
in~\eqref{eqconjY}.
\end{pf}
%
%%%%%%%%%%%%%%%%%%%%%%%%%%%%%%%%%
%
%le4.4 #&#
\begin{lemma}[(Growth bound for large initial values)] \label{lupperesti}
Assume that the above setting
is fulfilled.
Then we have
%
%e48 #&#
\begin{equation}
\label{equpperesti} \biggl(\frac{T}{N} \biggr)^{{1}/{4}} \bigl
\llvert y_n^{N,x} \bigr\rrvert \leq \biggl( \biggl(
\frac{T}{N} \biggr)^{{1}/{4}}\llvert x\rrvert \biggr)^{ (5^n )}
\end{equation}
for all $n\in\{0,1,\ldots,N\}$,
$|x|\geq (\frac{2N}{T} )^{{1}/{4}}$
and
all $N\in\N$.
\end{lemma}
\begin{pf}%[Proof of Lemma~\ref{lupperesti}]
Fix $N\in\N$ and $|x|\geq (\frac{2N}{T} )^{{1}/{4}}$.
We prove~\eqref{equpperesti} by induction
on $n\in\{0,1,\ldots,N\}$.
The base case $n=0$ is trivial.
For the induction step $n\to n+1$, note that
Lemma~\ref{linstability}
and
the induction hypothesis imply
%
%e49 #&#
\begin{eqnarray}
\qquad \biggl(\frac{T}{N} \biggr)^{{1}/{4}} \bigl\llvert
y_{n+1}^{N,x} \bigr\rrvert &=& \biggl(\frac{T}{N}
\biggr)^{{1}/{4}} \bigl\llvert y_n^{N,x} \bigr\rrvert
\biggl(\frac{T}{N} \bigl(y_n^{N,x}
\bigr)^4-1 \biggr) \leq \biggl( \biggl(\frac{T}{N}
\biggr)^{{1}/{4}} \bigl\llvert y_{n}^{N,x} \bigr
\rrvert \biggr)^5
\nonumber
\\[-8pt]
\\[-8pt]
\qquad \nonumber
&\leq& \biggl( \biggl( \biggl(\frac{T}{N}
\biggr)^{{1}/{4}}\llvert x\rrvert \biggr)^{ (5^n )}
\biggr)^5 = \biggl( \biggl(\frac{T}{N} \biggr)^{{1}/{4}}
\llvert x\rrvert \biggr)^{ ( 5^{ (n+1) }  )} %
\end{eqnarray}
for all $n\in\{0,1,\ldots,N-1\}$.
This completes the proof of Lemma~\ref{lupperesti}.
\end{pf}
%
%%%%%%%%%%%%%%%%%%%%%%%%%%%%%%%%%
%
%le4.5 #&#
\begin{lemma}[(Monotonicity)] \label{lmonotonicityy}
Assume that the above setting
is fulfilled.
Then we have
%
%e50 #&#
\begin{equation}
\label{eqmonotonicityy} \bigl\llvert y_n^{N,x} \bigr
\rrvert \geq \bigl\llvert y_n^{N,y} \bigr\rrvert
\end{equation}
for all $n\in\{0,1,\ldots,N\}$,
all $x,y\in\R$ satisfying
$|x|\geq|y|$, $|x|\geq (\frac{2N}{T} )^{{1}/{4}}$
and
all $N\in\N$.
\end{lemma}
\begin{pf}%[Proof of Lemma~\ref{lmonotonicityy}]
Fix $N\in\N$ and
$ x, y \in\mathbb{R} $
with
$|x|\geq|y|$, $|x|\geq (\frac{2N}{T} )^{{1}/{4}}$.
We prove
\eqref{eqmonotonicityy} by induction
on $ n \in\{ 0, 1, \ldots, N\} $.
The base case $n=0$ is trivial.
For the induction step $n\to n+1$, note that
Lemma~\ref{linstability}
and
the induction hypothesis
imply
%
%e51 #&#
\begin{eqnarray}
\bigl\llvert y_{n+1}^{N,x} \bigr\rrvert &= &\bigl
\llvert y_{n}^{N,x} \bigr\rrvert \biggl( \frac{T}{N}
\bigl\llvert y_{n}^{N,x} \bigr\rrvert^4 - 1
\biggr) \geq \bigl\llvert y_{n}^{N,y} \bigr\rrvert \biggl(
\frac{T}{N} \bigl\llvert y_{n}^{N,x} \bigr
\rrvert^4 - 1 \biggr) \
\nonumber
\\[-8pt]
\\[-8pt]
\nonumber
&\geq& \bigl\llvert y_{n}^{N,y}
\bigr\rrvert \biggl\llvert \frac{T}{N} \bigl\llvert y_{n}^{N,y}
\bigr\rrvert^4 - 1 \biggr\rrvert = \bigl\llvert y_{n+1}^{N,y}
\bigr\rrvert %
\end{eqnarray}
for all $n\in\{0,1,\ldots,N-1\}$.
This completes the proof of Lemma~\ref{lmonotonicityy}.
\end{pf}
%
%%%%%%%%%%%%%%%%%%%%%%%%%%%%%%%%%
%
%le4.6 #&#
\begin{lemma}[(Dynamics of multiples of the initial value)] \label
{linitialmultiple}
Assume that the above setting
is fulfilled.
Then we have
%
%e52 #&#
\begin{equation}
\label{eqinitialmultiple} \bigl\llvert y_n^{N,Mx} \bigr
\rrvert \geq M^{ (5^n )} \bigl\llvert y_n^{N,x} \bigr
\rrvert
\end{equation}
for all $n\in\{0,1,\ldots,N\}$,
$|x|\geq (\frac{2N}{T} )^{{1}/{4}}$,
$ M \in[1,\infty) $
and
all $N\in\N$.
\end{lemma}
\begin{pf}%[Proof of Lemma~\ref{linitialmultiple}]
Fix $N\in\N$.
We prove~\eqref{eqinitialmultiple}
by induction
on $n\in\{0,1,\ldots,N\}$.
The base case $ n = 0 $ is trivial.
For the induction step $n\to n+1$, note that
Lemma~\ref{linstability}
and
the induction hypothesis
imply
%
%e53 #&#
\begin{eqnarray}
\bigl\llvert y_{n+1}^{N,Mx} \bigr\rrvert &= &\bigl
\llvert y_{ n }^{N,Mx} \bigr\rrvert \biggl( \frac{T}{N}
\bigl\llvert y_{ n }^{N,Mx} \bigr\rrvert^4 - 1
\biggr)
\nonumber
\\
&\geq& M^{
(
5^n
)
} \bigl\llvert y_{ n }^{N,x}
\bigr\rrvert \biggl( \frac{T}{N} \bigl( M^{
(
5^n
)
} \bigl\llvert
y_{ n }^{N,x} \bigr\rrvert \bigr)^{ 4 } - 1 \biggr)
\\
& \geq& M^{
(
5^n
)
} \bigl\llvert y_{ n }^{N,x} \bigr
\rrvert \biggl( \frac{T}{N} \bigl( M^{
(
5^n
)
} \bigl\llvert
y_{ n }^{N,x} \bigr\rrvert \bigr)^{ 4 } -
M^{
(
4 \cdot
5^n
)
} \biggr) = M^{
( 5^{(n+1)}  )
} \bigl\llvert y_{n+1}^{N,x}
\bigr\rrvert \nonumber
\end{eqnarray}
for all $n\in\{0,1,\ldots,N-1\}$,
$
|x| \geq
(
\frac{ 2 N }{ T }
)^{
{ 1 }/{ 4 }
}
$
and
all $M\in[1,\infty)$.
This completes the proof of Lemma~\ref{linitialmultiple}.
\end{pf}

%co4.7 #&#
\begin{cor}
\label{cordoublegrowth}
Assume that the above setting
is fulfilled. Then we have
$
\llvert  y^{ N, x }_n \rrvert
\geq
M^{
( 5^n  )
}
(
\frac{ 2 N }{ T }
)^{{1}/{4} }
$
for all
$ n \in\{ 0, 1, \ldots, N \} $,
$
| x | \geq
M
(
\frac{ 2 N }{ T }
)^{{1}/{4} }
$,
$ M \in[ 1, \infty) $
and all
$ N \in\mathbb{N} $.
\end{cor}
\begin{pf}%[Proof of Corollary~\ref{cordoublegrowth}]
Lemmas~\ref{lmonotonicityy},~\ref{linitialmultiple}
and $|y^{N,({2N}/{T})^{{1}/{4}}}_n|=(\frac{2N}{T})^{{1}/{4}}$
imply
%
%e54 #&#
\begin{eqnarray}\quad
\bigl\llvert y^{ N, x }_n \bigr\rrvert \geq \bigl|
y^{ N,
M
(
{ 2 N }/{ T }
)^{{1}/{4} }
}_n \bigr| \geq M^{
( 5^n  )
}\bigl |y^{N,({2N}/{T})^{{1}/{4}}}_n
\bigr| = M^{
( 5^n  )
} \biggl( \frac{ 2 N }{ T } \biggr)^{{1}/{4} }
\end{eqnarray}
for all
$ n \in\{ 0, 1, \ldots, N \} $,
$
| x | \geq
M
(
\frac{ 2 N }{ T }
)^{{1}/{4} }
$,
$ M \in[ 1, \infty) $
and all
$ N \in\mathbb{N} $. This completes the proof of Corollary~\ref
{cordoublegrowth}.
\end{pf}

%
%%%%%%%%%%%%%%%%%%%%%%%%%%%%%%%%%
%
%le4.8 #&#
\begin{lemma} \label{lloweresti}
Assume that the above setting
is fulfilled.
Then we have
%
%e55 #&#
\begin{equation}
\label{eqloweresti} \bigl\llvert y_N^{N,x} \bigr\rrvert
\geq \biggl( \frac{ 2 N }{ T } \biggr)^{{1}/{4} } \sqrt{e}^{
( 5^{ (1-r) N }  )
}
\end{equation}
for all
$
|x| \geq
( \frac{ 2 N }{ T } )^{{1}/{4}}
(1+5^{(-r N)} )
$,
$
N \in\N
$
and all
$ r \in(0,\infty) $.
\end{lemma}
\begin{pf}%[Proof of Lemma~\ref{lloweresti}]
We apply the inequality $1+z\geq\exp(\frac{z}{2})$
for all $z\in[0,2]$.
Noting that $5^{(-r N)}\leq1\leq2$ for all $N\in\N$
and all $r\in(0,\infty)$,
we infer from Corollary~\ref{cordoublegrowth}
%
%e56 #&#
\begin{eqnarray}\quad
\bigl\llvert y_{N}^{N,x} \bigr\rrvert &\geq&
\biggl(\frac{2N}{T} \biggr)^{{1}/{4}} \bigl(1+5^{(-r N)}
\bigr)^{(5^N)} \geq \biggl(\frac{2N}{T} \biggr)^{{1}/{4}}
\biggl[\exp \biggl(\frac{1}{2}5^{(-r N)} \biggr)
\biggr]^{(5^N)}
\nonumber
\\[-8pt]
\\[-8pt]
\quad \nonumber
&=& \biggl(\frac{2N}{T} \biggr)^{{1}/{4}}
\sqrt{e}^{ (5^{(1-r) N} )} %
\end{eqnarray}
for all
$
|x| \geq
( \frac{ 2 N }{ T } )^{{1}/{4}}
(1+5^{(-r N)} )
$,
$
N \in\N
$
and all
$ r \in(0,\infty) $. This completes the proof of Lemma~\ref{lloweresti}.
\end{pf}

%%%%%%%%%%%%%%%%%%%%%%%%%%%%%%%%%%%%%%%%%%%%%%
%
%le4.9 #&#
\begin{lemma}[(Almost sure
finiteness of $ N_1 $)]
\label{lem_n1}
Assume that the above setting
is fulfilled. Then
$
\mathbb{P} [
N_1 < \infty
]
= 1
$.
\end{lemma}

The proof
of Lemma~\ref{lem_n1}
is postponed to the \hyperref[secappendix]{Appendix}.
We now present the proof
of Theorem~\ref{thmconj}.
It makes of use of Lemma~\ref{lem_n1}.

%%%%%%%%%%%%%%%%%%%%%
%
\begin{pf*}{Proof of
Theorem~\ref{thmconj}}
Fix $p\in(0,\infty)$ throughout this proof.
Our proof of
Theorem~\ref{thmconj} is then
divided into
four parts.
In the first part we analyze the
behavior of the
multilevel Monte Carlo Euler
approximations on the events
$
\{
\eta_N >
2^{ (L_N+1) / 4 }
T^{ - 1/4 }
\}
\cap
\{
N_1 \leq N
\}
=
\{
\omega\in\Omega\dvtx
\eta_N(\omega) >
2^{ (L_N(\omega)+1) / 4 }
T^{ - 1/4 } ,
N_1(\omega) \leq N
\}
$
for
$ N \in\{2^1,2^2,2^3,\ldots\} $; see inequality~\eqref{eqresultcase1}.
In the second part
of this proof we concentrate on\vadjust{\goodbreak}
the events
$
\{
\theta_N \geq
\eta_N
\}
\cap
\{
\eta_N \leq
2^{ (L_N+1) / 4 }
T^{ - 1/4 }
\}
\cap
\{
N_1 \leq N
\}
$
for
$ N \in\{2^1,2^2,2^3,\ldots\} $; see inequality~\eqref{eqresultcase2}.
In the third part
of this proof we investigate
the events
$\{2^{ L_N / 4 }
T^{ - 1/4 }<
\theta_N
\}\cap
\{\theta_N <
\eta_N \}
\cap\{ \eta_N \leq
2^{ (L_N+1) / 4 }
T^{ - 1/4 }
\}\cap \{N_1 \leq N\}$
for
$ N \in\{2^1,2^2,2^3,\ldots\} $
[see inequality~\eqref{eqresultcase3}]
and in the fourth part we analyze
the behavior of
the multilevel Monte Carlo
Euler approximations on the events
$
\{
2^{ L_N / 4 }
T^{ - 1/4 }
\geq
\theta_N
\}
\cap
\{
\theta_N <
\eta_N
\}
\cap
\{
\eta_N \leq
2^{ (L_N+1) / 4 }
T^{ - 1/4 }
\}
\cap
\{
N_1 \leq N
\}
$
for
$ N \in\{2^1,2^2,2^3,\ldots\} $
[see inequality~\eqref{eqresultcase4}].
Combining all four parts
[inequalities~\eqref{eqresultcase1},
\eqref{eqresultcase2},
\eqref{eqresultcase3} and
\eqref{eqresultcase4}]
and
Lemma~\ref{lem_n1} will then complete
the proof of Theorem~\ref{thmconj}
as we
will show below. In these four parts we will frequently use
%
%e57 #&#
\begin{equation}
\{N_1\leq N\} \subseteq \bigl(A_N^{(1)}
\bigr)^{c} \cap \bigl(A_N^{(2)}
\bigr)^{c} \cap \bigl(A_N^{(3)}
\bigr)^{c} \cap \bigl(A_N^{(4)}
\bigr)^{c}
\end{equation}
for all $ N \in\{ 2^1, 2^2, 2^3, \ldots\} $.

We begin with the first part
and consider the events
$
\{
\eta_N >\break
2^{ (L_N+1) / 4 }
T^{ - 1/4 }
\}
\cap
\{
N_1 \leq N
\}
$
for $ N \in\{ 2^1, 2^2, 2^3, \ldots\} $.
Note that Lemma~\ref{lmonotonicityy}, the inequalities
$
\eta_N
\geq
2^{ (L_N+1) / 4 }
T^{ - 1/4 }
(
1 +
5^{ ( -\dl\cdot2^{L_N} ) }
)
$
on
$
\{
\eta_N >\break
2^{ (L_N+1) / 4 }
T^{ - 1/4 }
\}
\cap
\{
N_1 \leq N
\}
$
[see \eqref{eqA3}]
and
$
| \xi^{ l,k} |
< 2^{ (l-1)/4 } T^{ -1/4 }
N
$ on
$  \{ N_1 \leq N  \} $ for all
$ k \in\{ 1, 2, \ldots, \frac{N}{2^l} \} $,
$ l \in\{ 1, 2, \ldots, \operatorname{ld}( N ) \} $
[see \eqref{eqA2}]
and the definition~\eqref{eqdefLN} of $ L_N $ imply
%
%e58 #&#
\begin{eqnarray}
&& \sum_{ l = 0 }^{ \ld( N ) }
\frac{ 2^l }{ N } \sum_{ k = 1 }^{{ N }/{ 2^l } } \bigl
\llvert y_{ 2^l }^{ 2^l, \xi^{l,k} } \bigr\rrvert^{ p } - \sum
_{ l = 1 }^{ \ld( N ) } \frac{ 2^{l} }{ N } \sum
_{ k = 1 }^{{ N }/{ 2^{l} } } \bigl\llvert y_{ 2^{ (l-1) } }^{
2^{ (l-1) }, \xi^{l,k}
}
\bigr\rrvert^{ p }
\nonumber\\
&&\qquad \geq \frac{ 2^{L_N} }{ N } \bigl\llvert y_{ 2^{L_N} }^{ 2^{L_N}, \eta_N } \bigr
\rrvert^{ p } - \sum_{ l = 1 }^{ L_N }
\frac{ 2^{l} }{ N } \sum_{ k = 1 }^{{ N }/{ 2^{l} } } \bigl
\llvert y_{ 2^{ (l-1) } }^{
2^{ (l-1) }, \xi^{l,k}
} \bigr\rrvert^{ p } - \sum
_{ l = L_N+1 }^{ \ld( N ) } \frac{ 2^{l} }{ N } \sum
_{ k = 1 }^{{ N }/{ 2^{l} } } \bigl\llvert y_{ 2^{ (l-1) } }^{
2^{ (l-1) }, \xi^{l,k}
}
\bigr\rrvert^{ p }
\nonumber
\\[-8pt]
\\[-8pt]
\nonumber
&&\qquad \geq \frac{ 2^{L_N} }{ N } \bigl\llvert y_{ 2^{L_N} }^{ 2^{L_N},
(
{
2^{(L_N + 1)}
}/{ T }
)^{
{1}/{4}
}
(
1 + 5^{ (-\dl2^{L_N}) }
)
} \bigr
\rrvert^{ p } \\
&&\qquad\quad{}- \sum_{ l = 1 }^{ L_N }
\bigl\llvert y_{ 2^{ (l-1) } }^{
2^{ (l-1) } ,
(
{ 2^{(l-1)} }/{ T }
)^{{1}/{4} }
N
} \bigr\rrvert^{ p } - \sum
_{ l = L_N+1 }^{ \ld( N ) } \bigl\llvert
y_{ 2^{ (l-1) } }^{
2^{ (l-1) },
(
{ 2^l }/{ T }
)^{
{1}/{4}
}
} \bigr\rrvert^{ p } \nonumber
\end{eqnarray}
on
$
\{
\eta_N >
2^{ (L_N + 1)/4} T^{ -1/4 }
\}
\cap
\{
N_1 \leq N
\}
$
and Lemmas~\ref{lloweresti},~\ref{lupperesti} and~\ref
{lstability} hence yield
\begin{eqnarray*}
&& \sum_{ l = 0 }^{ \ld( N ) }
\frac{ 2^l }{ N } \sum_{ k = 1 }^{{ N }/{ 2^l } } \bigl
\llvert y_{ 2^l }^{ 2^l, \xi^{l,k} } \bigr\rrvert^{ p } - \sum
_{ l = 1 }^{ \ld( N ) } \frac{ 2^{l} }{ N } \sum
_{ k = 1 }^{{ N }/{ 2^{l} } } \bigl\llvert y_{ 2^{ (l-1) } }^{
2^{ (l-1) }, \xi^{l,k}
}
\bigr\rrvert^{ p }
\\
&&\qquad \geq \frac{1}{N} \biggl\llvert \biggl( \frac{ 2\cdot2^{L_N} }{ T }
\biggr)^{
{1}/{4}
} \sqrt{e}^{
(
5^{ (1-\dl) 2^{L_N} }
)
} \biggr\rrvert^{ p }
\\
&&\qquad\quad{}- \sum_{ l = 1 }^{ L_N } \biggl(
\frac{ 2^{(l-1)} }{ T } \biggr)^{ {p}/{4} } N^{
(
p \cdot5^{
(
2^{ (l-1) }
)
}
)
} - \sum
_{ l = L_N+1 }^{ \ld( N ) } \biggl( \frac{ 2^l }{ T }
\biggr)^{
{p}/{4}
}
\\
&&\qquad \geq N^{ - 1 } T^{ - {p}/{4} } \cdot \exp \biggl( \frac{p}{2}
\cdot 5^{ (1-\dl) 2^{L_N} } \biggr)\\
&&\qquad\quad{} - L_N \biggl( \frac{ 2^{(L_N-1)} }{ T }
\biggr)^{ {p}/{4} } N^{
(
p \cdot5^{  ( 2^{(L_N-1)}  ) }
)
} - \ld( N ) \biggl( \frac{2^{\ld(N)}}{T}
\biggr)^{ {p}/{4} } %
\end{eqnarray*}
on
$
\{
\eta_N >
2^{ (L_N + 1)/4} T^{ -1/4 }
\}
\cap
\{
N_1 \leq N
\}
$
for all
$ N \in\{2^1,2^2,2^3,\ldots\} $.
Therefore, we obtain
%
%e59 #&#
\begin{eqnarray}
&& \sum_{ l = 0 }^{ \ld( N ) }
\frac{ 2^l }{ N } \sum_{ k = 1 }^{{ N }/{ 2^l } } \bigl
\llvert y_{ 2^l }^{ 2^l, \xi^{l,k} } \bigr\rrvert^{ p } - \sum
_{ l = 1 }^{ \ld( N ) } \frac{ 2^{l} }{ N } \sum
_{ k = 1 }^{{ N }/{ 2^{l} } } \bigl\llvert y_{ 2^{ (l-1) } }^{
2^{ (l-1) }, \xi^{l,k}
}
\bigr\rrvert^{ p }
\nonumber\\
&&\qquad\geq N^{ - 1 } T^{ - {p}/{4} } \cdot \exp \biggl( \frac{p}{2}
\cdot 5^{ (1-\dl) 2^{L_N} } \biggr)
\nonumber
\\
&&\qquad\quad{}- \ld( N ) N^{ {p}/{4} } T^{ - {p}/{4} }
\cdot N^{
(
p \cdot5^{  ( 2^{(L_N-1)}  ) }
)
} - \ld( N ) N^{ {p}/{4} } T^{ - {p}/{4} }
\\
&&\qquad \geq T^{ - {p}/{4} } \cdot \exp \biggl( \frac{p}{2} \cdot
5^{ (1-\dl) 2^{L_N} } - \operatorname{ln}( N ) \biggr)\nonumber \\
&&\qquad\quad{}- T^{ - {p}/{4} } \cdot
N^{
( 1 + {p}/{4} +
p \cdot5^{  ( 2^{(L_N-1)}  ) }
)
} \nonumber
\end{eqnarray}
on
$
\{
\eta_N >
2^{ (L_N + 1)/4} T^{ -1/4 }
\}
\cap
\{
N_1 \leq N
\}
$
and the estimate $2^{L_N}\geq\sigmab^2\sqrt{T}\ln(N)$ on $ \{
N_1 \leq N
\} $ [see~\eqref{eqA1}] hence shows
%
%e60 #&#
\begin{eqnarray}
\label{eqresultcase1}
 &&\sum_{ l = 0 }^{ \ld( N ) }
\frac{ 2^l }{ N } \sum_{ k = 1 }^{{ N }/{ 2^l } } \bigl
\llvert y_{ 2^l }^{ 2^l, \xi^{l,k} } \bigr\rrvert^{ p } - \sum
_{ l = 1 }^{ \ld( N ) } \frac{ 2^{l} }{ N } \sum
_{ k = 1 }^{{ N }/{ 2^{l} } } \bigl\llvert y_{ 2^{ (l-1) } }^{
2^{ (l-1) }, \xi^{l,k}
}
\bigr\rrvert^{ p }
\nonumber\\
&&\qquad \geq \inf_{
x \in[
\sigmab^2\sqrt{T} \operatorname{ln}(N) , %\bar{\sigma}^2 \sqrt{T} ,
\infty
)
} \biggl[ \exp \biggl( \frac{p}{2} \cdot
5^{  (  (1-\dl ) x  ) } - \operatorname{ln}( N ) \biggr)
\nonumber
\\[-8pt]
\\[-8pt]
\nonumber
&&\hspace*{84pt}\qquad{} - \exp \biggl(
\operatorname{ln}(N) \biggl( 1 + \frac{p}{4} + p \cdot 5^{ {x}/{2} }
\biggr) \biggr) \biggr] \cdot T^{ - {p}/{4} } \\
&&\qquad\geq r(N) \cdot T^{ - {p}/{4} }\nonumber
\end{eqnarray}
on
$
\{
\eta_N >
2^{ (L_N + 1)/4} T^{ -1/4 }
\}
\cap
\{
N_1 \leq N
\}
$
for all $N\in\{2^1,2^2,2^3,\ldots\}$
where
$ r \dvtx\mathbb{N} \rightarrow\mathbb{R} $
is a function defined by
\begin{eqnarray*}
r( N ) &:=& \inf_{
x \in [
\bar{\sigma}^2 \sqrt{T} \operatorname{ln}(N) ,
\infty
)
} \biggl[ \exp \biggl( \frac{p}{2} \cdot
5^{  (  (1-\dl ) x  ) } - \operatorname{ln}(2N) \biggr) \\
&&\hspace*{52pt}\qquad{}- \exp \biggl(
\operatorname{ln}(N) \biggl( 1 + \frac{p}{4} + p \cdot 5^{ {x}/{2} }
\biggr) \biggr) \biggr]
\end{eqnarray*}
for all $ N \in\mathbb{N} $.

In the next step we analyze the
behavior of the multilevel Monte Carlo
Euler approximations on the events\vadjust{\goodbreak}
$
\{
\theta_N \geq
\eta_N
\}
\cap
\{
\eta_N \leq
2^{ (L_N + 1)/4} T^{ -1/4 }
\}
\cap
\{
N_1 \leq N
\}
$
for
$ N \in\{2^1,2^2,2^3,\ldots\} $.
To this end note that Lemma~\ref{lmonotonicityy}, the inequalities
$
\theta_N \geq
(
1 + 4^{  ( -2^{ (L_N - 1) }  ) }
) \eta_N
$
on $  \{ \theta_N \geq\eta_N  \} \cap \{ N_1 \leq N
\} $ [see \eqref{eqA4}] and
$
| \xi^{ l,k} |
< 2^{ (l-1)/4 } T^{ -1/4 }
N
$ on
$  \{ N_1 \leq N  \} $ for all $ k \in\{ 1, 2, \ldots,
\frac
{N}{2^l} \} $, $ l \in\{ 1, 2, \ldots,\break \operatorname{ld}( N ) \} $
[see \eqref
{eqA2}] and the definition~\eqref{eqdefLN} of $ L_N $ imply
\begin{eqnarray*}
&& \sum_{ l = 0 }^{ \ld( N ) }
\frac{ 2^l }{ N } \sum_{ k = 1 }^{{ N }/{ 2^l } } \bigl
\llvert y_{ 2^l }^{ 2^l, \xi^{l,k} } \bigr\rrvert^{ p } - \sum
_{ l = 1 }^{ \ld( N ) } \frac{ 2^{l} }{ N } \sum
_{ k = 1 }^{{ N }/{ 2^{l} } } \bigl\llvert y_{ 2^{ (l-1) } }^{
2^{ (l-1) }, \xi^{l,k}
}
\bigr\rrvert^{ p }
\\
& &\qquad\geq \frac{1}{N} \bigl\llvert y_{ 2^{ (L_N - 1) }
}^{ 2^{ (L_N - 1) }, \theta_N
} \bigr
\rrvert^{ p } - \bigl\llvert y_{ 2^{ (L_N-1) }
}^{ 2^{ (L_N-1) }, \eta_N
} \bigr
\rrvert^{ p } \\
&&\qquad\quad{}- \sum_{ l = 1 }^{ L_N - 1 }
\frac{ 2^{l} }{ N } \sum_{ k = 1 }^{{ N }/{ 2^{l} } } \bigl
\llvert y_{ 2^{ (l-1) } }^{
2^{ (l-1) }, \xi^{l,k}
} \bigr\rrvert^{ p } - \sum
_{ l = L_N + 1 }^{ \ld(N) } \frac{ 2^{l} }{ N } \sum
_{ k = 1 }^{{ N }/{ 2^{l} } } \bigl\llvert
y_{ 2^{ (l-1) } }^{
2^{ (l-1) }, \xi^{l,k}
} \bigr\rrvert^{ p }
\\
&&\qquad \geq \frac{1}{N} \bigl\llvert y_{ 2^{ (L_N - 1) }
}^{ 2^{ (L_N - 1) },
(
1 +
4^{  ( - 2^{ (L_N-1) }  ) }
)
\eta_N
} \bigr
\rrvert^{ p } - \bigl\llvert y_{ 2^{ (L_N-1) }
}^{ 2^{ (L_N-1) }, \eta_N
} \bigr
\rrvert^{ p } \\
&&\qquad\quad{}- \sum_{ l = 1 }^{ L_N - 1 }
\bigl\llvert y_{ 2^{ (l-1) } }^{
2^{ (l-1) },
(
{ 2^{(l-1)} }/{ T }
)^{{1}/{4} }
N
} \bigr\rrvert^{ p } - \sum
_{ l = L_N+1 }^{ \ld( N ) } \bigl\llvert
y_{ 2^{ (l-1) } }^{
2^{ (l-1) },
(
{2^l}/{T}
)^{{1}/{4} }
} \bigr\rrvert^{ p } %
\end{eqnarray*}
on
$
\{
\theta_N \geq
\eta_N
\}
\cap
\{
\eta_N \leq
2^{ (L_N + 1)/4} T^{ -1/4 }
\}
\cap
\{
N_1 \leq N
\}
$
for all
$ N \in\{2^1,2^2,2^3,\ldots\} $.
%%%%
Lemmas~\ref{lmonotonicityy},
\ref{lupperesti} and~\ref{lstability} therefore show
%
%e61 #&#
\begin{eqnarray}
&&\sum_{ l = 0 }^{ \ld( N ) }
\frac{ 2^l }{ N } \sum_{ k = 1 }^{{ N }/{ 2^l } } \bigl
\llvert y_{ 2^l }^{ 2^l, \xi^{l,k} } \bigr\rrvert^{ p } - \sum
_{ l = 1 }^{ \ld( N ) } \frac{ 2^{l} }{ N } \sum
_{ k = 1 }^{{ N }/{ 2^{l} } } \bigl\llvert y_{ 2^{ (l-1) } }^{
2^{ (l-1) }, \xi^{l,k}
}
\bigr\rrvert^{ p }
\nonumber\\
&&\qquad \geq \frac{1}{2 N} \bigl\llvert y_{ 2^{ (L_N - 1) }
}^{ 2^{ (L_N - 1) },
(
1 +
4^{  ( - 2^{ (L_N-1) }  ) }
)
\eta_N
} \bigr
\rrvert^{ p } - \bigl\llvert y_{ 2^{ (L_N-1) }
}^{ 2^{ (L_N-1) }, \eta_N
} \bigr
\rrvert^{ p } + \frac{ 1 }{ 2 N } \bigl\llvert y_{ 2^{ (L_N-1) }
}^{ 2^{ (L_N-1) }, \eta_N
}
\bigr\rrvert^{ p }
\\
&&\qquad\quad{} - \sum_{ l = 1 }^{ L_N - 1 } \biggl(
\frac{ 2^{(l-1)} }{ T } \biggr)^{ {p}/{4} } N^{
(
p \cdot5^{
(
2^{ (l-1) }
)
}
)
} - \sum
_{ l = L_N+1 }^{ \ld( N ) } \biggl( \frac{ 2^l }{ T }
\biggr)^{
{p}/{4}
} \nonumber
\end{eqnarray}
on
$
\{
\theta_N \geq
\eta_N
\}
\cap
\{
\eta_N \leq
2^{ (L_N + 1)/4} T^{ -1/4 }
\}
\cap
\{
N_1 \leq N
\}
$
for all
$ N \in\{2^1,2^2,2^3,\ldots\} $.
By definition of $\eta_N$ and of $L_N$
we have
$\eta_N\geq2^{L_N/4} T^{ -1/4 }$ on $ \{ N_1 \leq N \} $ [see~\eqref{eqA1}]
for all
$ N \in\{2^1,2^2,2^3,\ldots\} $.
Consequently we get
the inequality $ \eta_N \geq2^{L_N/4} T^{ -1/4 }
( 1 + 5^{  ( -\delta\cdot2^{ (L_N-1) }  ) }  ) $
on $  \{ N_1 \leq N  \} $ [see \eqref{eqA3}]
for all
$ N \in\{2^1,2^2,2^3,\ldots\} $.
Lemmas~\ref{linitialmultiple}
and~\ref{lmonotonicityy}
hence yield
%
%e62 #&#
\begin{eqnarray}
%
% \lefteqn{
&&\sum_{ l = 0 }^{ \ld( N ) }
\frac{ 2^l }{ N } \sum_{ k = 1 }^{{ N }/{ 2^l } } \bigl
\llvert y_{ 2^l }^{ 2^l, \xi^{l,k} } \bigr\rrvert^{ p } - \sum
_{ l = 1 }^{ \ld( N ) } \frac{ 2^{l} }{ N } \sum
_{ k = 1 }^{{ N }/{ 2^{l} } } \bigl\llvert y_{ 2^{ (l-1) } }^{
2^{ (l-1) }, \xi^{l,k}
}
\bigr\rrvert^{ p } % }
\nonumber\\
&&\qquad \geq \biggl(
\frac{ 1 }{ 2 N } \bigl( 1 + 4^{  ( - 2^{ (L_N-1) }  ) } \bigr)^{
(
p \cdot5^{  (
2^{ ( L_N - 1 ) }
) }
)
} - 1 \biggr)
\bigl\llvert y_{ 2^{ (L_N-1) }
}^{ 2^{ (L_N-1) }, \eta_N
} \bigr\rrvert^{ p }
\nonumber
\\[-8pt]
\\[-8pt]
\nonumber
&&\qquad\quad{} + \frac{ 1 }{ 2 N } \bigl\llvert y_{ 2^{ (L_N-1) }
}^{ 2^{ (L_N-1) },
(
{ 2^{ L_N } }/{ T }
)^{{1}/{4} }
(
1
+ 5^{  (- \delta2^{ (L_N - 1) }  )
}
)
} \bigr
\rrvert^{ p }
\\
&&\qquad\quad{} - L_N \biggl( \frac{ 2^{ (L_N - 2) } }{ T } \biggr)^{ {p}/{4} }
N^{
(
p \cdot5^{  ( 2^{ (L_N-2) }  ) }
)
} - \ld( N ) \biggl( \frac{2^{\ld(N)}}{T} \biggr)^{ {p}/{4} }\nonumber
\end{eqnarray}
on
$
\{
\theta_N \geq
\eta_N
\}
\cap
\{
\eta_N \leq
2^{ (L_N + 1)/4} T^{ -1/4 }
\}
\cap
\{
N_1 \leq N
\}
$
for all
$ N \in\{2^1,2^2,2^3,\ldots\} $.
Lemma~\ref{lloweresti} and $ L_N \leq\operatorname{ld}( N ) $
therefore imply
\begin{eqnarray*}
&& \sum_{ l = 0 }^{ \ld( N ) }
\frac{ 2^l }{ N } \sum_{ k = 1 }^{{ N }/{ 2^l } } \bigl
\llvert y_{ 2^l }^{ 2^l, \xi^{l,k} } \bigr\rrvert^{ p } - \sum
_{ l = 1 }^{ \ld( N ) } \frac{ 2^{l} }{ N } \sum
_{ k = 1 }^{{ N }/{ 2^{l} } } \bigl\llvert y_{ 2^{ (l-1) } }^{
2^{ (l-1) }, \xi^{l,k}
}
\bigr\rrvert^{ p }
\\
&&\qquad \geq \biggl( \frac{ 1 }{ 2 N } \bigl( 1 + 2^{  ( - 2^{ L_N }  ) }
\bigr)^{
(
p \cdot5^{  (
2^{ (L_N-1) }
) }
)
} - 1 \biggr) \bigl\llvert y_{ 2^{ (L_N-1) }
}^{ 2^{ (L_N-1) }, \eta_N
}
\bigr\rrvert^{ p }
\\
&&\qquad\quad{} + \frac{ 1 }{ 2 N } \biggl( \frac{ 2^{ L_N } }{ T } \biggr)^{ {p}/{4} } \exp
\biggl( \frac{ p }{ 2 } \cdot 5^{
(
(1-\dl) 2^{ (L_N - 1) }
)
} \biggr) % \\ &
\\
&&\qquad\quad{}- \ld( N )
\biggl( \frac{ 2^{ L_N } }{ T } \biggr)^{ {p}/{4} } N^{
(
p \cdot5^{  ( 2^{ (L_N-2) }  ) }
)
} - \ld( N )
N^{ { p }/{ 4 } } T^{ - { p }/{ 4 } } %
\end{eqnarray*}
on
$
\{
\theta_N \geq
\eta_N
\}
\cap
\{
\eta_N \leq
2^{ (L_N + 1)/4} T^{ -1/4 }
\}
\cap
\{
N_1 \leq N
\}
$
for all
$ N \in\{2^1,2^2,2^3,\ldots\} $.
The inequalities $ 1 \leq L_N \leq\ld(N) $
and $ 1 + 2^{ (-x) } \geq\exp( 2^{(-x-1) } ) $
for all $ x \in[0, \infty) $ hence give
%
%e63 #&#
\begin{eqnarray}
&& \sum_{ l = 0 }^{ \ld( N ) }
\frac{ 2^l }{ N } \sum_{ k = 1 }^{{ N }/{ 2^l } } \bigl
\llvert y_{ 2^l }^{ 2^l, \xi^{l,k} } \bigr\rrvert^{ p } - \sum
_{ l = 1 }^{ \ld( N ) } \frac{ 2^{l} }{ N } \sum
_{ k = 1 }^{{ N }/{ 2^{l} } } \bigl\llvert y_{ 2^{ (l-1) } }^{
2^{ (l-1) }, \xi^{l,k}
}
\bigr\rrvert^{ p }
\nonumber\\
&&\qquad \geq \biggl( \frac{ 1 }{ 2 N } \cdot \exp \bigl( 2^{  ( - 2^{ L_N } - 1  ) } \cdot p
\cdot5^{  (
2^{ (L_N-1) }
) } \bigr) - 1 \biggr) \bigl\llvert y_{ 2^{ (L_N-1) }
}^{ 2^{ (L_N-1) }, \eta_N
}
\bigr\rrvert^{ p }
\nonumber
\\[-8pt]
\\[-8pt]
\nonumber
&&\qquad\quad{} + \frac{ 1 }{ 2 N T^{ { p }/{ 4 } } } \cdot \exp \biggl( \frac{ p }{ 2 } \cdot
5^{
(
(1-\dl) 2^{ (L_N - 1) }
)
} \biggr) % \\ &
\\
&&\qquad\quad{}- 2 \ld( N ) N^{ { p }/{ 4 } }
T^{ - { p }/{ 4 } } \cdot N^{
(
p \cdot5^{  ( 2^{ (L_N-2) }  ) }
)
}\nonumber %
\end{eqnarray}
on
$
\{
\theta_N \geq
\eta_N
\}
\cap
\{
\eta_N \leq
2^{ (L_N + 1)/4} T^{ -1/4 }
\}
\cap
\{
N_1 \leq N
\}
$
for all
$ N \in\{2^1,2^2,2^3,\ldots\} $.
This shows
\begin{eqnarray*}
&& \sum_{ l = 0 }^{ \ld( N ) }
\frac{ 2^l }{ N } \sum_{ k = 1 }^{{ N }/{ 2^l } } \bigl
\llvert y_{ 2^l }^{ 2^l, \xi^{l,k} } \bigr\rrvert^{ p } - \sum
_{ l = 1 }^{ \ld( N ) } \frac{ 2^{l} }{ N } \sum
_{ k = 1 }^{{ N }/{ 2^{l} } } \bigl\llvert y_{ 2^{ (l-1) } }^{
2^{ (l-1) }, \xi^{l,k}
}
\bigr\rrvert^{ p }
\\
& &\qquad\geq \biggl( \frac{1}{ 2 N } \cdot \exp \biggl( \frac{ p }{ 2 } \cdot
\biggl( \frac{5}{4} \biggr)^{
(
2^{ ( L_N - 1 )}
)
} \biggr) - 1 \biggr) \bigl
\llvert y_{ 2^{ (L_N-1) }
}^{ 2^{ (L_N-1) }, \eta_N
} \bigr\rrvert^{ p }
\\
&&\qquad\quad{} + \frac{ 1 }{ 2 N
T^{ {p}/{4} }
} \cdot \exp \biggl( \frac{ p }{ 2 } \cdot
5^{
(
(1-\dl) 2^{ (L_N - 1) }
)
} \biggr) \\[-1pt]
&&\qquad\quad{}- T^{ - { p }/{ 4 } } \cdot N^{
(
1 + { p }/{ 4 } +
p \cdot5^{  ( 2^{ (L_N-2) }  ) }
)
}
\\[-1pt]
&&\qquad \geq \biggl( \frac{1}{ 2 N } \cdot \exp \biggl( \frac{ p }{ 2 } \cdot
\biggl( \frac{5}{4} \biggr)^{
(
2^{ ( L_N - 1 )}
)
} \biggr) - 1 \biggr) \bigl
\llvert y_{ 2^{ (L_N-1) }
}^{ 2^{ (L_N-1) }, \eta_N
} \bigr\rrvert^{ p }
\\[-1pt]
&&\qquad\quad{} + \inf_{
x \in [
2^{ (L_N - 1) } ,
\infty
)
} \biggl[ \exp \biggl( \frac{p}{2} \cdot
5^{  (  (1-\dl ) x  ) } - \operatorname{ln}(2N) \biggr)\\[-1pt]
&&\hspace*{83pt}\qquad{} - \exp \biggl(
\operatorname{ln}(N) \biggl( 1 + \frac{p}{4} + p \cdot 5^{ {x}/{2} }
\biggr) \biggr) \biggr] \cdot T^{ - {p}/{4} } %
\end{eqnarray*}
on
$
\{
\theta_N \geq
\eta_N
\}
\cap
\{
\eta_N \leq
2^{ (L_N + 1)/4} T^{ -1/4 }
\}
\cap
\{
N_1 \leq N
\}
$
for all
$ N \in\{2^1,2^2,2^3,\ldots\} $
and, using the estimate $ 2^{ (L_N -1 ) } \geq\sigmab^2\sqrt{T} \ln(N)
$ on $  \{ N_1 \leq N  \} $
[see \eqref{eqA1}],
%
%e64 #&#
\begin{eqnarray}\label{eqbeforegeqzero}
&& \sum_{ l = 0 }^{ \ld( N ) } \frac{ 2^l }{ N } \sum
_{ k = 1 }^{{ N }/{ 2^l } } \bigl\llvert
y_{ 2^l }^{ 2^l, \xi^{l,k} } \bigr\rrvert^{ p } - \sum
_{ l = 1 }^{ \ld( N ) } \frac{ 2^{l} }{ N } \sum
_{ k = 1 }^{{ N }/{ 2^{l} } } \bigl\llvert y_{ 2^{ (l-1) } }^{
2^{ (l-1) }, \xi^{l,k}
}
\bigr\rrvert^{ p }
\nonumber
\\[-1pt]
&&\qquad \geq \biggl( \frac{ 1 }{ 2 N } \cdot \exp \biggl(
\frac{ p }{ 2 } \cdot \biggl( \frac{5}{4} \biggr)^{
(
\bar{\sigma}^2 \sqrt{T}
\operatorname{ln}(N)
)
}
\biggr) - 1 \biggr) \bigl\llvert y_{ 2^{ (L_N-1) }
}^{ 2^{ (L_N-1) }, \eta_N
} \bigr
\rrvert^{ p }
\nonumber
\\[-8pt]
\\[-8pt]
\nonumber
& &\qquad\quad{}+ \inf_{
x \in [
\bar{\sigma}^2 \sqrt{T}
\operatorname{ln}(N) ,
\infty
)
} \biggl[ \exp \biggl( \frac{p}{2} \cdot
5^{  (  (1-\dl ) x  ) } - \operatorname{ln}(2N) \biggr) \\[-1pt]
&&\hspace*{98pt}\qquad{}- \exp \biggl(
\operatorname{ln}(N) \biggl( 1 + \frac{p}{4} + p \cdot 5^{ {x}/{2} }
\biggr) \biggr) \biggr] \cdot T^{ - {p}/{4} }
\nonumber
\end{eqnarray}
on
$
\{
\theta_N \geq
\eta_N
\}
\cap
\{
\eta_N \leq
2^{ (L_N + 1)/4} T^{ -1/4 }
\}
\cap
\{
N_1 \leq N
\}
$
for all
$ N \in\{2^1,2^2,2^3,\ldots\} $.
It follows from
%
%e65 #&#
\begin{equation}
\lim_{N\to\infty} \frac{1}{2N} \exp \biggl( \frac{ p }{ 2 } \cdot
N^{
(
\bar{\sigma}^2 \sqrt{T}
\operatorname{ln}( 5 / 4 )
)
} \biggr) =\infty
\end{equation}
that there exists an $N_2\in\{2^1,2^2,2^3,\ldots\}$ such that
%
%e66 #&#
\begin{equation}
\frac{1}{2N} \exp \biggl( \frac{ p }{ 2 } \cdot N^{
(
\bar{\sigma}^2 \sqrt{T}
\operatorname{ln}( 5 / 4 )
)
}
\biggr) -1 \geq0
\end{equation}
for all $N\in[N_2,\infty)$.
Using this, we deduce from~\eqref{eqbeforegeqzero}
%
%e67 #&#
\begin{eqnarray}
\label{eqresultcase2} %
&& \sum_{ l = 0 }^{ \ld( N ) }
\frac{ 2^l }{ N } \sum_{ k = 1 }^{{ N }/{ 2^l } } \bigl
\llvert y_{ 2^l }^{ 2^l, \xi^{l,k} } \bigr\rrvert^{ p } - \sum
_{ l = 1 }^{ \ld( N ) } \frac{ 2^{l} }{ N } \sum
_{ k = 1 }^{{ N }/{ 2^{l} } } \bigl\llvert y_{ 2^{ (l-1) } }^{
2^{ (l-1) }, \xi^{l,k}
}
\bigr\rrvert^{ p }
\nonumber\\[-1pt]
&&\qquad \geq \biggl( \frac{ 1 }{ 2 N } \cdot \exp \biggl( \frac{ p }{ 2 } \cdot
N^{
(
\bar{\sigma}^2 \sqrt{T}
\operatorname{ln}( 5 / 4 )
)
} \biggr) - 1 \biggr) \bigl\llvert y_{ 2^{ (L_N-1) }
}^{ 2^{ (L_N-1) }, \eta_N
}
\bigr\rrvert^{ p } + r( N ) \cdot T^{ - {p}/{4} }\\[-1pt]
&&\qquad \geq r( N ) \cdot
T^{ - {p}/{4} }\nonumber\vadjust{\goodbreak} %
\end{eqnarray}
on
$
\{
\theta_N \geq
\eta_N
\}
\cap
\{
\eta_N \leq
2^{ (L_N + 1)/4} T^{ -1/4 }
\}
\cap
\{
N_1 \leq N
\}
$
for all
$ N \in\{ N_2, 2^1 N_2,\break  2^2 N_2, \ldots\} $.

Next, we analyze the behavior
of the multilevel Monte Carlo
Euler approximations on the events
$
\{
2^{L_N/4} T^{-1/4}
<
\theta_N <
\eta_N \leq
2^{ (L_N + 1)/4} T^{ -1/4 }
\}
\cap
\{
N_1 \leq N
\}
$
for
$ N \in\{2^1,2^2,2^3,\ldots\} $.
% Here we exploit that $\xi^{L_N,k}\leq2^{ (L_N + 1)/4} T^{ -1/4 }$
% for all
% $ k \in\{ 1, 2, \ldots, \frac{N}{2^{L_N}} \} $
% on
% $
% \big\{
% \eta_N \leq
% 2^{ (L_N + 1)/4} T^{ -1/4 }
% \big\}
% $.
% Moreover we have
% $\xi^{l,k}\leq2^{l/4}T^{-1/4} \leq2^{(l+1)/4}T^{-1/4}$
% for all
% $ k \in\{ 1, 2, \ldots, \frac{N}{2^{l}} \} $,
% $ l \in\{ L_N+1, L_N+2, \ldots, \operatorname{ld}( N ) \} $
% by definition~\eqref{eqdefLN} of $L_N$
% for all
% $ N \in\{2^1,2^2,2^3,\ldots\} $.
Note that Lemma~\ref{lmonotonicityy} and the inequality
$
| \xi^{ l,k} |
< 2^{ (l-1)/4} T^{ -1/4 } N
$ on
$  \{ N_1 \leq N  \} $
for all $ k \in\{ 1, 2, \ldots, \frac{N}{2^l} \} $,
$ l \in\{ 0, 1, 2, \ldots, \operatorname{ld}( N ) \} $ [see~\eqref{eqA2}]
imply
%
%e68 #&#
\begin{eqnarray}
&& \Biggl\llvert \sum_{ l = 0 }^{ \ld( N ) }
\frac{ 2^l }{ N } \sum_{ k = 1 }^{{ N }/{ 2^l } } \bigl
\llvert y_{ 2^l }^{ 2^l, \xi^{l,k} } \bigr\rrvert^{ p } - \sum
_{ l = 1 }^{ \ld( N ) } \frac{ 2^{l} }{ N } \sum
_{ k = 1 }^{{ N }/{ 2^{l} } } \bigl\llvert y_{ 2^{ (l-1) } }^{
2^{ (l-1) }, \xi^{l,k}
}
\bigr\rrvert^{ p } \Biggr\rrvert
\nonumber\\
&&\qquad \geq \sum_{ l = 1 }^{ \ld( N ) } \frac{ 2^{l} }{ N }
\sum_{ k = 1 }^{{ N }/{ 2^{l} } } \bigl\llvert
y_{ 2^{ (l-1) } }^{
2^{ (l-1) }, \xi^{l,k}
} \bigr\rrvert^{ p } - \sum
_{ l = 0 }^{ \ld( N ) } \frac{ 2^l }{ N } \sum
_{ k = 1 }^{{ N }/{ 2^l } } \bigl\llvert y_{ 2^l }^{ 2^l, \xi^{l,k} }
\bigr\rrvert^{ p }
\nonumber\\
&&\qquad \geq \frac{1}{N} \bigl\llvert y_{ 2^{ (L_N - 1) }
}^{ 2^{ (L_N - 1) }, \eta_N
} \bigr
\rrvert^{ p } - \bigl\llvert y_{ 2^{ (L_N-1) }
}^{ 2^{ (L_N-1) }, \theta_N
} \bigr
\rrvert^{ p } - \sum_{ l = 0 }^{ L_N - 2 }
\frac{ 2^{l} }{ N } \sum_{ k = 1 }^{{ N }/{ 2^{l} } } \bigl
\llvert y_{ 2^{ l } }^{
2^{ l }, \xi^{l,k}
} \bigr\rrvert^{ p }
\nonumber
\\[-8pt]
\\[-8pt]
\nonumber
&&\qquad\quad{} - \sum
_{ l = L_N }^{ \ld(N) } \frac{ 2^{l} }{ N } \sum
_{ k = 1 }^{{ N }/{ 2^{l} } } \bigl\llvert
y_{ 2^{ l } }^{
2^{ l }, \xi^{l,k}
} \bigr\rrvert^{ p }
\\
&&\qquad \geq \frac{1}{2 N} \bigl\llvert y_{ 2^{ (L_N - 1) }
}^{ 2^{ (L_N - 1) }, \eta_N
} \bigr
\rrvert^{ p } - \bigl\llvert y_{ 2^{ (L_N-1) }
}^{ 2^{ (L_N-1) }, \theta_N
} \bigr
\rrvert^{ p } + \frac{1}{2 N} \bigl\llvert y_{ 2^{ (L_N - 1) }
}^{ 2^{ (L_N - 1) }, \eta_N
}
\bigr\rrvert^{ p } \nonumber\\
&&\qquad\quad{}- \sum_{ l = 0 }^{ L_N - 2 }
\bigl\llvert y_{ 2^{ l } }^{
2^{ l },
(
{ 2^{ (l-1) } }/{ T }
)^{{1}/{4} }
N
} \bigr\rrvert^{ p } -
\sum_{ l = L_N }^{ \ld( N ) } \bigl\llvert
y_{ 2^{ l } }^{
2^{ l },
(
{ 2^{(l+1)} }/{T}
)^{{1}/{4} }
} \bigr\rrvert^{ p }\nonumber %
\end{eqnarray}
on
$
\{
2^{ L_N/4} T^{ -1/4 } <
\theta_N <
\eta_N \leq
2^{ (L_N + 1)/4} T^{ -1/4 }
\}
\cap
\{
N_1 \leq N
\}
$
for all
$ N \in\{2^1,2^2,\break 2^3,\ldots\} $.
Therefore
Lemma~\ref{lmonotonicityy},
the inequality
$
\eta_N \geq
(
1 + 4^{  ( -2^{ (L_N - 1) }  ) }
) \theta_N
$
on $  \{ \theta_N < \eta_N  \} \cap \{ N_1 \leq N
\} $
[see \eqref{eqA4}] and
Lemma~\ref{lstability} result in
%
%e69 #&#
\begin{eqnarray}
&& \Biggl\llvert \sum_{ l = 0 }^{ \ld( N ) }
\frac{ 2^l }{ N } \sum_{ k = 1 }^{{ N }/{ 2^l } } \bigl
\llvert y_{ 2^l }^{ 2^l, \xi^{l,k} } \bigr\rrvert^{ p } - \sum
_{ l = 1 }^{ \ld( N ) } \frac{ 2^{l} }{ N } \sum
_{ k = 1 }^{{ N }/{ 2^{l} } } \bigl\llvert y_{ 2^{ (l-1) } }^{
2^{ (l-1) }, \xi^{l,k}
}
\bigr\rrvert^{ p } \Biggr\rrvert
\nonumber\\
&&\qquad \geq \frac{1}{2 N} \bigl\llvert y_{ 2^{ (L_N - 1) }
}^{ 2^{ (L_N - 1) },
(
1 +
4^{  ( - 2^{ (L_N-1) }  ) }
)
\theta_N
} \bigr
\rrvert^{ p } - \bigl\llvert y_{ 2^{ (L_N-1) }
}^{ 2^{ (L_N-1) }, \theta_N
} \bigr
\rrvert^{ p } + \frac{ 1 }{ 2 N } \bigl\llvert y_{ 2^{ (L_N-1) }
}^{ 2^{ (L_N-1) }, \eta_N
}
\bigr\rrvert^{ p }
\\
&&\qquad\quad{} - \sum_{ l = 0 }^{ L_N - 2 } \bigl\llvert
y_{ 2^{ l } }^{
2^{ l },
(
{ 2^{ l } }/{ T }
)^{{1}/{4} }
N
} \bigr\rrvert^{ p } - \sum
_{ l = L_N }^{ \ld( N ) } \biggl( \frac{ 2^{ ( l + 1 ) } }{ T }
\biggr)^{
{p}/{4}
} \nonumber%
\end{eqnarray}
on
$
\{
2^{ L_N/4} T^{ -1/4 } <
\theta_N <
\eta_N \leq
2^{ (L_N + 1)/4} T^{ -1/4 }
\}
\cap
\{
N_1 \leq N
\}
$
for all
$ N \in \{2^1,2^2,\break2^3,\ldots\} $.
Lemmas~\ref{linitialmultiple},
\ref{lmonotonicityy} and the estimate $ \eta_N \geq2^{L_N/4}
T^{ -1/4 }\times
( 1 + 5^{  ( -\delta\cdot2^{ (L_N-1) }  ) }  ) $
on $  \{ N_1 \leq N  \} $ [see \eqref{eqA1} and \eqref{eqA3}]
and Lemma~\ref{lupperesti}
hence yield
\begin{eqnarray*}
&& \Biggl\llvert \sum_{ l = 0 }^{ \ld( N ) }
\frac{ 2^l }{ N } \sum_{ k = 1 }^{{ N }/{ 2^l } } \bigl
\llvert y_{ 2^l }^{ 2^l, \xi^{l,k} } \bigr\rrvert^{ p } - \sum
_{ l = 1 }^{ \ld( N ) } \frac{ 2^{l} }{ N } \sum
_{ k = 1 }^{{ N }/{ 2^{l} } } \bigl\llvert y_{ 2^{ (l-1) } }^{
2^{ (l-1) }, \xi^{l,k}
}
\bigr\rrvert^{ p } \Biggr\rrvert
\\
&&\qquad \geq \biggl( \frac{1}{ 2 N } \bigl( 1 + 4^{  ( - 2^{ (L_N-1) }  ) }
\bigr)^{
(
p \cdot5^{  (
2^{ ( L_N - 1 ) }
) }
)
} - 1 \biggr) \bigl\llvert y_{ 2^{ (L_N-1) }
}^{ 2^{ (L_N-1) }, \theta_N
}
\bigr\rrvert^{ p }
\\
&&\qquad\quad{} + \frac{ 1 }{ 2 N } \bigl\llvert y_{ 2^{ (L_N-1) }
}^{ 2^{ (L_N-1) },
(
{ 2^{ L_N } }/{ T }
)^{{1}/{4} }
(
1
+ 5^{ - \delta\cdot
2^{ (L_N - 1) }
}
)
} \bigr
\rrvert^{ p } \\
&&\qquad\quad{}- \sum_{ l = 0 }^{ L_N - 2 }
\biggl( \frac{ 2^{ l } }{ T } \biggr)^{ {p}/{4} } N^{
(
p \cdot5^{
(
2^{ l }
)
}
)
} -
2^{ {p}/{4} } \ld(N) N^{ {p}/{4} } T^{ - {p}/{4} } %
\end{eqnarray*}
on
$
\{
2^{ L_N/4} T^{ -1/4 }
<
\theta_N <
\eta_N \leq
2^{ (L_N + 1)/4} T^{ -1/4 }
\}
\cap
\{
N_1 \leq N
\}
$
for all
$ N \in\{2^1,2^2,\break 2^3,\ldots\} $.
Therefore
Lemma~\ref{lloweresti} implies
%
%e70 #&#
\begin{eqnarray}
&& \Biggl\llvert \sum_{ l = 0 }^{ \ld( N ) }
\frac{ 2^l }{ N } \sum_{ k = 1 }^{{ N }/{ 2^l } } \bigl
\llvert y_{ 2^l }^{ 2^l, \xi^{l,k} } \bigr\rrvert^{ p } - \sum
_{ l = 1 }^{ \ld( N ) } \frac{ 2^{l} }{ N } \sum
_{ k = 1 }^{{ N }/{ 2^{l} } } \bigl\llvert y_{ 2^{ (l-1) } }^{
2^{ (l-1) }, \xi^{l,k}
}
\bigr\rrvert^{ p } \Biggr\rrvert
\nonumber\\
&&\qquad \geq \biggl( \frac{ 1 }{ 2 N } \bigl( 1 + 2^{  ( - 2^{ L_N }  ) }
\bigr)^{
(
p \cdot5^{  (
2^{ ( L_N - 1 ) }
) }
)
} - 1 \biggr) \bigl\llvert y_{ 2^{ (L_N-1) }
}^{ 2^{ (L_N-1) }, \theta_N
}
\bigr\rrvert^{ p }
\nonumber
\\[-8pt]
\\[-8pt]
\nonumber
&&\qquad\quad{} + \frac{ 1 }{ 2 N } \biggl( \frac{ 2^{ L_N } }{ T } \biggr)^{{p}/{4} } \exp
\biggl( \frac{ p }{ 2 } \cdot 5^{
(
(1-\dl) 2^{ (L_N - 1) }
)
} \biggr)\\
&&\quad\qquad{} - 2 \ld(N)
N^{{p}/{4} } T^{ - {p}/{4} } N^{
(
p \cdot5^{  ( 2^{ (L_N-2) }  ) }
)
} \nonumber%
\end{eqnarray}
on
$
\{
2^{ L_N/4 } T^{ -1/4 }
<
\theta_N <
\eta_N \leq
2^{ (L_N + 1)/4} T^{ -1/4 }
\}
\cap
\{
N_1 \leq N
\}
$
for all
$ N \in\{ 2^1, 2^2,\break  2^3, \ldots\} $.
The inequality
$ 1 + 2^{ (-x) } \geq\exp( 2^{(-x-1) } ) $
for all $ x \in[0, \infty) $ hence shows
%
%e71 #&#
\begin{eqnarray}
&& \Biggl\llvert \sum_{ l = 0 }^{ \ld( N ) }
\frac{ 2^l }{ N } \sum_{ k = 1 }^{{ N }/{ 2^l } } \bigl
\llvert y_{ 2^l }^{ 2^l, \xi^{l,k} } \bigr\rrvert^{ p } - \sum
_{ l = 1 }^{ \ld( N ) } \frac{ 2^{l} }{ N } \sum
_{ k = 1 }^{{ N }/{ 2^{l} } } \bigl\llvert y_{ 2^{ (l-1) } }^{
2^{ (l-1) }, \xi^{l,k}
}
\bigr\rrvert^{ p } \Biggr\rrvert
\nonumber\\
&&\qquad \geq \biggl( \frac{1}{2N} \cdot \exp \bigl( 2^{  ( - 2^{ L_N } - 1  ) } \cdot p
\cdot5^{  (
2^{ (L_N-1) }
) } \bigr) - 1 \biggr) \bigl\llvert y_{ 2^{ (L_N-1) }
}^{ 2^{ (L_N-1) }, \theta_N
}
\bigr\rrvert^{ p }
\\
&&\qquad\quad{} + \frac{ 1 }{ 2 N
T^{ { p }/{ 4 } }
} \cdot \exp \biggl( \frac{ p }{ 2 } \cdot
5^{
(
(1-\dl) 2^{ (L_N - 1) }
)
} \biggr) - T^{ - {p}/{4} } \cdot N^{
(
1 + { p }/{ 4 } +
p \cdot5^{  ( 2^{ (L_N-2) }  ) }
)
} \nonumber%
\end{eqnarray}
on
$
\{
2^{ L_N/4} T^{ -1/4 }
<
\theta_N <
\eta_N \leq
2^{ (L_N + 1)/4} T^{ -1/4 }
\}
\cap
\{
N_1 \leq N
\}
$
for all
$ N \in\{ 2^1, 2^2,\break  2^3, \ldots\} $.
Consequently
\begin{eqnarray*}
&& \Biggl\llvert \sum_{ l = 0 }^{ \ld( N ) }
\frac{ 2^l }{ N } \sum_{ k = 1 }^{{ N }/{ 2^l } } \bigl
\llvert y_{ 2^l }^{ 2^l, \xi^{l,k} } \bigr\rrvert^{ p } - \sum
_{ l = 1 }^{ \ld( N ) } \frac{ 2^{l} }{ N } \sum
_{ k = 1 }^{{ N }/{ 2^{l} } } \bigl\llvert y_{ 2^{ (l-1) } }^{
2^{ (l-1) }, \xi^{l,k}
}
\bigr\rrvert^{ p } \Biggr\rrvert
\\
& &\qquad\geq \biggl( \frac{1}{2 N} \cdot \exp \biggl( \frac{ p }{ 2 } \cdot
\biggl( \frac{5}{4} \biggr)^{
(
2^{ ( L_N - 1 )}
)
} \biggr) - 1 \biggr) \bigl
\llvert y_{ 2^{ (L_N-1) }
}^{ 2^{ (L_N-1) }, \theta_N
} \bigr\rrvert^{ p }
\\
&&\qquad\quad{} + \inf_{
x \in [
2^{ (L_N-1) } ,
\infty
)
} \biggl[ \exp \biggl( \frac{p}{2} \cdot
5^{  (  (1-\dl ) x  ) } - \operatorname{ln}(2N) \biggr)\\
&&\hspace*{83pt}\qquad{} - \exp \biggl(
\operatorname{ln}(N) \biggl( 1 + \frac{p}{4} + p \cdot 5^{{x}/{2} }
\biggr) \biggr) \biggr] \cdot T^{ - {p}/{4} } %
\end{eqnarray*}
on
$
\{
2^{ L_N /4} T^{ -1/4 }
<
\theta_N <
\eta_N \leq
2^{ (L_N + 1)/4} T^{ -1/4 }
\}
\cap
\{
N_1 \leq N
\}
$
for all
$ N \in\{ 2^1,\break 2^2,  2^3, \ldots\} $.
The estimate $ 2^{ (L_N -1 ) }
\geq\bar{ \sigma}^2
\sqrt{T} \ln(N) $
on $  \{ N_1 \leq N  \} $
[see \eqref{eqA1}] therefore implies
\begin{eqnarray*}
&& \Biggl\llvert \sum_{ l = 0 }^{ \ld( N ) }
\frac{ 2^l }{ N } \sum_{ k = 1 }^{{ N }/{ 2^l } } \bigl
\llvert y_{ 2^l }^{ 2^l, \xi^{l,k} } \bigr\rrvert^{ p } - \sum
_{ l = 1 }^{ \ld( N ) } \frac{ 2^{l} }{ N } \sum
_{ k = 1 }^{{ N }/{ 2^{l} } } \bigl\llvert y_{ 2^{ (l-1) } }^{
2^{ (l-1) }, \xi^{l,k}
}
\bigr\rrvert^{ p } \Biggr\rrvert
\\
&&\qquad \geq \biggl( \frac{1}{2 N} \cdot \exp \biggl( \frac{ p }{ 2 } \cdot
\biggl( \frac{5}{4} \biggr)^{
(
\bar{ \sigma}^2 \sqrt{T}
\operatorname{ln}(N)
)
} \biggr) - 1 \biggr) \bigl
\llvert y_{ 2^{ (L_N-1) }
}^{ 2^{ (L_N-1) }, \theta_N
} \bigr\rrvert^{ p }
\\
& &\qquad\quad{}+ \inf_{
x \in [
\bar{ \sigma}^2 \sqrt{T}
\operatorname{ln}(N) ,
\infty
)
} \biggl[ \exp \biggl( \frac{p}{2} \cdot
5^{  (  (1-\dl ) x  ) } - \operatorname{ln}(2N) \biggr)\\
&&\hspace*{98pt}\qquad{} - \exp \biggl(
\operatorname{ln}(N) \biggl( 1 + \frac{p}{4} + p \cdot 5^{{x}/{2} }
\biggr) \biggr) \biggr] \cdot T^{ - {p}/{4} } %
\end{eqnarray*}
on
$
\{
2^{ L_N/4} T^{ -1/4 }
<
\theta_N <
\eta_N \leq
2^{ (L_N + 1)/4} T^{ -1/4 }
\}
\cap
\{
N_1 \leq N
\}
$
for all
$ N \in\{ 2^1,\break 2^2 ,  2^3 \ldots\} $. Finally, we obtain
%
%e72 #&#
\begin{eqnarray}
\label{eqresultcase3} %
&& \Biggl\llvert \sum
_{ l = 0 }^{ \ld( N ) } \frac{ 2^l }{ N } \sum
_{ k = 1 }^{{ N }/{ 2^l } } \bigl\llvert y_{ 2^l }^{ 2^l, \xi^{l,k} }
\bigr\rrvert^{ p } - \sum_{ l = 1 }^{ \ld( N ) }
\frac{ 2^{l} }{ N } \sum_{ k = 1 }^{{ N }/{ 2^{l} } } \bigl
\llvert y_{ 2^{ (l-1) } }^{
2^{ (l-1) }, \xi^{l,k}
} \bigr\rrvert^{ p } \Biggr
\rrvert
\nonumber\\
&&\qquad \geq \biggl( \frac{ 1 }{ 2 N } \cdot \exp \biggl( \frac{ p }{ 2 } \cdot
N^{
(
\bar{ \sigma}^2
\sqrt{ T }
\operatorname{ln}( 5 / 4 )
)
} \biggr) - 1 \biggr) \bigl\llvert y_{ 2^{ (L_N-1) }
}^{ 2^{ (L_N-1) }, \theta_N
}
\bigr\rrvert^{ p } + r( N ) \cdot T^{ - {p}/{4} } \\
&&\qquad\geq r( N ) \cdot
T^{ - {p}/{4} } \nonumber%
\end{eqnarray}
on
$
\{
2^{ L_N/4} T^{ -1/4 }
<
\theta_N <
\eta_N \leq
2^{ (L_N + 1)/4} T^{ -1/4 }
\}
\cap
\{
N_1 \leq N
\}
$
for all
$ N \in \{ N_2,\break 2^1 N_2, 2^2 N_2, \ldots\} $.

Finally, we analyze the behavior
of the multilevel Monte Carlo
Euler \mbox{approximations} on the events
$
\{
\theta_N
\leq
2^{ L_N/4} T^{ -1/4 }
\}
\cap
\{
\theta_N <
\eta_N \leq\break
2^{ (L_N + 1)/4} T^{ -1/4 }
\}
\cap
\{
N_1 \leq N
\}
$
for
$ N \in\{2^1,2^2,2^3,\ldots\} $.
% Here we exploit that
% $\xi^{l,k}\leq2^{(l+1)/4} T^{ -1/4 } $
% for all
% $ k \in\{ 1, 2, \ldots, \frac{N}{2^{l}} \} $
% and all
% $ l \in\{ L_N-1, L_N, L_N+1, \ldots, \ln(N) \} $
% on
% $
% \big\{
% \theta_N
% \leq
% 2^{ L_N /4} T^{ -1/4 }
% \big\}
% \cap
% \big\{
% \eta_N \leq
% 2^{ (L_N + 1)/4} T^{ -1/4 }
% \big\}
% $.
% Moreover, by definition of $\eta_N$ and of $L_N$
% we have
% $\eta_N\geq2^{ L_N/4} T^{ -1/4 }$
% for all
% $ N \in\{2^1,2^2,2^3,\ldots\} $.
% Consequently we get
% the inequality $ \eta_N \geq2^{ L_N/4} T^{ -1/4 }
% \big( 1 + 5^{ \left( -\delta\cdot2^{ (L_N-1) } \right) } \big) $
% on $ \left\{ N_1 \leq N \right\} $ (see \eqref{eqA1} and
% for all
% $ N \in\{2^1,2^2,2^3,\ldots\} $.
Note that Lemma~\ref{lmonotonicityy} and the inequality
$
\llvert  \xi^{ l,k} \rrvert
< 2^{ (l-1)/4} T^{ -1/4 } N
$ on
$  \{ N_1 \leq N  \} $ for all
$ k \in\{ 1, 2, \ldots, \frac{N}{2^l} \} $,
$ l \in\{ 0, 1, \ldots, \operatorname{ld}( N ) \} $
[see~\eqref{eqA2}]
imply
\begin{eqnarray*}
&& \Biggl\llvert \sum_{ l = 0 }^{ \ld( N ) }
\frac{ 2^l }{ N } \sum_{ k = 1 }^{{ N }/{ 2^l } } \bigl
\llvert y_{ 2^l }^{ 2^l, \xi^{l,k} } \bigr\rrvert^{ p } - \sum
_{ l = 1 }^{ \ld( N ) } \frac{ 2^{l} }{ N } \sum
_{ k = 1 }^{{ N }/{ 2^{l} } } \bigl\llvert y_{ 2^{ (l-1) } }^{
2^{ (l-1) }, \xi^{l,k}
}
\bigr\rrvert^{ p } \Biggr\rrvert
\\
&&\qquad \geq \sum_{ l = 1 }^{ \ld( N ) } \frac{ 2^{l} }{ N }
\sum_{ k = 1 }^{{ N }/{ 2^{l} } } \bigl\llvert
y_{ 2^{ (l-1) } }^{
2^{ (l-1) }, \xi^{l,k}
} \bigr\rrvert^{ p } - \sum
_{ l = 0 }^{ \ld( N ) } \frac{ 2^l }{ N } \sum
_{ k = 1 }^{{ N }/{ 2^l } } \bigl\llvert y_{ 2^l }^{ 2^l, \xi^{l,k} }
\bigr\rrvert^{ p }
\\[-1pt]
&&\qquad \geq \frac{ 1 }{ N } \bigl\llvert y_{ 2^{(L_N-1)} }^{ 2^{(L_N-1)}, \eta_N } \bigr
\rrvert^{ p } - \sum_{ l = 0 }^{ L_N-2 }
\frac{ 2^{l} }{ N } \sum_{ k = 1 }^{{ N }/{ 2^{l} } } \bigl
\llvert y_{ 2^{ l } }^{
2^{ l }, \xi^{l,k}
} \bigr\rrvert^{ p } - \sum
_{ l = L_N-1 }^{ \ld( N ) } \frac{ 2^{l} }{ N } \sum
_{ k = 1 }^{{ N }/{ 2^{l} } } \bigl\llvert y_{ 2^{ l } }^{
2^{ l }, \xi^{l,k}
}
\bigr\rrvert^{ p }
\\[-2pt]
&&\qquad \geq \frac{ 1 }{ N } \bigl\llvert y_{ 2^{(L_N-1)} }^{ 2^{(L_N-1)}, ( { 2^{L_N} }/{ T } )^{{1}/{4} }
( 1 + 5^{ ( -\delta\cdot2^{ (L_N-1) } ) }  )
} \bigr
\rrvert^{ p } - \sum_{ l = 0 }^{ L_N-2 }
\bigl\llvert y_{ 2^{ l } }^{
2^{ l } ,
(
{ 2^{(l-1)} }/{ T }
)^{{1}/{4} }
N
} \bigr\rrvert^{ p } \\[-2pt]
&&\qquad\quad{}- \sum
_{ l = L_N-1 }^{ \ld( N ) } \bigl\llvert
y_{ 2^{ l } }^{
2^{ l },
(
{ 2\cdot2^l }/{ T }
)^{
{1}/{4}
}
} \bigr\rrvert^{ p } %
\end{eqnarray*}
on
$
\{
\theta_N
\leq
2^{ L_N /4} T^{ -1/4 }
\}
\cap
\{
\theta_N <
\eta_N \leq
2^{ (L_N + 1)/4} T^{ -1/4 }
\}
\cap
\{
N_1 \leq N
\}
$
for all
$ N \in\{2^1,2^2,2^3,\ldots\} $
and, applying Lemmas~\ref{lloweresti},
\ref{lmonotonicityy},
\ref{lupperesti} and~\ref{lstability},
\begin{eqnarray*}
&& \Biggl\llvert \sum_{ l = 0 }^{ \ld( N ) }
\frac{ 2^l }{ N } \sum_{ k = 1 }^{{ N }/{ 2^l } } \bigl
\llvert y_{ 2^l }^{ 2^l, \xi^{l,k} } \bigr\rrvert^{ p } - \sum
_{ l = 1 }^{ \ld( N ) } \frac{ 2^{l} }{ N } \sum
_{ k = 1 }^{{ N }/{ 2^{l} } } \bigl\llvert y_{ 2^{ (l-1) } }^{
2^{ (l-1) }, \xi^{l,k}
}
\bigr\rrvert^{ p } \Biggr\rrvert
\\[-2pt]
&&\qquad \geq \frac{1}{N} \biggl\llvert \biggl( \frac{ 2^{L_N} }{ T }
\biggr)^{
{1}/{4}
} \sqrt{e}^{
(
5^{ (1-\dl) 2^{(L_N-1)} }
)
} \biggr\rrvert^{ p }
- \sum_{ l = 0 }^{ L_N-2 } \biggl(
\frac{ 2^{l} }{ T } \biggr)^{{p}/{4} } N^{
(
p \cdot5^{
(
2^{ l }
)
}
)
} \\[-2pt]
&&\qquad\quad{}- \sum
_{ l = L_N-1 }^{ \ld( N ) } \biggl( \frac{ 2\cdot2^l }{ T }
\biggr)^{
{p}/{4}
}
\\[-2pt]
&&\qquad \geq N^{ - 1 } T^{ - {p}/{4} } \cdot \exp \biggl( \frac{p}{2}
\cdot 5^{ (1-\dl) 2^{(L_N-1)} } \biggr) - L_N \biggl( \frac{ 2^{(L_N-2)} }{ T }
\biggr)^{{p}/{4} } N^{
(
p \cdot5^{  ( 2^{(L_N-2)}  ) }
)
} \\[-2pt]
&&\qquad\quad{}- \ld( N ) \biggl( \frac{2\cdot2^{\ld(N)}}{T}
\biggr)^{{p}/{4} } %
\end{eqnarray*}
on
$
\{
\theta_N
\leq
2^{ L_N/4} T^{ -1/4 }
\}
\cap
\{
\theta_N <
\eta_N \leq
2^{ (L_N + 1)/4} T^{ -1/4 }
\}
\cap
\{
N_1 \leq N
\}
$
for all
$ N \in\{2^1,2^2,2^3,\ldots\} $.
Therefore, we obtain
\begin{eqnarray*}
&& \Biggl\llvert \sum_{ l = 0 }^{ \ld( N ) }
\frac{ 2^l }{ N } \sum_{ k = 1 }^{{ N }/{ 2^l } } \bigl
\llvert y_{ 2^l }^{ 2^l, \xi^{l,k} } \bigr\rrvert^{ p } - \sum
_{ l = 1 }^{ \ld( N ) } \frac{ 2^{l} }{ N } \sum
_{ k = 1 }^{{ N }/{ 2^{l} } } \bigl\llvert y_{ 2^{ (l-1) } }^{
2^{ (l-1) }, \xi^{l,k}
}
\bigr\rrvert^{ p } \Biggr\rrvert
\\[-2pt]
&&\qquad \geq N^{ - 1 } T^{ - {p}/{4} } \cdot \exp \biggl( \frac{p}{2}
\cdot 5^{ (1-\dl) 2^{(L_N-1)} } \biggr)\\[-2pt]
&&\qquad \quad{} - \ld( N ) N^{{p}/{4} } T^{ - {p}/{4} }
\cdot N^{
(
p \cdot5^{  ( 2^{(L_N-2)}  ) }
)
}\\[-2pt]
&&\qquad\quad{} - \ld( N ) \biggl(\frac{2N}{T} \biggr)^{{p}/{4} }
\\[-2pt]
&&\qquad \geq T^{ - {p}/{4} } \cdot \exp \biggl( \frac{p}{2} \cdot
5^{ (1-\dl) 2^{(L_N-1)} } - \operatorname{ln}( N ) \biggr) - T^{ - {p}/{4} } \cdot
N^{
( 1 + {p/}{4} +
p \cdot5^{  ( 2^{(L_N-2)}  ) }
)
} %
\end{eqnarray*}
on
$
\{
\theta_N
\leq
2^{ L_N /4} T^{ -1/4 }
\}
\cap
\{
\theta_N <
\eta_N \leq
2^{ (L_N + 1)/4} T^{ -1/4 }
\}
\cap
\{
N_1 \leq N
\}
$
and hence, using $2^{(L_N-1)}\geq\sigmab^2\sqrt{T}\ln(N)$ on $ \{
N_1 \leq N
\} $,
%
%e73 #&#
\begin{eqnarray}
\label{eqresultcase4} && \Biggl\llvert \sum
_{ l = 0 }^{ \ld( N ) } \frac{ 2^l }{ N } \sum
_{ k = 1 }^{{ N }/{ 2^l } } \bigl\llvert y_{ 2^l }^{ 2^l, \xi^{l,k} }
\bigr\rrvert^{ p } - \sum_{ l = 1 }^{ \ld( N ) }
\frac{ 2^{l} }{ N } \sum_{ k = 1 }^{{ N }/{ 2^{l} } } \bigl
\llvert y_{ 2^{ (l-1) } }^{
2^{ (l-1) }, \xi^{l,k}
} \bigr\rrvert^{ p } \Biggr
\rrvert
\nonumber\\
&&\qquad \geq \inf_{
x \in[
\sigmab^2\sqrt{T} \operatorname{ln}(N) ,
\infty
)
} \biggl[ \exp \biggl( \frac{p}{2} \cdot
5^{  (  (1-\dl ) x  ) } - \operatorname{ln}(2 N ) \biggr)
\nonumber
\\[-8pt]
\\[-8pt]
\nonumber
&&\hspace*{85pt}\qquad{} - \exp \biggl(
\operatorname{ln}(N) \biggl( 1 + \frac{p}{4} + p \cdot 5^{{x}/{2} }
\biggr) \biggr) \biggr] \cdot T^{ - {p}/{4} } \\
&&\qquad= r(N) \cdot T^{ - {p}/{4}
}\nonumber
\end{eqnarray}
on
$
\{
\theta_N
\leq
2^{ L_N /4} T^{ -1/4 }
\}
\cap
\{
\theta_N <
\eta_N \leq
2^{ (L_N + 1)/4} T^{ -1/4 }
\}
\cap
\{
N_1 \leq N
\}
$
for all $N\in\{2^1,2^2,2^3,\ldots\}$.

Combining \eqref{eqresultcase1}, \eqref{eqresultcase2},
\eqref{eqresultcase3} and \eqref{eqresultcase4} then shows
%
%e74 #&#
\begin{eqnarray}
\label{eqendineq} \Biggl\llvert \sum_{ l = 0 }^{ \ld( N ) }
\frac{ 2^l }{ N } \sum_{ k = 1 }^{{ N }/{ 2^l } } \bigl
\llvert y_{ 2^l }^{ 2^l, \xi^{l,k} } \bigr\rrvert^{ p } - \sum
_{ l = 1 }^{ \ld( N ) } \frac{ 2^{l} }{ N } \sum
_{ k = 1 }^{{ N }/{ 2^{l} } } \bigl\llvert y_{ 2^{ (l-1) } }^{
2^{ (l-1) }, \xi^{l,k}
}
\bigr\rrvert^{ p } \Biggr\rrvert \geq r(N) \cdot T^{ - { p }/{ 4 } }
\end{eqnarray}
on
$ \{ N_1 \leq N \} $
for all
$ N \in\{ N_2, 2^1N_2, 2^2N_2, \ldots\} $.
Equation~\eqref{eqrewrite2}
and inequality~\eqref{eqendineq}
imply
%
%e75 #&#
\begin{eqnarray}
&& \Biggl\llvert \frac{1}{N} \sum
_{ k = 1 }^{ N } \bigl| Y_1^{ 1, 0, k }(\omega)
\bigr|^p + \sum_{ l = 1 }^{ \ld( N ) }
\frac{ 2^l }{ N } \sum_{ k = 1 }^{{ N }/{ 2^l } } \bigl( \bigl|
Y_{ 2^l }^{ 2^l, l, k }(\omega) \bigr|^p - \bigl| Y_{ 2^{ (l-1) } }^{
2^{ (l-1) }, l, k
}(
\omega) \bigr|^p \bigr) \Biggr\rrvert
\nonumber\\
&&\qquad = \Biggl\llvert \sum_{ l = 0 }^{ \ld( N ) }
\frac{ 2^l }{ N } \sum_{ k = 1 }^{{ N }/{ 2^l } } \bigl
\llvert y_{ 2^l }^{ 2^l, \xi^{l,k}(\omega) } \bigr\rrvert^{ p } - \sum
_{ l = 1 }^{ \ld( N ) } \frac{ 2^{l} }{ N } \sum
_{ k = 1 }^{
{ N }/{ 2^{l} }
} \bigl\llvert y_{ 2^{ (l-1) } }^{
2^{ (l-1) }, \xi^{l,k}(\omega)
}
\bigr\rrvert^{ p } \Biggr\rrvert \\
&&\qquad\geq r(N) \cdot T^{ - { p }/{ 4 } }\nonumber
\end{eqnarray}
for all
$ N \in\{ N_1(\omega),
2^1 \cdot N_1(\omega), 2^2
\cdot N_1(\omega), \ldots\}
\cap
[
N_2
,
\infty
)
$ and all
$ \omega
\in
\{ N_1 < \infty \}
$.
The fact
$ \lim_{ N \rightarrow\infty}
r(N) = \infty$
therefore shows
%
%e76 #&#
\begin{eqnarray}
&&\mathop{\lim_{
N \rightarrow\infty}}_{
\operatorname{ld}(N) \in\mathbb{N}
} \Biggl\llvert \frac{1}{N}
\sum_{ k = 1 }^{ N } \bigl| Y_1^{ 1, 0, k }(
\omega) \bigr|^p\nonumber
\\[-8pt]
\\[-8pt]
\nonumber
&&\quad\qquad{} + \sum_{ l = 1 }^{ \ld( N ) }
\frac{ 2^l }{ N } \sum_{ k = 1 }^{{ N }/{ 2^l } } \bigl( \bigl|
Y_{ 2^l }^{ 2^l, l, k }(\omega) \bigr|^p
 - \bigl| Y_{ 2^{ (l-1) } }^{
2^{ (l-1) }, l, k
}(
\omega) \bigr|^p \bigr) \Biggr\rrvert = \infty
\end{eqnarray}
for all
$ \omega
\in
\{ N_1 < \infty \}
$.
Hence, Lemma~\ref{lem_n1} finally yields
%
%e77 #&#
\begin{eqnarray}
\mathop{\lim_{
N \rightarrow\infty}}_{
\operatorname{ld}(N) \in\mathbb{N}
} \Biggl\llvert \frac{1}{N}
\sum_{ k = 1 }^{ N } \bigl| Y_1^{ 1, 0, k }
\bigr|^p + \sum_{ l = 1 }^{ \ld( N ) }
\frac{ 2^l }{ N } \sum_{ k = 1 }^{{ N }/{ 2^l } } \bigl( \bigl|
Y_{ 2^l }^{ 2^l, l, k } \bigr|^p - \bigl| Y_{ 2^{ (l-1) } }^{
2^{ (l-1) }, l, k
}
\bigr|^p \bigr) \Biggr\rrvert = \infty
\end{eqnarray}
$ \mathbb{P} $-almost surely.
This completes
the proof of
Theorem~\ref{thmconj}.
\end{pf*}

%s5 #&#
\section{Divergence of the multilevel Monte Carlo Euler method}
\label{secdivergencemlMCE}
\begin{figure}

\includegraphics{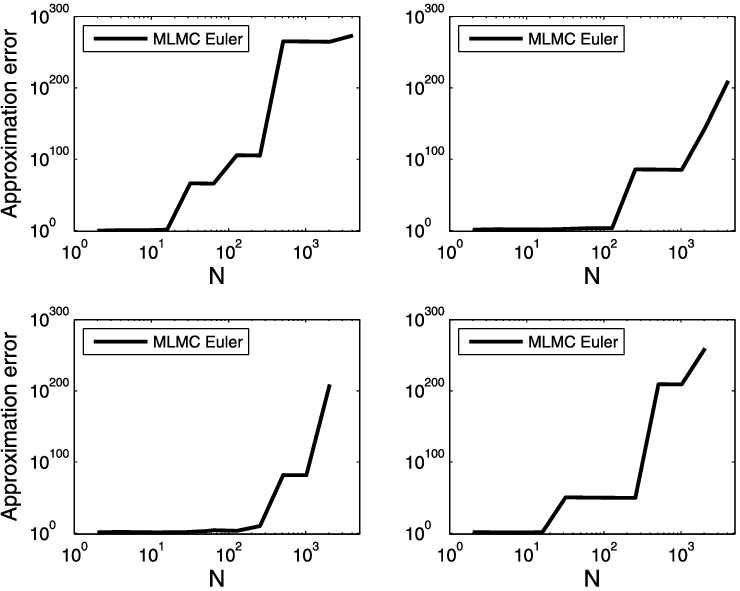}

\caption{Four sample paths of the approximation error of the multilevel
Monte Carlo Euler approximation
for the Ginzburg--Landau equation~\protect\eqref{eqginzburg}.}
\label{fginzburg}
\end{figure}

Motivated by Figure~\ref{fginzburg} below and
by the divergence result of the
multilevel Monte Carlo Euler
method in Section~\ref{secmultilevelpathwisedivergence},
we conjecture in this section that
the multilevel Monte Carlo Euler method
diverges with probability
one whenever
one of the coefficients of the SDE
grows superlinearly; see Conjecture~\ref{conj}.
Whereas divergence with probability
one seems to be quite difficult
to establish,
strong divergence is a rather
immediate consequence
of the divergence
of the Euler method
in Theorem~\ref{thmeEdivergence}
above.
We derive this strong divergence in Corollary~\ref{thmlimmean} below.
For practical simulations the
much more important
question is, however, consistency
and inconsistency, respectively;
see, for example,
Nikulin~\cite{n01},
Cram{\'e}r~\cite{c99},
Appendix~A.1 in
Glasserman~\cite{g04}
and also
Theorem~\ref{thmconj} above and
Conjecture~\ref{conj} below.

Throughout this section assume that
the following setting is fulfilled.
Let $ T \in( 0, \infty) $,
let $ ( \Omega, \mathcal{F},
\mathbb{P} ) $ be a probability space
with a normal filtration
$ ( \mathcal{F}_t )_{ t \in[0,T] } $,
let $ W^{ l, k } \dvtx[0,T] \times
\Omega\rightarrow\mathbb{R} $,
$ l \in\N_0 $,
$ k \in\N$, be a family
of independent
one-dimensional standard
$ ( \mathcal{F}_t )_{ t \in[0,T] }
$-Brownian motions and let
$ \xi^{ l, k } \dvtx\Omega
\rightarrow\mathbb{R} $,
$ l \in\N_0 $,
$ k \in\N$, be a
family
of independent identically
distributed
$ \mathcal{F}_0 /
\mathcal{B}( \mathbb{R} )
$-measurable mappings
with
$
\mathbb{E} [ |\xi^{ 0, 1 } |^p  ]
< \infty
$
for all $ p \in[1, \infty) $.
Moreover, let $ \mu, \sigma\dvtx\mathbb{R}
\rightarrow\mathbb{R} $ be two
continuous mappings such that
there exists
a predictable stochastic process
$ X \dvtx[0,T] \times\Omega\rightarrow\mathbb{R} $
which satisfies
$
\int_0^T
\llvert  \mu( X_s ) \rrvert
+
\llvert  \sigma( X_s ) \rrvert^2 \,ds
<
\infty
$
$ \mathbb{P} $-almost surely\ and
%
%e78 #&#
\begin{equation}
\label{eqsolution1} X_t = \xi^{ 0, 1 } + \int
_0^t \mu( X_s ) \,ds + \int
_0^t \sigma( X_s )
\,dW_s^{ 0, 1 }
\end{equation}
$ \mathbb{P} $-almost surely
for all $ t \in[0,T] $. The drift coefficient $ \mu$ is the
infinitesimal mean of the process $ X $ and the diffusion coefficient $
\sigma$ is the infinitesimal standard deviation of the process $ X $.
We then define a family of Euler approximations $ Y_n^{ N, l, k }
\dvtx\Omega\rightarrow\mathbb{R} $, $ n \in\{ 0, 1, \ldots, N \}
$, $ N \in\N$,
$ l \in\N_0 $,
$ k \in\N$, by
$ Y_0^{ N, l, k } := \xi^{ l, k } $ and
\[
Y_{ n + 1 }^{ N, l, k } := Y_{ n }^{ N, l, k } + \mu
\bigl( Y_{ n }^{ N, l, k } \bigr) \cdot \frac{ T }{ N } + \sigma
\bigl( Y_{ n }^{ N, l, k } \bigr) \cdot \bigl( W_{ { ( n + 1 ) T }/{ N } }^{ l, k }
- W_{ { n T }/{ N } }^{ l, k } \bigr)
\]
for all $ n \in\{ 0, 1, \ldots, N - 1 \} $,
$ N \in\N$,
$ l \in\N_0 $
and all $ k \in\N$.
%
%co5.1 #&#
\begin{conj}[(Divergence
with probability one
of the multilevel Monte Carlo
Euler method)]
\label{conj}
Assume that the above setting
is fulfilled
and let
$ \alpha, c \in( 1, \infty) $
be real numbers such that\vadjust{\goodbreak}
$
\frac{ \llvert  x \rrvert^{ \alpha} }{ c }
\leq
\llvert  \mu( x ) \rrvert
+
\llvert  \sigma( x ) \rrvert
\leq
c \llvert  x \rrvert^c
$
for all $ x \in\mathbb{R} $
with $ | x | \geq c $.
Moreover, assume that
$ \mathbb{P} [ \sigma( \xi^{0,1} ) \neq0  ] > 0 $
or that
there exists a real number
$\beta\in(1,\infty)$ such that
$
\P [
\llvert \xi^{0,1}\rrvert \geq x
]
\geq\beta^{  ( - x^\beta ) }$
for all $x\in[1,\infty)$.
Moreover,
let $ f \dvtx\mathbb{R} \rightarrow\mathbb{R} $
be $ \mathcal{B}( \mathbb{R} ) / \mathcal{B}( \mathbb{R} )
$-measurable
with
$
\frac{ 1 }{ c } | x |^{ { 1 }/{ c } }
- c \leq f( x )
\leq
c ( 1 + | x |^c )
$
for all $ x \in\mathbb{R} $.
Then we conjecture
%
%e79 #&#
\begin{eqnarray}
\label{eqlimmean2} && \mathop{\lim_{
N \rightarrow\infty}}_{ \ld( N ) \in\N
} \Biggl\llvert
\mathbb{E} \bigl[ f( X_T ) \bigr]\nonumber\hspace*{-30pt}
\\[-9pt]
\\[-9pt]
\nonumber
&&\qquad{}  - \frac{ 1 }{ N } \sum
_{ k = 1 }^{ N } f \bigl( Y_1^{ 1, 0, k }
\bigr)
- \sum_{ l = 1 }^{ \ld( N ) } \frac{ 2^l }{ N }
\Biggl( \sum_{ k = 1 }^{{ N }/{ 2^l } } f \bigl(
Y_{ 2^l }^{ 2^l, l, k } \bigr) - f \bigl( Y_{ 2^{ ( l - 1 ) } }^{ 2^{ ( l - 1 ) }, l, k }
\bigr) \Biggr) \Biggr\rrvert = \infty\hspace*{-30pt}
\end{eqnarray}
$\mathbb{P}$-almost surely.
\end{conj}
To support this conjecture, we ran simulations
for the stochastic Ginzburg--Landau equation
given by
the solution $(X_t)_{t\in[0,1]}$ of
%
%e80 #&#
\begin{equation}
\label{eqginzburg} d X_t = \bigl( 2X_t -
X_t^3 \bigr) \,dt + 2 X_t
\,dW_t,\qquad  X_0 = 1
\end{equation}
for all $ t \in[0,1] $.
Its solution is known
explicitly (e.g., Section 4.4 in~\cite{kp92})
and is given by
%
%e81 #&#
\begin{equation}
\label{eqginzburgexplicit} X_t=\frac{\exp (2W_t )} {
\sqrt{1+2 \int_0^t\exp (4 W_s )\,ds}}
\end{equation}
for $t \in[0,1]$.
We used this explicit solution to compute
$\E[ (X_1)^2 ] \approx0.8114 $.
%More precisely,
%we approximated the Lebesgue integral in the denominator
%of~\eqref{eqginzburgexplicit} with a Riemann sum
%with $10^4$
%summands.
%Moreover, we approximate the second moment at time $1$
%by a Monte Carlo simulation with $10^8$ independent approximations of
%$X_1$.
%This results in the approximate value
%$\mathbb{E}\bigl[(X_1)^2\bigr]\approx TODO$.
Figure~\ref{fginzburg} shows
four sample paths
of the approximation error of the multilevel Monte Carlo Euler method
for the Ginzburg--Landau equation~\eqref{eqginzburg}.
Only finite values of the sample paths are plotted.
The next corollary is an immediate
consequence of Theorem~\ref{thmeEdivergence}
above.
%%%%%%%%%%%%%%%%%%%%
%
%f4 #&#

%
%%%%%%%%%%%%%%%%%%%%

%%%%%%%%%%%%%%%%%%%%%%%%%%%%%%%%%
%
%co5.2 #&#
\begin{cor}[(Strong divergence
of the multilevel Monte Carlo
Euler method)]
\label{thmlimmean}
Assume that the above setting
is fulfilled
and let
$ \alpha, c \in( 1, \infty) $
be real numbers such that
$
\frac{ \llvert  x \rrvert^{ \alpha} }{ c }
\leq
\llvert  \mu( x ) \rrvert
+
\llvert  \sigma( x ) \rrvert
\leq
c \llvert  x \rrvert^c
$
for all $ x \in\mathbb{R} $
with $ | x | \geq c $.
Moreover, assume that
$ \mathbb{P} [ \sigma( \xi^{0,1} ) \neq0  ] > 0 $
or that
there exists a
real number $\beta\in(1,\infty)$
such that
$
\P [
|\xi^{0,1}| \geq
x
]
\geq\beta^{ (-x^\beta )}$
for all $x\in[1,\infty)$.
Additionally, let $ f \dvtx\mathbb{R} \rightarrow\mathbb{R} $
be $ \mathcal{B}( \mathbb{R} ) / \mathcal{B}( \mathbb{R} ) $-measurable
with
$
\frac{ 1 }{ c } | x |^{ { 1 }/{ c } }
- c \leq f( x )
\leq
c ( 1 + | x |^c )
$
for all $ x \in\mathbb{R} $.
Then we obtain
%
%e82 #&#
\begin{eqnarray}
\label{eqlimmeanexpectation}&& \mathop{\lim_{
N \rightarrow\infty}}_{ \ld( N ) \in\mathbb{N}
}
\mathbb{E} \Biggl[ \Biggl\llvert \mathbb{E} \bigl[ f( X_T ) \bigr] -
\frac{ 1 }{ N } \sum_{ k = 1 }^{ N } f \bigl(
Y_1^{ 1, 0, k } \bigr)
\nonumber
\\[-9pt]
\\[-9pt]
\nonumber
&&\hspace*{28pt}\qquad{}- \sum_{ l = 1 }^{ \ld( N ) }
\frac{ 2^l }{ N } \Biggl( \sum_{ k = 1 }^{{ N }/{ 2^l } } f
\bigl( Y_{ 2^l }^{ 2^l, l, k } \bigr) - f \bigl( Y_{ 2^{ ( l - 1 ) } }^{ 2^{ ( l - 1 ) }, l, k }
\bigr) \Biggr) \Biggr\rrvert^p \Biggr] = \infty
\end{eqnarray}
for all $ p \in[ 1, \infty) $.\vadjust{\goodbreak}
\end{cor}
\begin{pf}%[(Proof of Corollary~\ref{thmlimmean})]
First of all, note that the
assumption
$
\mathbb{E} [ |\xi^{ 0, 1 } |^p  ]
< \infty
$
for all $ p \in[ 1, \infty) $,
the continuity of
$
\mu, \sigma\dvtx\mathbb{R}
\rightarrow\mathbb{R}
$,
the inequality
$
\llvert  \mu( x ) \rrvert
+
\llvert  \sigma( x ) \rrvert
\leq
c \llvert  x \rrvert^c
$
for all $ x \in\mathbb{R} $ with $ | x | \geq c $
and the estimate
$ | f( x ) | \leq c( 1 + | x |^c ) $
for all $ x \in\mathbb{R} $ imply
$
\mathbb{E} [
| f( Y_N^{ N, 0, 1 } ) |
]
<
\infty
$
for all $ N \in\mathbb{N} $.
Therefore, we obtain
\begin{eqnarray*}
\mathbb{E} \Biggl[ \frac{ 1 }{ N } \sum_{ k = 1 }^N
f \bigl( Y_1^{ 1, 0, k } \bigr) + \sum
_{ l = 1 }^{ \ld( N ) } \frac{ 2^l }{ N } \Biggl( \sum
_{ k = 1 }^{{ N }/{ 2^l } } f \bigl( Y_{ 2^l }^{ 2^l, l, k }
\bigr) - f \bigl( Y_{ 2^{ ( l - 1 ) } }^{ 2^{ ( l - 1 ) }, l, k } \bigr) \Biggr) \Biggr] =
\mathbb{E} \bigl[ f \bigl( Y_N^{ N, 0, 1 } \bigr) \bigr]
\end{eqnarray*}
for all $ N \in\{ 2^1, 2^2, 2^3, \ldots\} $.
The estimate
$
f( x )
\geq
\frac{ 1 }{ c }
| x |^{ { 1 }/{ c } } - c
$
for all $ x \in\mathbb{R} $ and
Theorem~\ref{thmeEdivergence}
hence give
\begin{eqnarray*}
&&\mathop{\lim_{
N \rightarrow\infty}}_{\ld( N ) \in\mathbb{N}
} \mathbb{E} \Biggl[
\frac{ 1 }{ N } \sum_{ k = 1 }^N f \bigl(
Y_1^{ 1, 0, k } \bigr) + \sum_{ l = 1 }^{ \ld( N ) }
\frac{ 2^l }{ N } \Biggl( \sum_{ k = 1 }^{{ N }/{ 2^l } } f
\bigl( Y_{ 2^l }^{ 2^l, l, k } \bigr) - f \bigl( Y_{ 2^{ ( l - 1 ) } }^{ 2^{ ( l - 1 ) }, l, k }
\bigr) \Biggr) \Biggr]\\
&&\qquad \geq \frac{ 1 }{ c } \Bigl( \lim_{ N \rightarrow\infty}
\mathbb{E} \bigl[ \bigl\llvert Y_N^{ N, 0, 1 } \bigr
\rrvert^{ { 1 }/{ c } } \bigr] \Bigr) - c = \infty.
\end{eqnarray*}
In the case
$ \mathbb{E} [ | f( X_T ) |  ]
< \infty$,
the triangle inequality
and Jensen's inequality
then yield
%
%e83 #&#
\begin{eqnarray}
\label{eqfnotinf}
&&\mathop{\lim_{
N \rightarrow\infty}}_{ \ld( N ) \in\mathbb{N}
} \Biggl\llVert
\mathbb{E} \bigl[ f( X_T ) \bigr]\nonumber\hspace*{-35pt}\\
&&\hspace*{13pt}\qquad{} - \frac{ 1 }{ N } \sum
_{ k = 1 }^{ N } f \bigl( Y_1^{ 1, 0, k }
\bigr) - \sum_{ l = 1 }^{ \ld( N ) } \frac{ 2^l }{ N }
\Biggl( \sum_{ k = 1 }^{{ N }/{ 2^l } } f \bigl(
Y_{ 2^l }^{ 2^l, l, k } \bigr) - f \bigl( Y_{ 2^{ ( l - 1 ) } }^{ 2^{ ( l - 1 ) }, l, k }
\bigr) \Biggr) \Biggr\rrVert_{ L^p( \Omega; \mathbb{R} ) }
\nonumber\hspace*{-35pt}\\
&&\quad\geq \mathop{\lim_{
N \rightarrow\infty}}_{ \ld( N ) \in\mathbb{N}
} \mathbb{E} \Biggl[
\frac{ 1 }{ N } \sum_{ k = 1 }^{ N } f \bigl(
Y_1^{ 1, 0, k } \bigr) - \sum_{ l = 1 }^{ \ld( N ) }
\frac{ 2^l }{ N } \Biggl( \sum_{ k = 1 }^{{ N }/{ 2^l } } f
\bigl( Y_{ 2^l }^{ 2^l, l, k } \bigr) - f \bigl( Y_{ 2^{ ( l - 1 ) } }^{ 2^{ ( l - 1 ) }, l, k }
\bigr) \Biggr) \Biggr]\hspace*{-35pt}\\
&&\qquad{} - \mathbb{E} \bigl[ \bigl| f( X_T )\bigr | \bigr]\nonumber\hspace*{-35pt}\\
&&\quad =
\infty\nonumber\hspace*{-35pt}
\end{eqnarray}
for all $ p \in[ 1, \infty) $.
This shows \eqref{eqlimmeanexpectation}
in the case $ \mathbb{E} [ | f( X_T ) |  ] < \infty
$.
In the case
$ \mathbb{E} [ | f( X_T ) |  ] = \infty
$, the estimate
$ f(x) \geq-c $ for all $ x \in\mathbb{R} $
shows
$ \mathbb{E} [ f( X_T )  ] = \infty
$
and this implies
\eqref{eqlimmeanexpectation}
in the case
$ \mathbb{E} [ | f( X_T ) |  ] = \infty
$.
The proof of
Corollary~\ref{thmlimmean}
is thus completed.
\end{pf}
%
% %%%%
% %%%%%%%%%%%%%%%%%%%%%%%%%%%%%%%%%%%%%%%%%%%%%%%%%%%
%
% %%%%
% %%%%%%%%%%%%%%%%%%%%%%%%%%%%%%%%%%%%%%%%%%%%%%%%%%%
%
%s6 #&#
\section{Convergence of the multilevel
Monte Carlo tamed Euler method}
\label{sectamedconvergence}
In this section we combine the multilevel Monte Carlo method
with a tamed Euler method.
We aim at path-dependent payoff functions.
Therefore, we consider
piecewise linear time
interpolations of the numerical
approximations, which have
continuous sample paths
and which are implementable.
Theorem~\ref{thmtamedconvergence2} shows that
these piecewise linear interpolations
of the tamed Euler approximations
converge in the strong sense
with the optimal convergence order
according to M\"{u}ller-Gronbach's
lower bound in the Lipschitz
case in~\cite{m02}.
Theorem~\ref{thmtamedconvergence}
then establishes almost sure
and strong convergence of the
multilevel Monte Carlo method
combined with the tamed Euler method.
The payoff function is allowed to depend on
the whole path.
We assume the payoff function only
to be locally Lipschitz continuous
and the local Lipschitz constant to grow at most polynomially.

Throughout this section assume
that the following setting is fulfilled.
Let $ T \in(0, \infty) $,
let $ ( \Omega, \mathcal{F}, \mathbb{P} ) $
be a probability space with a
normal filtration
$ ( \mathcal{F}_t )_{ t \in[0,T] } $,
let $ d, m \in\N$, let
$ W^{ l, k } \dvtx[0,T] \times\Omega
\rightarrow\mathbb{R}^m $,
$ l \in\N_0 $,
$ k \in\N$,
be a family of independent standard
$ ( \mathcal{F}_t )_{ t \in[0,T] } $-Brownian
motions and let
$ \xi^{ l, k } \dvtx\Omega\rightarrow\mathbb{R}^d $,
$ l \in\N_0 $,
$ k \in\N$,
be a family of independent identically
distributed $ \mathcal{F}_0
/ \mathcal{B}( \mathbb{R}^d )
$-measurable mappings with
$ \mathbb{E} [ \| \xi^{ 0, 1 }
\|^p_{ \mathbb{R}^d }
] < \infty$ for all $ p \in[1, \infty) $.
Here and below we use the Euclidean norm
$\llVert x\rrVert_{\mathbb{R}^n}:=\sqrt{x_1^2+x_2^2+\cdots+x_n^2}$
for all $x=(x_1,x_2,\ldots,x_n)\in\mathbb{R}^n$ and all $n\in\N$.
Moreover, let
$ \mu\dvtx\mathbb{R}^d
\rightarrow\mathbb{R}^d $
be a continuously differentiable
and globally one-sided Lipschitz continuous function whose derivative
grows at most polynomially and
let $ \sigma=\break ( \sigma_{ i, j } )_{ i \in\{ 1, 2, \ldots, d \}, j
\in
\{ 1, 2, \ldots, m \} }
\dvtx\mathbb{R}^d \rightarrow
\mathbb{R}^{ d \times m } $ be a globally Lipschitz continuous
function. More formally, suppose that
there exists a real number $ c \in[0, \infty) $
such that
$
\anglel  x - y, \mu( x ) - \mu( y )
\angler_{ \mathbb{R}^d }
\leq c \| x - y \|^2_{ \mathbb{R}^d }
$,
$
\| \mu'( x ) \|_{ L( \mathbb{R}^d ) }
\leq c ( 1 + \| x \|_{ \mathbb{R}^d }^c )
$
and
$
\| \sigma( x ) - \sigma( y )
\|_{ L( \mathbb{R}^m, \mathbb{R}^d ) }
\leq
c
\| x - y \|_{ \mathbb{R}^d }
$
for all $ x, y \in\mathbb{R}^d $.
Here and below we use
$ \| x \|:=
( \sum_{i=1}^d | x_i |^2
)^{ {1}/{2} }
$
and
$\anglel x,y\angler_{\R^d}:=
\sum_{i=1}^d x_i\cdot y_i$
for all $x=(x_1,x_2,\ldots,x_d),y=(y_1,\ldots,y_d)\in\mathbb{R}^d$.
Then consider the SDE
%
%e84 #&#
\begin{equation}
\label{eqSDE} dX_t = \mu( X_t ) \,dt + \sigma(
X_t ) \,dW_t^{ 0, 1 }, \qquad X_0 =
\xi^{ 0, 1 }
\end{equation}
for $ t \in[0,T] $.
Under the assumptions above, the SDE~\eqref{eqSDE} is known to have a
unique solution. More formally, there exists an
up to indistinguishability
unique adapted stochastic
process
$ X \dvtx[0,T] \times\Omega\rightarrow
\mathbb{R}^d $ with continuous
sample paths fulfilling
%
%e85 #&#
\begin{equation}
X_t = \xi^{ 0, 1 } + \int_0^t
\mu( X_s ) \,ds + \int_0^t
\sigma( X_s ) \,dW_s^{ 0, 1 }
\end{equation}
$ \mathbb{P} $-almost surely for all $ t \in[0,T] $; see, for example,
Theorem~2.4.1 in Mao~\cite{m97}. The drift
coefficient $ \mu$ is the
infinitesimal mean of the process
$ X $ and the diffusion
coefficient $ \sigma$ is
the infinitesimal standard
deviation of the process $ X $.
In the next step we define
a family of tamed Euler approximations
$ Y_n^{ N, l, k } \dvtx\Omega
\rightarrow\mathbb{R}^d $,
$ n \in\{ 0, 1, \ldots, N \} $,
$ N \in\N$,
$ l \in\N_0 $,
$ k \in\N$,
by
$ Y_0^{ N, l, k } := \xi^{ l, k } $
and
%
%e86 #&#
\begin{eqnarray}
\label{eqtamedEuler} Y_{ n + 1 }^{ N, l, k } &:=& Y_n^{ N, l, k }
+ \frac{
\mu( Y_n^{ N, l, k } )
\cdot
{ T }/{ N }
}{
1 +
\llVert \mu( Y_n^{ N, l, k } )
\cdot
{ T }/{ N }
\rrVert_{ \mathbb{R}^d }
}
\nonumber
\\[-8pt]
\\[-8pt]
\nonumber
&&{} + \sigma \bigl( Y_n^{ N, l, k } \bigr)
\bigl( W_{ { ( n + 1 ) T }/{ N } }^{ l, k } - W_{ { n T }/{ N } }^{ l, k } \bigr)
\end{eqnarray}
for all $ n \in\{ 0, 1, \ldots, N - 1 \} $,
$ N \in\N$,
$ l \in\N_0 $
and all $ k \in\N$.
In order to formulate our convergence
theorem for the multilevel Monte Carlo
tamed Euler approximations, we now
introduce piecewise continuous time
interpolations of the time discrete numerical approximations~\eqref
{eqtamedEuler}.
More formally, let
$ \bar{Y}^{ N, l, k } \dvtx
[0,T] \times\Omega\rightarrow
\mathbb{R}^d $,
$ N \in\N$,
$ l \in\N_0 $,
$ k \in\N$,
be a family of stochastic processes
with continuous sample paths defined by
%
%e87 #&#
\begin{eqnarray}
%
% \label{eqcontinuoustimeinterpolation}
\bar{Y}_t^{ N, l, k } &:=&
Y_n^{ N, l, k } + \frac{
( t - { n T }/{ N }  )
}{
{T}/{N}
} \bigl( Y_{ n + 1 }^{ N, l, k }
- Y_n^{ N, l, k } \bigr)
\nonumber
\\[-8pt]
\\[-8pt]
\nonumber
&=& \biggl( \frac{ t N }{ T } - n \biggr)
Y_{ n + 1 }^{ N, l, k } + \biggl( n + 1 - \frac{ t N }{ T } \biggr)
Y_n^{ N, l, k } %
\end{eqnarray}
for all $ t \in[ \frac{ n T }{ N },
\frac{ ( n + 1 ) T }{ N } ] $,
$ n \in\{ 0, 1, \ldots, N - 1 \} $,
$ N \in\N$,
$ l \in\N_0 $
and all $ k \in\N$.

The following corollary is a direct
consequence of
Hutzenthaler, Jentzen and Kloeden~\cite{hjk10b}
and M\"{u}ller-Gron\-bach~\cite{m02}; see also Ritter~\cite{r90b}.
It asserts that the piecewise linear approximations
$\bar{Y}^N$, $N\in\N$,
converge in the strong sense to the exact solution.
The convergence order is $\frac{1}{2}$
except for a logarithmic term.
%
%%%%%%%%%%%%%%%%%%%%%%%%%%%%%%%%%%%%%%%%%%5
% Theorem
%%%%%%%%%%%%%%%%%%%%%%%%%%%%%%%%%%%%%%%%%%5
%
%co6.1 #&#
\begin{cor}[(Strong convergence
of the tamed Euler method)]
\label{thmtamedconvergence2}
Assume that the above setting is fulfilled. Then there exists a family
$ R_p \in[0,\infty) $, $ p \in[1, \infty) $, of real numbers such that
%
%e88 #&#
\begin{equation}
\label{eqtamedconvergence3} \Bigl( \mathbb{E} \Bigl[ \sup_{ t \in[0,T] } \bigl
\llVert X_t - \bar{Y}_t^{ N, 0, 1 } \bigr
\rrVert^p_{
\mathbb{R}^d
} \Bigr] \Bigr)^{ { 1 }/{ p } } \leq
R_p \cdot \frac{ \sqrt{ 1 + \ld( N ) } }{ \sqrt{ N } }
\end{equation}
for all $ N \in\N$ and all $ p \in[ 1, \infty) $.
\end{cor}
The convergence rate
$
N^{ - {1}/{2} }
( 1 + \ld( N ) )^{ {1}/{2} }
$
for $ N \in\N$
obtained in \eqref{eqtamedconvergence3}
is sharp according to
M\"{u}ller-Gronbach's
lower bound established
in Theorem~3
in~\cite{m02} in the
case of globally Lipschitz
continuous coefficients; see also Hofmann,
M\"{u}ller-Gronbach
and Ritter~\cite{hmr00a}.

\begin{pf*}{Proof of Corollary~\ref{thmtamedconvergence2}}
Let
$ \tilde{Y}^N \dvtx[0,T] \times\Omega
\rightarrow\mathbb{R}^d $,
$ N \in\N$,
be
stochastic processes
defined by
\[
\tilde{Y}_t^N := Y_n^{ N, 0, 1 } +
\frac{
\mu( Y_n^{ N, 0, 1 } )
\cdot
( t - { n T }/{ N }  )
}{
1 +
\|
\mu( Y_n^{ N, 0, 1 } )
\cdot
{ T }/{ N }
\|_{ \mathbb{R}^d }
} + \sigma \bigl( Y_n^{ N, 0, 1 } \bigr) \bigl(
W_t^{0,1} - W_{ { n T }/{ N } }^{0,1} \bigr)
\]
for all $ t \in[ \frac{ n T }{ N } $,
$ \frac{ ( n + 1 ) T }{ N } ] $,
$ n \in\{ 0, 1, \ldots, N - 1 \} $
and all $ N \in\N$.
Theorem~1.1 in~\cite{hjk10b} then
shows the existence of a family
$ \tilde{R}_p \in[0, \infty) $,
$ p \in[1, \infty) $, of real numbers such that
$
\|
\sup_{ t \in[0,T] }
\|
X_t - \tilde{Y}_t^N
\|_{ \mathbb{R}^d }
\|_{ L^p( \Omega; \mathbb{R} ) }
\leq
\frac{ \tilde{R}_p }{ \sqrt{N} }
$
for all $ N \in\N$
and all $ p\in[1, \infty) $.
The triangle inequality hence yields
%
%e89 #&#
\begin{eqnarray}
\label{eqest1}&&  \Bigl\llVert \sup_{ t \in[0,T] } \bigl\llVert X_t
- \bar{Y}_t^{ N, 0, 1 } \bigr\rrVert_{ \mathbb{R}^d } \Bigr
\rrVert_{ L^p( \Omega; \mathbb{R} ) }
\nonumber
\\[-8pt]
\\[-8pt]
\nonumber
&&\qquad\leq \frac{ \tilde{R}_p }{ \sqrt{N} } + \Bigl\llVert
\sup_{ t \in[0,T] } \bigl\llVert \tilde{Y}_t^N -
\bar{Y}_t^{ N, 0, 1 } \bigr\rrVert_{ \mathbb{R}^d } \Bigr
\rrVert_{ L^p( \Omega; \mathbb{R} ) }
\end{eqnarray}
for all $ N \in\N$
and all $ p \in[1, \infty) $.
Moreover, we have
%
%e90 #&#
\begin{eqnarray}
\label{eqest2}&& \bigl\llVert \tilde{Y}_t^N -
\bar{Y}_t^{ N, 0, 1 } \bigr\rrVert_{ \mathbb{R}^d }\nonumber\\
 &&\qquad= \biggl
\llVert \sigma \bigl( Y_n^{ N, 0, 1 } \bigr) \bigl(
W_t^{0,1} - W_{ { n T }/{ N } }^{0,1} \bigr)\nonumber\\
&&\hspace*{15pt}\qquad {}-
\biggl( \frac{ t N }{ T } - n \biggr) \sigma \bigl( Y_n^{ N, 0, 1 }
\bigr) \bigl( W_{ { ( n + 1 ) T }/{ N } }^{0,1} - W_{ { n T }/{ N } }^{0,1}
\bigr) \biggr\rrVert_{ \mathbb{R}^d }
\\
&&\qquad\leq
\bigl\llVert \sigma \bigl( Y_n^{ N, 0, 1 }
\bigr) \bigr\rrVert_{
L( \mathbb{R}^m, \mathbb{R}^d )
} \biggl\llVert W_t^{0,1}
- W_{ { n T }/{ N } }^{0,1} \nonumber\\
&&\hspace*{108pt}\qquad{}- \biggl( \frac{ t N }{ T } - n \biggr)
\bigl( W_{ { ( n + 1 ) T }/{ N } }^{0,1} - W_{ { n T }/{ N } }^{0,1}
\bigr) \biggr\rrVert_{ \mathbb{R}^m }\nonumber
\end{eqnarray}
for all
$ t \in[ \frac{ n T }{ N }, \frac{ (n+1) T }{ N } ] $,
$ n \in \{ 0, 1, \ldots, N-1  \} $
and all
$ N \in\N$.
Combining \eqref{eqest1}, \eqref{eqest2}
and H{\"o}lder's inequality then gives
%
%e91 #&#
%e92 #&#
%e93 #&#
%e94 #&#
\begin{eqnarray}
\label{eqtamedconvergenceA}
&&\hspace*{-4pt} \Bigl\llVert \sup_{ t \in[0,T] } \bigl
\llVert X_t - \bar{Y}_t^{ N, 0, 1 } \bigr
\rrVert_{ \mathbb{R}^d } \Bigr\rrVert_{ L^p( \Omega; \mathbb{R} ) } \nonumber\\
&&\hspace*{-4pt}\qquad\leq \frac{ \tilde{R}_p }{ \sqrt{N} } +
\Bigl\llVert \max_{ n \in\{ 0, 1, \ldots, N - 1 \} } \bigl\llVert \sigma \bigl( Y_n^{ N, 0, 1 }
\bigr) \bigr\rrVert_{ L( \mathbb{R}^m, \mathbb{R}^d ) } \Bigr\rrVert_{ L^{ 2p }( \Omega; \mathbb{R} ) }
\nonumber\\
&&\hspace*{-4pt}\qquad\quad{} \times\biggl\llVert \max_{ n \in\{ 0, 1, \ldots, N - 1 \} }
\sup_{
t \in
[
{ n T }/{N},
{ ( n + 1 ) T }/{ N }
]
} \biggl\llVert
W_t^{0,1} - W_{ { n T }/{ N } }^{0,1}
\nonumber
\\[-8pt]
\\[-8pt]
\nonumber
&&\hspace*{-2pt}\hspace*{114pt}\qquad{}- \biggl(
\frac{ t N }{ T } - n \biggr) \bigl( W_{ { ( n + 1 ) T }/{ N } }^{0,1} -
W_{ { n T }/{ N } }^{0,1} \bigr) \biggr\rrVert_{ \mathbb{R}^m } \biggr
\rrVert_{ L^{ 2p }( \Omega; \mathbb{R} ) }
\nonumber\\
&&\hspace*{-4pt}\qquad \leq
\frac{ \tilde{R}_p }{ \sqrt{N} } + \sqrt{ \frac{ T }{ N } } \Bigl( c \cdot
\sup_{ M \in\N} \Bigl\llVert \max_{ n \in\{ 0, 1, \ldots, M \} } \bigl\llVert
Y_n^{ M, 0, 1 } \bigr\rrVert_{ \mathbb{R}^d } \Bigr
\rrVert_{ L^{ 2p }( \Omega; \mathbb{R} ) } + \bigl\llVert \sigma( 0 ) \bigr\rrVert_{ L( \mathbb{R}^m, \mathbb{R}^d ) }
\Bigr) \nonumber\\
&&\hspace*{-4pt}\qquad\quad{}\times\Bigl\llVert \max_{ n \in\{ 1, 2, \ldots, N \} } \sup_{ t \in[ 0, 1 ] } \bigl\llvert
\beta_t^n - t \cdot\beta_1^n
\bigr\rrvert \Bigr\rrVert_{ L^{ 2p }( \Omega; \mathbb{R} ) }\nonumber
\end{eqnarray}
for all $ N \in\N$ and all $ p \in[ 1, \infty) $
where $ \beta^n \dvtx[0,1] \times\Omega\rightarrow\mathbb{R} $, $
n \in\N$, is a sequence of independent
one-dimensional standard Brownian motions. Moreover, Theorem~1.1 in
\cite{hjk10b},
in particular,
implies
%
%e95 #&#
\begin{equation}
\label{eqtamedconvergenceB} \sup_{ M \in\N} \Bigl\llVert
\max_{ n \in\{ 0, 1, \ldots, M \} } \bigl\llVert Y_n^{ M, 0, 1 } \bigr
\rrVert_{ \mathbb{R}^d } \Bigr\rrVert_{ L^p( \Omega; \mathbb{R} ) } < \infty
\end{equation}
for all $ p \in[ 1, \infty) $. Additionally,
Corollary~2 in
M\"{u}ller-Gronbach~\cite{m02}
(see also Ritter~\cite{r90b})
shows
%
%e96 #&#
\begin{equation}
\label{eqtamedconvergenceC} \sup_{ N \in\N} \Bigl( \bigl( 1 + \ld( N )
\bigr)^{ -{1}/{2} } \Bigl\llVert \max_{ n \in\{ 0, 1, \ldots, N \} } \sup_{ t \in[0,1] }
\bigl\llvert \beta_t^n - t \cdot\beta_1^n
\bigr\rrvert \Bigr\rrVert_{ L^p( \Omega; \mathbb{R} ) } \Bigr) < \infty
\end{equation}
for all $ p \in[ 1, \infty) $. Combining~\eqref
{eqtamedconvergenceA}, \eqref{eqtamedconvergenceB} and \eqref
{eqtamedconvergenceC} finally completes the proof of Corollary~\ref
{thmtamedconvergence2}.
\end{pf*}
%
%pr6.2 #&#
\begin{prop}[(Strong
consistency, converence
with probability one and strong
convergence of the multilevel Monte
Carlo tamed Euler method)]
\label{thmtamedconvergence}
Assume that the above setting is fulfilled,
let $ c \in[0, \infty) $ and let
$ f \dvtx
C( [0,T],\break \mathbb{R}^d ) \rightarrow
\mathbb{R} $ be a function from the space
of continuous functions
$ C( [0,T], \mathbb{R}^d ) $
into the real numbers $ \mathbb{R} $
satisfying
%
%e97 #&#
\begin{eqnarray}
\label{eqtamedconvergence1} &&\bigl\llVert f( v ) - f( w ) \bigr
\rrVert_{ C( [0,T], \mathbb{R}^d ) }
\nonumber
\\[-8pt]
\\[-8pt]
\nonumber
&&\qquad\leq c \bigl( 1 + \llVert v \rrVert_{ C( [0,T], \mathbb{R}^d ) }^c
+ \llVert w \rrVert_{ C( [0,T], \mathbb{R}^d ) }^c \bigr) \llVert v - w
\rrVert_{ C( [0,T], \mathbb{R}^d ) }
\end{eqnarray}
for all $ v, w \in C( [0,T], \mathbb{R}^d ) $. Then there exists a
family $ C_p \in[0, \infty) $, $ p \in[1, \infty) $, of real numbers
such that
%
%e98 #&#
\begin{eqnarray}
\label{eqtamedconvergence2} &&\Biggl( \mathbb{E} \Biggl[ \Biggl\llvert \mathbb{E}
\bigl[ f( X ) \bigr] \nonumber \\
&&\qquad{}- \frac{ 1 }{ N } \sum_{ k = 1 }^{ N }
f \bigl( \bar{Y}^{ 1, 0, k } \bigr)- \sum_{ l = 1 }^{ \ld( N ) }
\frac{ 2^l }{ N } \Biggl( \sum_{ k = 1 }^{{ N }/{ 2^l } } f
\bigl( \bar{Y}^{ 2^l, l, k } \bigr) - f \bigl( \bar{Y}^{ 2^{ ( l - 1 ) }, l, k } \bigr)
\Biggr) \Biggr\rrvert^p \Biggr] \Biggr)^{ { 1 }/{ p } }
\\
&&\qquad\leq
C_p \cdot \frac{
( 1 + \ld( N )  )^{ { 3 }/{ 2 } }
}{ \sqrt{ N } }\nonumber
\end{eqnarray}
for all $ N \in\{ 2^1, 2^2, 2^3, \ldots\} $
and all $ p \in[1, \infty) $.
In particular, there are finite
$ \mathcal{F} / \mathcal{B}( [0,\infty) )
$-measurable mappings
$ \tilde{C}_{ \varepsilon} \dvtx\Omega
\rightarrow[0,\infty) $,
$ \varepsilon\in(0,\frac{1}{2}) $,
such that
%
%e99 #&#
\begin{eqnarray}
\label{eqtamedconvergence4}\qquad && \Biggl\llvert \mathbb{E} \bigl[ f( X ) \bigr] -
\frac{ 1 }{ N } \sum_{ k = 1 }^{ N } f \bigl(
\bar{Y}^{ 1, 0, k } \bigr) - \sum_{ l = 1 }^{ \ld( N ) }
\frac{ 2^l }{ N } \Biggl( \sum_{ k = 1 }^{{ N }/{ 2^l } } f
\bigl( \bar{Y}^{ 2^l, l, k } \bigr) - f \bigl( \bar{Y}^{ 2^{ ( l - 1 ) }, l, k } \bigr)
\Biggr) \Biggr\rrvert
\nonumber
\\[-8pt]
\\[-8pt]
\nonumber
&&\qquad \leq \frac{ \tilde{C}_{ \varepsilon} }{
N^{
(
{1}/{2} - \varepsilon
)
}
}
\end{eqnarray}
for all $ N \in\N$
and all $ \varepsilon\in(0,\frac{1}{2}) $
$ \mathbb{P} $-almost surely.
\end{prop}
The convergence rate
$
N^{ - {1}/{2} }
(
1 + \ld(N)
)^{ {3}/{2} }
$
for $ N \in\N$
obtained in
\eqref{eqtamedconvergence2}
is the same as in
Remark~8
in
Creutzig et al.~\cite{cdmr09}.
For numerical approximation results
for SDEs with globally Lipschitz
continuous coefficients but
under less restrictive smoothness
assumption on the payoff function,
the reader is referred to
Giles, Higham and
Mao~\cite{ghm09}
and D\"{o}rsek
and Teichmann~\cite{dt10}.
Moreover, numerical approximation
results for SDEs with
nonglobally Lipschitz continuous
and at most linearly growing
coefficients can be found
in Yan~\cite{y02}, for instance.
\begin{pf*}{Proof of Proposition~\ref{thmtamedconvergence}}
The triangle inequality gives
\begin{eqnarray*}
&&\hspace*{-4pt} \Biggl\llVert \mathbb{E} \bigl[ f( X ) \bigr] - \frac{ 1 }{ N } \sum
_{ k = 1 }^{ N } f \bigl( \bar{Y}^{ 1, 0, k }
\bigr) - \sum_{ l = 1 }^{ \ld( N ) } \frac{ 2^l }{ N }
\Biggl( \sum_{ k = 1 }^{{ N }/{ 2^l } } f \bigl(
\bar{Y}^{ 2^l, l, k } \bigr) - f \bigl( \bar{Y}^{ 2^{ ( l - 1 ) }, l, k } \bigr) \Biggr)
\Biggr\rrVert_{ L^p( \Omega; \mathbb{R} ) }
\\
&&\hspace*{-4pt}\qquad\leq \bigl\llvert \mathbb{E} \bigl[ f( X ) \bigr] - \mathbb{E} \bigl[ f \bigl(
\bar{Y}^{ N, 0, 1 } \bigr) \bigr] \bigr\rrvert + \frac{ 1 }{ N } \Biggl
\llVert \sum_{ k = 1 }^{ N } \bigl( \mathbb{E}
\bigl[ f \bigl( \bar{Y}^{ 1, 0, 1 } \bigr) \bigr] - f \bigl(
\bar{Y}^{ 1, 0, k } \bigr) \bigr) \Biggr\rrVert_{ L^p( \Omega; \mathbb{R} ) }
\\
&&\hspace*{-4pt}\qquad\quad{} + \sum_{ l = 1 }^{ \ld( N ) } \frac{ 2^l }{ N }
\Biggl\llVert \sum_{ k = 1 }^{{ N }/{ 2^l } } \bigl(
\mathbb{E} \bigl[ f \bigl( \bar{Y}^{ 2^l, 0, 1 } \bigr) \bigr] - \mathbb{E}
\bigl[ f \bigl( \bar{Y}^{ 2^{ ( l - 1 ) }, 0, 1 } \bigr) \bigr] \\
&&\hspace*{96pt}\qquad{}- f \bigl(
\bar{Y}^{ 2^l, l, k } \bigr) + f \bigl( \bar{Y}^{ 2^{ ( l - 1 ) }, l, k } \bigr) \bigr)
\Biggr\rrVert_{ L^p( \Omega; \mathbb{R} ) }
\end{eqnarray*}
for all $ N \in\{ 2^1, 2^2, 2^3, \ldots\}
$
and all $ p \in[1,\infty) $
and the Burkholder--Davis--Gundy
inequality in
Theorem~6.3.10
in Stroock~\cite{s93}
shows the
existence of real numbers
$ K_p \in[0, \infty) $, $ p \in[1, \infty) $,
such that
\begin{eqnarray*}
&&\hspace*{-4pt} \Biggl\llVert \mathbb{E} \bigl[ f( X ) \bigr] - \frac{ 1 }{ N } \sum
_{ k = 1 }^{ N } f \bigl( \bar{Y}^{ 1, 0, k }
\bigr) - \sum_{ l = 1 }^{ \ld( N ) } \frac{ 2^l }{ N }
\Biggl( \sum_{ k = 1 }^{{ N }/{ 2^l } } f \bigl(
\bar{Y}^{ 2^l, l, k } \bigr) - f \bigl( \bar{Y}^{ 2^{ ( l - 1 ) }, l, k } \bigr) \Biggr)
\Biggr\rrVert_{ L^p( \Omega; \mathbb{R} ) }
\\
&&\hspace*{-4pt}\qquad\leq \mathbb{E} \bigl[ \bigl| f( X )
 - f \bigl( \bar{Y}^{ N, 0, 1 } \bigr) \bigr|
\bigr] + \frac{ K_p }{ \sqrt{ N } } \bigl\llVert \mathbb{E} \bigl[ f \bigl(
\bar{Y}^{ 1, 0, 1 } \bigr) \bigr] - f \bigl( \bar{Y}^{ 1, 0, 1 } \bigr)
\bigr\rrVert_{ L^p( \Omega; \mathbb{R} ) }
\\
&&\hspace*{-4pt}\qquad\quad{} + \sum_{ l = 1 }^{ \ld( N ) } \frac{ 2^{ { l }/{ 2 } } K_p }{ \sqrt{ N } }
\bigl\llVert \mathbb{E} \bigl[ f \bigl( \bar{Y}^{ 2^l, 0, 1 } \bigr) \bigr] -
\mathbb{E} \bigl[ f \bigl( \bar{Y}^{ 2^{ ( l - 1 ) }, 0, 1 } \bigr) \bigr] \\
&&\hspace*{98pt}\qquad{}- f \bigl(
\bar{Y}^{ 2^l, 0, 1 } \bigr) + f \bigl( \bar{Y}^{ 2^{ ( l - 1 ) }, 0, 1 } \bigr) \bigr
\rrVert_{ L^p( \Omega; \mathbb{R} ) }
\end{eqnarray*}
for all $ N \in\{ 2^1, 2^2, 2^3, \ldots\}
$ and all $ p \in[1,\infty) $.
In the next step estimate~\eqref{eqtamedconvergence1},
H{\"o}lder's inequality and the
triangle inequality show
\begin{eqnarray*}
&&\hspace*{-4pt} \Biggl\llVert \mathbb{E} \bigl[ f( X ) \bigr] - \frac{ 1 }{ N } \sum
_{ k = 1 }^{ N } f \bigl( \bar{Y}^{ 1, 0, k }
\bigr) - \sum_{ l = 1 }^{ \ld( N ) } \frac{ 2^l }{ N }
\Biggl( \sum_{ k = 1 }^{{ N }/{ 2^l } } f \bigl(
\bar{Y}^{ 2^l, l, k } \bigr) - f \bigl( \bar{Y}^{ 2^{ ( l - 1 ) }, l, k } \bigr) \Biggr)
\Biggr\rrVert_{ L^p( \Omega; \mathbb{R} ) }
\\
&&\hspace*{-4pt}\qquad\leq c \bigl( 1 + \llVert X \rrVert^c_{ L^{ 2 c }( \Omega; C( [0,T], \mathbb{R}^d ) ) } + \bigl
\llVert \bar{Y}^{ N, 0, 1 } \bigr\rrVert^c_{ L^{ 2 c }( \Omega; C( [0,T], \mathbb{R}^d ) ) }
\bigr)\\
&&\hspace*{-4pt}\qquad\quad{}\times\bigl\llVert X - \bar{Y}^{ N, 0, 1 } \bigr\rrVert_{ L^2( \Omega; C( [0,T], \mathbb{R}^d ) ) }
%&&\hspace*{-4pt}\qquad\quad{}
+ \frac{ 2 K_p }{ \sqrt{ N } } \bigl\llVert f \bigl( \bar{Y}^{ 1, 0, 1 } \bigr)
\bigr\rrVert_{ L^p( \Omega; \mathbb{R} ) }\\
&&\hspace*{-4pt}\qquad\quad{} + \sum_{ l = 1 }^{ \ld( N ) }
\frac{
2^{
( { l }/{ 2 } + 1 )
}
K_p
}{ \sqrt{ N } } \bigl\llVert f \bigl( \bar{Y}^{ 2^l, 0, 1 } \bigr) - f
\bigl( \bar{Y}^{ 2^{ ( l - 1 ) }, 0, 1 } \bigr) \bigr\rrVert_{ L^p( \Omega; \mathbb{R} ) },
\end{eqnarray*}
and
Corollary~\ref{thmtamedconvergence2}
and again
estimate~\eqref{eqtamedconvergence1}
hence give\vspace*{-1pt}
\begin{eqnarray*}
&&\hspace*{-4pt} \Biggl\llVert \mathbb{E} \bigl[ f( X ) \bigr] - \frac{ 1 }{ N } \sum
_{ k = 1 }^{ N } f \bigl( \bar{Y}^{ 1, 0, k }
\bigr) - \sum_{ l = 1 }^{ \ld( N ) } \frac{ 2^l }{ N }
\Biggl( \sum_{ k = 1 }^{{ N }/{ 2^l } } f \bigl(
\bar{Y}^{ 2^l, l, k } \bigr) - f \bigl( \bar{Y}^{ 2^{ ( l - 1 ) }, l, k } \bigr) \Biggr)
\Biggr\rrVert_{ L^p( \Omega; \mathbb{R} ) }
\\
&&\hspace*{-4pt}\qquad\leq 2 c R_2 \Bigl( 1 + \sup_{ M \in\N} \bigl\llVert
\bar{Y}^{ M, 0, 1 } \bigr\rrVert^c_{ L^{ 2 c }( \Omega; C( [0,T], \mathbb{R}^d ) ) }
\Bigr)
\frac{ \sqrt{ 1 + \ld(N) } }{ \sqrt{ N } } \\
&&\hspace*{-4pt}\qquad\quad{}+ \frac{ 2 K_p }{ \sqrt{ N } } \bigl\llVert f
\bigl(\bar{Y}^{ 1, 0, 1 } \bigr) \bigr\rrVert_{ L^p( \Omega; \mathbb{R} ) }
\\
&&\hspace*{-4pt}\qquad\quad{} + \sum_{ l = 1 }^{ \ld( N ) } \frac{
2^{ ( { l }/{ 2 } + 2 ) }
c K_p
}{
\sqrt{ N }
}
\Bigl( 1 + \sup_{ M \in\N} \bigl\llVert \bar{Y}^{ M, 0, 1 } \bigr
\rrVert^c_{ L^{ 2 p c }( \Omega; C( [0,T], \mathbb{R}^d ) ) } \Bigr)\\
&&\hspace*{-4pt}\quad\qquad{}\times \bigl\llVert \bar{Y}^{ 2^l, 0, 1 }
- \bar{Y}^{ 2^{ ( l - 1 ) }, 0, 1 } \bigr\rrVert_{ L^{ 2 p }( \Omega; C( [0,T], \mathbb{R}^d ) ) }
\end{eqnarray*}
for all $ N \in\{ 2^1, 2^2, 2^3, \ldots\} $
and all $ p \in[1,\infty) $.
The triangle inequality, again Corollary~\ref{thmtamedconvergence2}
and the estimate
$ \| f( v ) \|_{
C( [0,T], \mathbb{R}^d )
} \leq(
2 c + \| f( 0 )
\|_{ C( [0,T], \mathbb{R}^d ) }
)( 1 + \| v \|^{ ( c + 1 )
}_{ C( [0,T], \mathbb{R}^d ) } ) $
for all $ v \in C( [0,T],
\mathbb{R}^d ) $ then yield
\begin{eqnarray*}
&&\hspace*{-4pt} \Biggl\llVert \mathbb{E} \bigl[ f( X ) \bigr] - \frac{ 1 }{ N } \sum
_{ k = 1 }^{ N } f \bigl( \bar{Y}^{ 1, 0, k }
\bigr) - \sum_{ l = 1 }^{ \ld( N ) } \frac{ 2^l }{ N }
\Biggl( \sum_{ k = 1 }^{{ N }/{ 2^l } } f \bigl(
\bar{Y}^{ 2^l, l, k } \bigr) - f \bigl( \bar{Y}^{ 2^{ ( l - 1 ) }, l, k } \bigr) \Biggr)
\Biggr\rrVert_{ L^p( \Omega; \mathbb{R} ) }
\\
&&\hspace*{-4pt}\qquad\leq 2 c R_2 \Bigl( 1 + \sup_{ M \in\N} \bigl\llVert
\bar{Y}^{ M, 0, 1 } \bigr\rrVert^c_{ L^{ 2 c }( \Omega; C( [0,T], \mathbb{R}^d ) ) } \Bigr)
\frac{ \sqrt{ 1 + \ld(N) } }{ \sqrt{ N } }
\\
&&\hspace*{-4pt}\qquad\quad{} + 2 K_p \bigl( 2 c + \bigl\llVert f( 0 ) \bigr
\rrVert_{ C( [0,T], \mathbb{R}^d ) } \bigr) \bigl( 1 + \bigl\llVert \bar{Y}^{ 1, 0, 1 }
\bigr\rrVert^{
( c + 1 )
}_{
L^{ p ( c + 1 )
}( \Omega; C( [0,T], \mathbb{R}^d ) )
} \bigr) \frac{ 1 }{ \sqrt{ N } }
\\
&&\hspace*{-4pt}\qquad\quad{} + c K_p R_{ 2 p } \Bigl( 1 + \sup_{ M \in\N} \bigl
\llVert \bar{Y}^{ M, 0, 1 } \bigr\rrVert^c_{ L^{ 2 p c }( \Omega; C( [0,T], \mathbb{R}^d ) ) } \Bigr)
\\
&&\hspace*{-4pt}\qquad\quad{}\times\sum_{ l = 1 }^{ \ld( N ) } \frac{
2^{ ( { l }/{ 2 } + 3 ) }
\sqrt{ 1 + \ld( 2^{ l } ) }
}{
2^{ { ( l - 1 ) }/{ 2 } }
\sqrt{ N }
}
\end{eqnarray*}
and finally
\begin{eqnarray*}
&&\hspace*{-4pt} \Biggl\llVert \mathbb{E} \bigl[ f( X ) \bigr] - \frac{ 1 }{ N } \sum
_{ k = 1 }^{ N } f \bigl( \bar{Y}^{ 1, 0, k }
\bigr) - \sum_{ l = 1 }^{ \ld( N ) } \frac{ 2^l }{ N }
\Biggl( \sum_{ k = 1 }^{{ N }/{ 2^l } } f \bigl(
\bar{Y}^{ 2^l, l, k } \bigr) - f \bigl( \bar{Y}^{ 2^{ ( l - 1 ) }, l, k } \bigr) \Biggr)
\Biggr\rrVert_{ L^p( \Omega; \mathbb{R} ) }
\\
&&\hspace*{-4pt}\qquad\leq 2 c R_2 \Bigl( 1 + \sup_{ M \in\N} \bigl\llVert
\bar{Y}^{ M, 0, 1 } \bigr\rrVert^c_{ L^{ 2 c }( \Omega; C( [0,T], \mathbb{R}^d ) ) } \Bigr)
\frac{
( 1 + \ld(N)  )^{ { 3 }/{ 2 } }
}{ \sqrt{ N } }
\\
&&\hspace*{-4pt}\qquad\quad{} + 2 K_p \bigl( 2 c + \bigl\llVert f( 0 ) \bigr
\rrVert_{ C( [0,T], \mathbb{R}^d ) } \bigr) \bigl( 1 + \bigl\llVert \bar{Y}^{ 1, 0, 1 }
\bigr\rrVert^{ ( c + 1 ) }_{ L^{ p ( c + 1 ) }( \Omega; C( [0,T], \mathbb
{R}^d ) ) } \bigr)\\
&&\hspace*{-4pt}\quad\qquad{}\times \frac{
( 1 + \ld(N)  )^{ { 3 }/{ 2 } }
}{ \sqrt{ N } }
\\
&&\hspace*{-4pt}\qquad\quad{} + 12 c K_p R_{ 2 p } \Bigl( 1 + \sup_{ M \in\N} \bigl
\llVert \bar{Y}^{ M, 0, 1 } \bigr\rrVert^c_{ L^{ 2 p c }( \Omega; C( [0,T], \mathbb{R}^d ) ) } \Bigr)
\frac{
( 1 + \ld(N)  )^{ { 3 }/{ 2 } }
}{ \sqrt{ N } }
\end{eqnarray*}
for all
$ N \in\{ 2^1, 2^2, 2^3, \ldots\} $
and all $ p \in[1,\infty) $.
This shows
\eqref{eqtamedconvergence2}.
Inequality~\eqref{eqtamedconvergence4}
then immediately follows from
Lemma 2.1 in Kloeden und Neuenkirch~\cite{kn07}.
This completes the proof of Proposition~\ref{thmtamedconvergence}.
\end{pf*}

It is well known that the multilevel Monte Carlo method
combined with the (fully)
implicit Euler method
converges
%in the strong sense
too.
The following simulation
indicates that this
multilevel Monte Carlo implicit
Euler method
is considerably slower than the
multilevel Monte Carlo tamed
Euler method.
We choose a multi-dimensional
Langevin equation as an example.
More precisely, we consider the
motion of a Brownian particle of
unit mass in
the $d$-dimensional
potential
$ \frac{1}{4} \|x\|^4 - \frac{1}{2} \| x \|^2
$, $ x \in\mathbb{R}^d $,
with $ d=10 $.
The corresponding force
on the particle is then
$ x - \|x\|^2 \cdot x $
for $x\in\R^d$.
More formally,
let $T=1$, $m=d=10$,
$\xi= (0, 0, \ldots, 0  ) $,
$\mu(x)=x-\llVert x\rrVert^2\cdot x$,
for all $x\in\mathbb{R}^d$,
and
let
$\sigma(x)= \mathbf{I} $
be the identity matrix
for all $x\in\mathbb{R}^d$.
Thus the SDE~\eqref{eqSDE} reduces to
the Langevin equation
%
%e100 #&#
\begin{equation}
\label{eqLangevin} dX_t = \bigl( X_t - \llVert
X_t\rrVert^2\cdot X_t \bigr) \,dt +
dW_t^{0,1}, \qquad X_0 = \xi
\end{equation}
for $ t \in[0,1] $.
Then the
implicit Euler scheme for the SDE~\eqref{eqLangevin} is given by
mappings
$ \tilde{\tilde{Y}}{}^N_n \dvtx\Omega\rightarrow\mathbb{R}^d $,
$ n \in\{ 0, 1, \ldots, N \} $, $ N \in\mathbb{N} $,
satisfying $ \tilde{\tilde{Y}}{}^N_0 = \xi$ and
%
%e101 #&#
\begin{eqnarray}
\label{eqimplicitEulerscheme} \qquad\tilde{\tilde{Y}}{}^N_{ n + 1 } =
\tilde{\tilde{Y}}{}^N_n + \frac{ T }{ N } \cdot \bigl(
\tilde{\tilde{Y}}{}^N_{ n + 1 } - \bigl\llVert \tilde{
\tilde{Y}}{}^N_{ n + 1 } \bigr\rrVert^2 \cdot \tilde{
\tilde{Y}}{}^N_{ n + 1 } \bigr) + \bigl( W_{ { (n + 1) T }/{ N } }^{0,1}
- W_{ { n T }/{ N } }^{0,1} \bigr)
\end{eqnarray}
for all $ n \in\{ 0, 1, \ldots, N - 1 \} $ and all $ N \in\mathbb{N}
$.
Note that we used
the \textsc{Matlab} function
$\operatorname{fsolve}(\ldots)$ in our implementation
of the implicit Euler
scheme~\eqref{eqimplicitEulerscheme}.
%%%%%%%%%%%%%%%%%%%%
%
%f5 #&#
\begin{figure}

\includegraphics{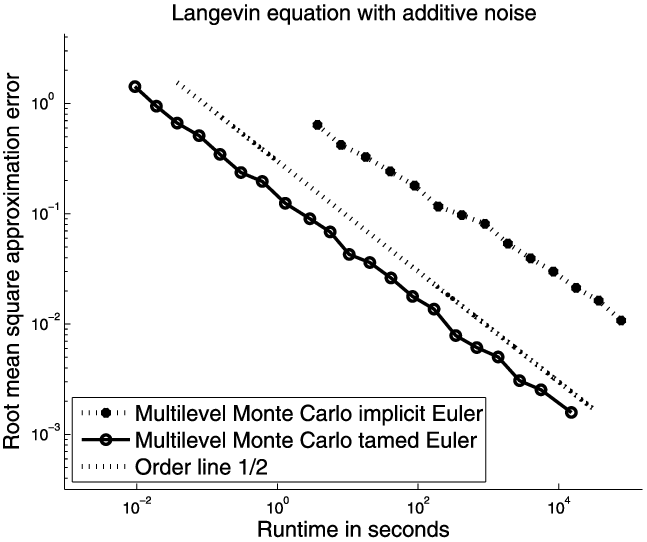}

\caption{Root mean square
approximation error
for
the uniform second moment
$\E [\sup_{t\in[0,1]}\llVert X_t\rrVert^2 ]$
of the exact solution
of~\protect\eqref{eqLangevin}
as function of the runtime
both for
the
multilevel Monte Carlo implicit Euler method
and for the
multilevel Monte Carlo tamed Euler method.}
\label{fcompareml}
\end{figure}

%
%%%%%%%%%%%%%%%%%%%%
Figure~\ref{fcompareml} displays
the root mean square
approximation error
of the
multilevel Monte Carlo implicit Euler method
for
the uniform second moment
$\E [\sup_{t\in[0,1]}\llVert X_t\rrVert^2 ]$
of the exact solution
of~\eqref{eqLangevin}
as function of the runtime when
$N\in\{2^5,2^6,\ldots,2^{18}\}$.
In addition
Figure~\ref{fcompareml} shows
the root mean square
approximation error
of the
multilevel Monte Carlo tamed Euler method
for the uniform second moment
$\E [\sup_{t\in[0,1]}\llVert X_t\rrVert^2 ]$
of the exact solution
of~\eqref{eqLangevin}
as function of the runtime when
$N\in\{2^5,2^6,\ldots,2^{25}\}$.
We see that both numerical approximations of the SDE~\eqref{eqLangevin}
apparantly converge with rate close to $\frac{1}{2}$.
Moreover the
multilevel Monte Carlo implicit Euler method
was considerably slower than the
multilevel Monte Carlo tamed Euler method.
This is presumably due to the additional computational effort
which is required to determine the zero of a nonlinear equation
in each time step of the implicit Euler method~\eqref{eqimplicitEulerscheme}.
More results on implicit
%and other strongly convergent
numerical methods
for SDEs can be found
in~\cite{h96,hms02,t02,hmps11,s10,MaoSzpruch2012pre,ms12},
for instance.
%

%sA #&#
\begin{appendix}
\section*{\texorpdfstring{Appendix: Proof of Lemma~\lowercase{\protect\ref{lem_n1}}}
{Appendix: Proof of Lemma 4.9}}
\label{secappendix}

Before
the proof of
Lemma~\ref{lem_n1}
is presented,
a few auxiliary results
(Lemmas~\ref{lAN1}--\ref{lAN4})
are established.\vadjust{\goodbreak}
%
%leA.1 #&#
\begin{lemmas} \label{lAN1}
Assume that the setting
described
in
Sections~\ref{secmultilevelpathwisedivergence}
and~\ref{secproofconj}
is fulfilled.
Then we have
%
%eA.1 #&#
\begin{equation}
\sum_{ n = 1 }^{ \infty} \mathbb{ P } \bigl[
A_{ 2^n }^{ (1) } \bigr] < \infty.
\end{equation}
\end{lemmas}
\begin{pf}
The definition~\eqref{eqdefLN} of
$ L_N $,
$ N \in\{ 2^1, 2^2, 2^3, \ldots\} $,
and independence of $\xi^{l,k}$,
$l\in\N_0$, $k\in\N$,
imply
\begin{eqnarray*}
\mathbb{P} \bigl[A_N^{(1)} \bigr] &=& \P
\biggl[ \forall l \in\mathbb{N} , \bigl\lfloor 2 \operatorname{ld} \bigl(
\sigmab^2 T^{
{1}/{2}
} \operatorname{ln}(N) \bigr) \bigr\rfloor
\leq l \leq\operatorname{ld}(N)\ \forall k \in \biggl\{ 1, 2, \ldots,
\frac{N}{2^l} \biggr\} \dvtx\\
&&\hspace*{188pt}\qquad{} \bigl| \xi^{l,k}\bigr | \leq 2^{ {l}/{4} }
T^{ -{1}/{4} } \biggr]
\\
&=& \prod_{
l =
\lfloor
2
\operatorname{ld}(
\sigmab^2
\sqrt{T}
\operatorname{ln}(N)
)
\rfloor
}^{\ld(N)} \prod
_{k=1}^{{N}/{2^l}} \P \bigl[ \bigl| \xi^{l,k}\bigr | \leq
2^{{ l }/{ 4 } } T^{ - { 1 }/{ 4 } } \bigr]\\
& =& \prod_{
l =
\lfloor
2
\operatorname{ld}(
\sigmab^2
\sqrt{T}
\operatorname{ln}(N)
)
\rfloor
}^{\ld(N)}
\bigl( \P \bigl[ \bigl| \xi^{0,1}\bigr | \leq 2^{{ l }/{ 4 } } T^{ - { 1 }/{ 4 } }
\bigr] \bigr)^{ {N}/{2^l}}
\\
& =& \prod_{
l =
\lfloor
2
\operatorname{ld}(
\sigmab^2
\sqrt{T}
\operatorname{ln}(N)
)
\rfloor
}^{\ld(N)} \bigl( 1 - \P \bigl[
\sigmab^{-1} \bigl| \xi^{0,1}\bigr | > 2^{{ l }/{ 4 } }
\sigmab^{ - 1 } T^{ - { 1 }/{ 4 } } \bigr] \bigr)^{ {N}/{2^l}} %
\end{eqnarray*}
for all $ N \in\{ N_0, 2^1 N_0,
2^2 N_0, \ldots\} $.
The inequality
%
%eA.2 #&#
%eA.3 #&#
\begin{eqnarray}
\P \bigl[ \sigmab^{-1} \bigl| \xi^{0,1}\bigr | > x
\bigr] &=& 2 \cdot\P \bigl[ \sigmab^{-1} \xi^{0,1} > x \bigr]\nonumber\\
& =&
2 \int_x^{ \infty} \frac{ 1 }{ \sqrt{ 2 \pi} }
e^{ - { y^2 }/{ 2 } } \,dy \geq 2 \int_x^{
x \sqrt{ { 3 }/{ 2 } }
}
\frac{ 1 }{ \sqrt{ 2 \pi} } e^{ - { y^2 }/{ 2 } } \,dy
\\
& \geq& \frac{ 2 }{ \sqrt{ 2 \pi} } \biggl( x \sqrt{ \frac{ 3 }{ 2 } } - x \biggr)
e^{ - {3 x^2}/{(2\cdot2)} } = \frac{ x }{ \sqrt{ \pi} } ( \sqrt{ 3 } - \sqrt{ 2 } )
e^{ - ({3 }/{4})x^2 }\nonumber\\
& =& \frac{
x
e^{ -( {3 }/{4})x^2 }
}{
\sqrt{\pi}  ( \sqrt{3} + \sqrt{2}  )
} \geq \frac{
1
}{ 6 } x
e^{ - ({3}/{4})x^2 } \nonumber %
\end{eqnarray}
for all $x\in[0,\infty)$
therefore yields
%
%eA.4 #&#
%eA.5 #&#
%eA.6 #&#
%eA.7 #&#
%eA.8 #&#
%eA.9 #&#
\begin{eqnarray}
\mathbb{P} \bigl[A_N^{(1)} \bigr] & \leq& \prod
_{
l =
\lfloor
2
\operatorname{ld}(
\sigmab^2
\sqrt{T}
\operatorname{ln}(N)
)
\rfloor
}^{\ld(N)} \biggl( 1 - \frac{1}{6}
\frac{
2^{{ l }/{ 4 } }
}{
\sigmab
T^{ { 1 }/{ 4 } }
} \cdot \exp \biggl( - \frac{3}{4}\cdot \frac{2^{ {l}/{2} }
}{
\sigmab^2 \sqrt{ T }
}
\biggr) \biggr)^{ {N}/{2^l}}
\nonumber\\
& \leq& \biggl( 1 - \frac{
2^{
\lfloor
2
\operatorname{ld}(
\sigmab^2
\sqrt{T}
\operatorname{ln}(N)
)
\rfloor
/ 4
}
}{
6 \sigmab
T^{ { 1 }/{ 4 } }
} \nonumber\\
&&\hspace*{3pt}{}\times \exp \biggl( - \frac{
3 \cdot
2^{
\lfloor
2
\operatorname{ld}(
\sigmab^2
\sqrt{T}
\operatorname{ln}(N)
)
\rfloor
/ 2
}
}{
4 \sigmab^2 \sqrt{ T }
}
\biggr) \biggr)^{
{ N }/{
2^{
\lfloor
2
\operatorname{ld}(
\sigmab^2
\sqrt{T}
\operatorname{ln}(N)
)
\rfloor
}
}
}
\nonumber\\
& \leq& \biggl( 1 - \frac{
1
}{
6 \sigmab
T^{ { 1 }/{ 4 } }
} \cdot \exp \biggl( - \frac{
3 \cdot
2^{
\lfloor
2
\operatorname{ld}(
\sigmab^2
\sqrt{T}
\operatorname{ln}(N)
)
\rfloor
/ 2
}
}{
4 \sigmab^2 \sqrt{ T }
}
\biggr) \biggr)^{
{ N }/{
2^{
\lfloor
2
\operatorname{ld}(
\sigmab^2
\sqrt{T}
\operatorname{ln}(N)
)
\rfloor
}
}
}
\nonumber
\\[-8pt]
\\[-8pt]
\nonumber
& \leq& \biggl( 1 - \frac{
1
}{
6 \sigmab
T^{ { 1 }/{ 4 } }
} \cdot \exp \biggl( - \frac{
3 \cdot
2^{
\operatorname{ld}(
\sigmab^2
\sqrt{T}
\operatorname{ln}(N)
)
}
}{
4 \sigmab^2 \sqrt{ T }
}
\biggr) \biggr)^{
{ N }/{
2^{
\lfloor
2
\operatorname{ld}(
\sigmab^2
\sqrt{T}
\operatorname{ln}(N)
)
\rfloor
}
}
}
\\
& \leq& \biggl( 1 - \frac{
1
}{
6 \sigmab
T^{ { 1 }/{ 4 } }
} \cdot \exp \biggl( - \frac{
3 \cdot
\sigmab^2
\sqrt{T}
\operatorname{ln}(N)
}{
4 \sigmab^2 \sqrt{ T }
}
\biggr) \biggr)^{
{ N }/{
2^{
2
\operatorname{ld}(
\sigmab^2
\sqrt{T}
\operatorname{ln}(N)
)
}
}
}
\nonumber\\
& =& \biggl( 1 - \frac{
1
}{
6 \sigmab
T^{ { 1 }/{ 4 } }
} \cdot \exp \biggl( - \frac{ 3 }{ 4 }
\operatorname{ln}(N) \biggr) \biggr)^{
{ N }/{
2^{
\operatorname{ld}(
(
\sigmab^2
\sqrt{T}
\operatorname{ln}(N)
)^2
)
}
}
}\nonumber\\
& = &\biggl( 1 -
\frac{
1
}{
6 \sigmab
T^{ { 1 }/{ 4 } }
} \cdot N^{-{3}/{4}} \biggr)^{
{ N }/{
(
\sigmab^2
\sqrt{T}
\operatorname{ln}(N)
)^2
}
} \nonumber%
\end{eqnarray}
for all $N\in\{N_0,2^1N_0,2^2N_0,\ldots\}$.
Next we estimate $1-x\leq\exp(-x)$ for all $x\in\R$ to get
\begin{eqnarray*}
\sum_{ n = 1 }^{ \infty} \mathbb{ P }
\bigl[ A_{ 2^n }^{ (1) } \bigr] &=& \mathop{\sum
_{
N \in\{2^1,2^2,2^3,\ldots\}}}_{ N< N_0
} \P \bigl[A_N^{(1)}
\bigr] +\mathop{\sum_{
N \in\{2^1,2^2,2^3,\ldots\} }}_{ N\geq N_0
} \P
\bigl[A_N^{(1)} \bigr]\\
& \leq& N_0 +\mathop{\sum
_{
N \in\{2^1,2^2,2^3,\ldots\}}}_{ N\geq N_0
} \biggl( 1 -
\frac{
1
}{
6 \sigmab
T^{ { 1 }/{ 4 } }
} \cdot N^{-{3}/{4}} \biggr)^{
{ N }/{
(
\sigmab^2
\sqrt{T}
\operatorname{ln}(N)
)^2
}
}
\\
&\leq& N_0 +\sum_{N=N_0}^\infty
\biggl[ \exp \biggl( - \frac{
1
}{
6 \sigmab
T^{ { 1 }/{ 4 } }
} \cdot N^{-{3}/{4}} \biggr)
\biggr]^{
{ N }/{
(
\sigmab^2
\sqrt{T}
\operatorname{ln}(N)
)^2
}
}\\
& =& N_0 +\sum_{N=N_0}^\infty
\exp \biggl( - \frac{
N^{{1}/{4}}
}{
6 \sigmab^{5}
T^{ { 5 }/{ 4 } }
(
\operatorname{ln}(N)
)^2
} \biggr) <\infty. %
\end{eqnarray*}
This completes the proof
of Lemma~\ref{lAN1}.
\end{pf}

%leA.2 #&#
\begin{lemmas} \label{lAN2}
Assume that the setting
described
in
Sections~\ref{secmultilevelpathwisedivergence}
and~\ref{secproofconj}
is fulfilled.
Then we have
%
%eA.10 #&#
\begin{equation}
\label{eqAN2} \sum_{ n = 1 }^{ \infty} \mathbb{ P
} \bigl[ A_{ 2^n }^{ (2) } \bigr] < \infty.
\end{equation}
\end{lemmas}
\begin{pf}%[Proof of Lemma~\ref{lAN2}]
Subadditivity of
the probability measure
$ \mathbb{P} $ and
the
inequality
$\P[\sigmab^{-1}|\xi^{0,1}|\geq x]
\leq
\frac{1}{x}\exp (
- \frac{x^2}{2}  )$
for all $x\in(0,\infty)$
(e.g., Lemma 22.2 in~\cite{klenke2008})
imply
%
%eA.11 #&#
%eA.12 #&#
%eA.13 #&#
%eA.14 #&#
%eA.15 #&#
\begin{eqnarray}
\mathbb{P} \bigl[A_N^{(2)} \bigr] &=& \P
\biggl[ \exists l \in \bigl\{ 0, 1, 2, \ldots, \ld(N) \bigr\}\ \exists k \in \biggl
\{ 1, 2, \ldots, \frac{N}{2^l} \biggr\} \dvtx\nonumber\\
&&\hspace*{82pt}\qquad{}\bigl | \xi^{l,k} \bigr| \geq
2^{
{ (l-1) }/{ 4 }
} T^{ - {1}/{4} } N \biggr]\nonumber\\
& \leq& \sum
_{l=0}^{\ld(N)} \sum_{k=1}^{{N}/{2^l}}
\P \bigl[ \bigl| \xi^{l,k} \bigr| \geq 2^{
{ (l-1) }/{ 4 }
} T^{ - {1}/{4} } N \bigr]
\nonumber\\
& =& \sum_{l=0}^{\ld(N)} \frac{N}{2^l}
\cdot \P \bigl[ \bigl| \xi^{0,1}\bigr | \geq 2^{
{ (l-1) }/{ 4 }
} T^{ - {1}/{4} } N
\bigr]
\nonumber
\\[-8pt]
\\[-8pt]
\nonumber
&= &\sum_{l=0}^{\ld(N)} \frac{N}{2^l}
\cdot \P \biggl[ \sigmab^{-1}\bigl| \xi^{0,1}\bigr | \geq
\frac{
2^{
{ (l-1) }/{ 4 }
}
N
}{
\sigmab T^{{1}/{4} }
} \biggr]
\\
& \leq &\sum_{l=0}^{\ld(N)}
\frac{
N
}{
2^l
} \cdot \frac{ \sigmab T^{{1}/{4} }
}{
2^{
{ (l-1) }/{ 4 }
}
N
} \exp \biggl( - \frac{
2^{
{ (l-1) }/{ 2 }
}
N^2
}{
2\sigmab^2 T^{ {1}/{2} }
}
\biggr)
\nonumber\\
& \leq& \sum_{l=0}^{\ld(N)}
\frac{ \sigmab T^{{1}/{4} }
}{
2^{
{ -1 }/{ 4 }
}
} \exp \biggl( - \frac{
2^{
{ -1 }/{ 2 }
}
N^2
}{
2\sigmab^2 T^{ {1}/{2} }
} \biggr) \nonumber\\
&=& \bigl(\ld(N)+1
\bigr) \sigmab2^{{1}/{4}}T^{{1}/{4} } \exp \biggl( - \frac{
N^2
}{
2^{{3}/{2}}\sigmab^2 T^{ {1}/{2} }
}
\biggr) \nonumber%
\end{eqnarray}
for all $N\in\{2^1,2^2,2^3,\ldots\}$.
Summing over $N\in\{2^1,2^2,2^3,\ldots\}$
results in
%
%eA.16 #&#
\begin{eqnarray}
\qquad
&&\sum_{ n = 1 }^{ \infty} \mathbb{ P }
\bigl[ A_{ 2^n }^{ (2) } \bigr]
\nonumber
\\[-8pt]
\\[-8pt]
\nonumber
&&\qquad \leq \sum
_{
N \in\{2^1,2^2,2^3,\ldots\}
} \bigl(\ld(N)+1 \bigr) \sigmab2^{{1}/{4}}T^{{1}/{4} }
\exp \biggl( - \frac{
N^2
}{
2^{{3}/{2}}\sigmab^2 T^{ {1}/{2} }
} \biggr) <\infty %
,
\end{eqnarray}
and this completes the proof
of Lemma~\ref{lAN2}.
\end{pf}

%leA.3 #&#
\begin{lemmas} \label{lnormalconditionalestimate}
Let
$
( \Omega, \mathcal{F}, \P
)
$
be a probability space, and
let
$
Z \dvtx\Omega\rightarrow\mathbb{R}
$
be a standard normally distributed
$ \mathcal{F} / \mathcal{B}( \mathbb{R} ) $-measurable
mapping. Then
%
%eA.17 #&#
\begin{eqnarray}
\P \bigl[|Z|<x+y \big| |Z|\geq x \bigr] \leq5xy
\end{eqnarray}
for all $x\in[\frac12,\infty)$ and all $y\in[0,\infty)$.
\end{lemmas}
\begin{pf}%[Proof of Lemma~\ref{lnormalconditionalestimate}]
Monotonicity of the exponential function yields
%
%eA.18 #&#
\begin{eqnarray}
\label{eqleqZleq}
\P \bigl[x\leq|Z|<x+y \bigr]& =&2\cdot \P [x\leq Z<x+y ]
\nonumber
\\[-8pt]
\\[-8pt]
\nonumber
&=&2\int
_x^{x+y}\frac{1}{\sqrt{2\pi}}e^{-{z^2}/{2}}
\,dz
\leq\frac{2}{\sqrt{2\pi}}y e^{-{x^2}/{2}} %
\end{eqnarray}
for all $x,y\in[0,\infty)$.
Apply the standard estimate
$
\mathbb{P}[|Z|\geq x]
\geq
\frac{x}{1+x^2}\frac{2}{\sqrt{2\pi}}\times
\exp (-\frac{x^2}{2} )
$
for all $x\in(0,\infty)$
(e.g., Lemma 22.2 in~\cite{klenke2008}),
inequality~\eqref{eqleqZleq}
and $\frac{x^2}{1+x^2}\geq\frac{1}{5}$ for all $x\in[\frac
{1}{2},\infty)$ to get
\begin{eqnarray*}
\P \bigl[|Z|<x+y | |Z|\geq x \bigr] &=& \frac{
\P [x\leq|Z|<x+y ]
}{
\P [|Z|\geq x ]
} \\
&\leq&
\frac{
({2}/{\sqrt{2\pi}})y e^{-{x^2}/{2}}
}{
({x}/{(1+x^2)})({2}/{\sqrt{2\pi}})\exp (-{x^2}/{2} )
} \\
&=& \frac{xy}{ ({x^2}/{(1+x^2)} )} \leq5xy %
\end{eqnarray*}
for all $x\in[\frac{1}{2},\infty)$ and all $y\in[0,\infty)$.
This completes the proof
of Lemma~\ref{lnormalconditionalestimate}.
\end{pf}

%leA.4 #&#
\begin{lemmas} \label{lAN3}
Assume that the setting
described
in
Sections~\ref{secmultilevelpathwisedivergence}
and~\ref{secproofconj}
is fulfilled. Then we have
%
%eA.19 #&#
\begin{equation}
\label{eqAN3} \sum_{ n = 1 }^{ \infty} \mathbb{ P
} \bigl[ A_{ 2^n }^{ (3) } \bigr] < \infty.
\end{equation}
\end{lemmas}

\begin{pf}%[Proof of Lemma~\ref{lAN3}]
Let
$
K \dvtx
\N_0 \times
\N_0 \times
\Omega\to\N
\cup\{ \infty\}
$
be defined as
%
%eA.20 #&#
\begin{equation}
K(v,l) := \min \bigl( \bigl\{ k \in\N \dvtx \bigl\llvert \xi^{ v, k }
\bigr\rrvert \geq 2^{{ l }/{ 4 } } T^{ - { 1 }/{ 4 } } \bigr\} \cup\{ \infty\} \bigr)
\end{equation}
for all
$ v, l \in\N_0 $.
Inserting
definition~\eqref{eqA3}
we get
%
%eA.21 #&#
%eA.22 #&#
%eA.23 #&#
%eA.24 #&#
%eA.25 #&#
%eA.26 #&#
\begin{eqnarray}
\mathbb{ P } \bigl[ A_{ N }^{ (3) } \bigr]& =&
\mathbb{ P } \bigl[ \exists l \in\N, \bigl\lfloor 2 \operatorname{ld} \bigl(
\sigmab^2 T^{
{1}/{2}
} \operatorname{ln}(N) \bigr) \bigr\rfloor
\leq l\leq \operatorname{ld}(N) + 1 \dvtx \nonumber\\
&&\hspace*{19pt} 2^{ {l}/{4} } T^{ - {1}/{4} }\leq
\eta_N < 2^{ {l}/{4} } T^{ - {1}/{4} } \bigl( 1 +
5^{ (- \dl\cdot2^{ (l-1) } ) } \bigr) \bigr]
\nonumber\\
&\leq& \sum_{l= \lfloor2\ld (\sigmab^2\sqrt{T}\ln(N) )
\rfloor}^{\ld(N)+1} \mathbb{ P } \bigl[
2^{ {l}/{4} } T^{ - {1}/{4} } \leq\eta_N < 2^{ {l}/{4} }
T^{ - {1}/{4} } \bigl( 1 + 5^{ (- \dl\cdot2^{ (l-1) } ) } \bigr) \bigr]
\nonumber\\
&\leq& \sum_{l= \lfloor2\ld (\sigmab^2\sqrt{T}\ln(N) )
\rfloor}^{\ld(N)+1} \mathbb{ P } \Bigl[
\exists v\in \bigl\{1,2,\ldots,\ld(N) \bigr\}\dvtx 2^{ {l}/{4} }
T^{ - {1}/{4} } \nonumber\\
&&\hspace*{88pt}\leq\max_{k\in\{1,2,\ldots,{N}/{2^v}\}} \bigl\llvert \xi^{v,k} \bigr
\rrvert < 2^{ {l}/{4} } T^{ - {1}/{4} } \bigl( 1 + 5^{ (- \dl\cdot2^{ (l-1) } ) } \bigr)
\Bigr] %% \end{eqnarray}
%and
\\
&\leq& \sum_{l= \lfloor2\ld (\sigmab^2\sqrt{T}\ln(N) )
\rfloor}^{\ld(N)+1} \sum
_{v=1}^{\ld(N)} \mathbb{ P } \biggl[ \biggl\{ \exists k
\in \biggl\{1,\ldots,\frac{N}{2^v} \biggr\}\dvtx \bigl\llvert
\xi^{v,k} \bigr\rrvert \geq2^{ {l}/{4} } T^{ - {1}/{4} } \biggr\}
\nonumber\\
&&{} \cap\biggl\{ \forall k\in \biggl\{1,\ldots,\frac{N}{2^v} \biggr\}\dvtx \bigl
\llvert \xi^{v,k} \bigr\rrvert < 2^{ {l}/{4} } T^{ - {1}/{4} }
\bigl( 1 + 5^{ (- \dl\cdot2^{ (l-1) } ) } \bigr) \biggr\} \biggr]
\nonumber\\
&\leq& \sum_{l= \lfloor2\ld (\sigmab^2\sqrt{T}\ln(N) )
\rfloor}^{\ld(N)+1} \sum
_{v=1}^{\ld(N)} \mathbb{ P } \biggl[ \biggl\{ K(v,l)\leq
\frac{N}{2^v} \biggr\}\nonumber\\
&&\hspace*{89pt}\qquad{} \cap \bigl\{ \bigl| \xi^{ v, K(v,l) } \bigr| <
2^{ {l}/{4} } T^{ - {1}/{4} } \bigl( 1 + 5^{ (- \dl\cdot2^{ (l-1) } ) } \bigr) \bigr\}
\biggr]\nonumber %
\end{eqnarray}
for all
$
N \in
\{ N_0, 2^1 N_0, 2^2 N_0, \ldots\}
$.
The method of rejection
sampling hence
results in
%
%eA.27 #&#
%eA.28 #&#
%eA.29 #&#
\begin{eqnarray}
\mathbb{ P } \bigl[ A_{ N }^{ (3) } \bigr]
& \leq &\sum_{
l =
\lfloor
2
\ld (
\sigmab^2 \sqrt{T} \ln( N )
)
\rfloor
}^{ \ld(N) + 1 } \sum
_{ v = 1 }^{ \ld( N ) } \mathbb{ P } \bigl[ \bigl| \xi^{ v, K( v, l) }
\bigr| < 2^{{ l }/{ 4 } } T^{ - { 1 }/{ 4 } }
\bigl( 1 + 5^{ (- \dl\cdot2^{ (l-1) } ) } \bigr) ,\nonumber\\
&&\hspace*{207pt}\qquad K(v,l)<\infty \bigr]
\nonumber\\
& =& \sum_{
l =
\lfloor
2
\ld (
\sigmab^2 \sqrt{ T } \ln(N)
)
\rfloor
}^{ \ld(N) + 1 } \sum
_{ v = 1 }^{ \ld( N ) } \mathbb{ P } \bigl[ \bigl| \xi^{0,1}
\bigr| < 2^{ {l}/{4} } T^{ - {1}/{4} }
 \bigl( 1 + 5^{ (- \dl\cdot2^{ (l-1) } ) } \bigr) |
\nonumber\\
&&\hspace*{165pt}\qquad \bigl|\xi^{ 0, 1} \bigr|\geq 2^{ {l}/{4} } T^{ - {1}/{4} } \bigr]
\\
& =& \sum_{
l =
\lfloor
2 \ld (
\sigmab^2 \sqrt{T} \ln(N)
)
\rfloor
}^{ \ld(N) + 1 } \ld(N) \cdot \mathbb{
P } \biggl[ \bigl| \sigmab^{-1} \xi^{ 0, 1} \bigr| < \frac{
2^{{ l }/{ 4 } }
(
1 +
5^{ ( - \dl\cdot2^{ (l-1) } ) }
)
}{
\sigmab T^{{1}/{4} }
}
\Big|\nonumber\\
&&\hspace*{180pt}\quad{}  \bigl\llvert \sigmab^{-1}\xi^{0,1} \bigr\rrvert \geq
\frac{
2^{ {l}/{4} }
}{
\sigmab T^{{1}/{4} }
} \biggr] \nonumber%
\end{eqnarray}
for all
$
N \in
\{ N_0, 2^1 N_0, 2^2 N_0, \ldots\}
$.
In order to apply
Lemma~\ref{lnormalconditionalestimate},
we note that
%
%eA.30 #&#
%eA.31 #&#
\begin{eqnarray}
\frac{2^{{l}/{4}}}{\sigmab T^{{1}/{4}}} &\geq& \frac{
2^{ \lfloor2\ld (\sigmab^2\sqrt{T}\ln(N) )
\rfloor/4}
}{
\sigmab T^{{1}/{4}}
} \geq
\frac{2^{ (2\ld (\sigmab^2\sqrt{T}\ln(N) )-1
)/4}}{\sigmab T^{{1}/{4}}} \nonumber\\
&= &\frac{2^{\ld (\sigmab^2\sqrt{T}\ln(N) )/2}}{\sigmab
T^{{1}/{4}}2^{{1}/{4}}}
\\
&= &\frac{\sqrt{\sigmab^2\sqrt{T}\ln(N)}}{\sigmab T^{{1}/{4}}2^{{1}/{4}}} = \frac{\sqrt{\ln(N)}}{ 2^{{1}/{4}} } \geq \frac{\sqrt{\ln(2)}}{ 2^{{1}/{4}} } \geq
\frac12\nonumber %
\end{eqnarray}
for all $l \in\mathbb{N} \cap [  \lfloor2 \operatorname
{ld} (
\sigmab^2 T^{{1}/{2}}\operatorname{ln}(N) ) \rfloor,
\infty )$
and
all $N\in\{2^1, 2^2, 2^3, \ldots\}$.
Lemma~\ref{lnormalconditionalestimate}
applied
to the standard normally distributed variable $\sigmab^{-1}\xi^{0,1}$
thus leads to
%
%eA.32 #&#
%eA.33 #&#
%eA.34 #&#
%eA.35 #&#
\begin{eqnarray}
\mathbb{ P } \bigl[ A_{ N }^{ (3) } \bigr] & \leq&
\ld(N) \Biggl[ \sum_{
l =
\lfloor
2
\ld (
\sigmab^2 \sqrt{T} \ln(N)
)
\rfloor
}^{ \ld(N) + 1 } 5 \cdot
\frac{
2^{{ l }/{ 4 } }
}{
\sigmab T^{ { 1 }/{ 4 } }
} \cdot \frac{
2^{{ l }/{ 4 } }
\cdot
5^{ ( - \dl\cdot2^{ (l-1) } ) }
}{
\sigmab T^{ { 1 }/{ 4 } }
} \Biggr]
\nonumber\\
& =& \ld(N) \Biggl[ \sum_{
l =
\lfloor
2
\ld (
\sigmab^2 \sqrt{T} \ln(N)
)
\rfloor
}^{ \ld(N) + 1 } 5
\cdot \frac{
2^{ {l}/{2} }
}{
\sigmab^2 T^{ {1}/{2} }
} \cdot 5^{ (- ({\dl}/{2}) \cdot2^{ l } ) } \Biggr]
\nonumber\\
& \leq& \ld(N) \Biggl[ \sum_{
l =
\lfloor
2
\ld (
\sigmab^2
\sqrt{T} \ln(N)
)
\rfloor
}^{ \ld(N) + 1 } 5
\cdot \frac{
2^{ ( \ld(N) + 1 ) }
}{
\sigmab^2 T^{ {1}/{2} }
} \cdot 5^{
(
- { \dl}/{ 2 }
\cdot
2^{
(
2
\ld(
\sigmab^2 \sqrt{T} \ln(N)
)
- 1
)
}
)
} \Biggr]
\\
& \leq& \bigl( \ld( N ) \bigr)^2 \cdot \frac{
10N
}{
\sigmab^2 T^{ {1}/{2} }
} \cdot
5^{
(
- { \dl}/{ 4 } \cdot
(
\sigmab^2 \sqrt{T} \ln(N)
)^2
)
} \nonumber\\[-2pt]
&\leq& \frac{ 10 N^3
}{
\bar{ \sigma}^2
T^{ { 1 }/{ 2 } }
} \cdot \exp \biggl( -
\frac{
\delta\bar{ \sigma}^4 T
( \ln(N)  )^2
}{ 4 } \biggr)\nonumber %
\end{eqnarray}
for all
$
N \in\{ N_0, 2^1 N_0, 2^2 N_0, \ldots\}
$.
Summing over
$
N \in\{ 2^1, 2^2, 2^3, \ldots\}
$
results~in
%
%eA.36 #&#
%eA.37 #&#
\begin{eqnarray}
\sum_{ n = 1 }^{ \infty} \mathbb{ P }
\bigl[ A_{ 2^n }^{ (3) } \bigr] & =& \mathop{\sum
_{
N \in\{ 2^1, 2^2, 2^3, \ldots\}}}_{ N< N_0
} \P \bigl[ A_N^{ (3) }
\bigr] + \mathop{\sum_{
N \in\{ 2^1, 2^2, 2^3, \ldots\}}}_{
N \geq N_0
} \P
\bigl[ A_N^{ (3) } \bigr]
\nonumber
\\[-9pt]
\\[-9pt]
\nonumber
& \leq& N_0 + \sum_{ N = N_0 }^{ \infty}
\frac{
10
}{
\bar{ \sigma}^2
T^{ { 1 }/{ 2 } }
} \cdot N^{
(
3 - \delta\bar{\sigma}^4 T \ln(N) / 4
)
} < \infty. %
\end{eqnarray}
This completes the proof
of Lemma~\ref{lAN3}.
\end{pf}

%leA.5 #&#
\begin{lemmas} \label{lAN4}
Assume that the setting
described
in
Sections~\ref{secmultilevelpathwisedivergence}
and~\ref{secproofconj}
is fulfilled.
Then
we have
%
%eA.38 #&#
\begin{equation}
\label{eqAN4} \sum_{ n = 1 }^{ \infty} \mathbb{ P
} \bigl[ A_{ 2^n }^{ (4) } \bigr] < \infty.
\end{equation}
\end{lemmas}
\begin{pf}%[Proof of Lemma~\ref{lAN4}]
First of all, define
a filtration
$
\tilde{ \mathcal{F} }_{ l }^N
$,
$
l \in\{ 0, 1, \ldots, \ld(N) - 1 \}
$,
through
%
%eA.39 #&#
\begin{eqnarray}
\label{eqdeffiltration} \tilde{\mathcal{F}}_{l}^N :=
\sigma_{ \Omega} \bigl( \xi^{ v, k } , k \in\mathbb{N}, v \in \bigl\{
\operatorname{ld}( N ) - l, \operatorname{ld}(N) - l + 1, \ldots,
\operatorname{ld}(N) \bigr\} \bigr)
\end{eqnarray}
for all
$
l \in\{ 0, 1, \ldots, \ld(N) - 1
\}
$
and
every
$
N \in\{ 2^1, 2^2, 2^3, \ldots\}
$
where
$ \sigma_{ \Omega}( \cdot) $ denotes
the smallest sigma-algebra generated
by its argument.
Moreover,
define an
$
\mathcal{F} / \mathcal{B}( \mathbb{R} )
$-measurable
mapping
$
\tilde{L}_N \dvtx
\Omega\to
\{ 0, 1, \ldots, \ld(N) - 1 \}
$
through
$
\tilde{ L }_N := \ld(N) - L_N
$
for every
$
N \in
\{ 2^1, 2^2, 2^3, \ldots\}
$.
Next observe that the identity
\begin{eqnarray*}
 \tilde{L}_N&= &\ld(N)-L_N \\[-2pt]
 &=& \ld(N)-\max
\biggl( \{ 1 \} \cup \biggl\{ l \in \bigl\{ 1, 2, \ldots, \ld(N) \bigr\}
\dvtx\\[-2pt]
&&\hspace*{97pt}\exists k \in \biggl\{ 1, 2, \ldots, \frac{N}{2^l} \biggr\} \dvtx \bigl|
\xi^{l,k} \bigr| > 2^{{ l }/{ 4 } } T^{ - {1}/{4} } \biggr\} \biggr)
\\[-2pt]
&=& \min \biggl( \bigl\{\ld(N)-1 \bigr\}\cup \biggl\{l\in \bigl\{0,1,\ldots,
\ld(N)-1 \bigr\}\dvtx\\[-2pt]
&&\hspace*{97pt}\exists k \in \biggl\{ 1, 2, \ldots, \frac{ N }{ 2^{ ( \ld(N) - l ) } }
\biggr\} \dvtx\\[-2pt]
&&\hspace*{97pt} \bigl| \xi^{ \ld(N) - l , k } \bigr| > 2^{
{ ( \ld(N) - l ) }/{ 4 }
} T^{ - 1/4 }
\biggr\} \biggr) %
\end{eqnarray*}
for every
$ N \in\{ 2^1, 2^2, 2^3, \ldots\} $
shows that
$ \tilde{L}_N $
is a stopping time with respect to
the filtration
$
\tilde{ \mathcal{F} }_{ l }^N
$,
$ l \in\{ 0, 1, \ldots, \ld(N) - 1 \}
$,
for every
$
N \in\{ 2^1, 2^2, 2^3, \ldots\}
$.
Consequently,
the
sigma-algebras
$
\tilde{\mathcal{F}}_{ \tilde{L}_N }^N
:=
\{
A \in\mathcal{F}
\dvtx
(
\forall
l \in\{ 0, 1, \ldots, \ld( N ) \}
\dvtx
A \cap\{ \tilde{L}_N = l \}
\in\tilde{\mathcal{ F }}^N_{ l }
)
\}
$
for
$ N \in\{ 2^1, 2^2, 2^3, \ldots\} $
are well-defined.
By definition
\eqref{eqdeffiltration}
the random variables
$ \xi^{ L_N - 1, k } $,
$ k \in\N$,
are independent of
$
\tilde{ \mathcal{F} }_{ \tilde{L}_N }^N
$
for every
$
N \in\{ 2^1, 2^2, 2^3, \ldots\}
$.
Indeed, observe that
\eqref{eqdeffiltration}
shows that
%
%eA.40 #&#
%eA.41 #&#
%eA.42 #&#
\begin{eqnarray}
\label{eqabc} %
 \P \bigl[ \bigl\{ \xi^{ L_N - 1, k } \in A \bigr\}
\cap B \bigr] &=& \sum_{ l = 0 }^{ \ld(N) - 1 } \P \bigl[
\bigl\{ \xi^{ L_N - 1, k } \in A \bigr\} \cap B \cap \{ \tilde{L}_N =
l \} \bigr]
\nonumber\\
& = &\sum_{ l = 0 }^{ \ld(N) - 1 } \P \bigl[ \bigl\{
\xi^{ \ld( N ) - l - 1, k } \in A \bigr\} \cap \underbrace{ \bigl( B \cap \{
\tilde{L}_N = l \} \bigr) }_{
\in\tilde{ \mathcal{F} }^N_{ l }
} \bigr]\\
& =& \sum
_{ l = 0 }^{ \ld(N) - 1 } \P \bigl[ \xi^{ \ld( N ) - l - 1, k } \in A
\bigr] \cdot \P \bigl[ B \cap \{ \tilde{L}_N = l \} \bigr]
\nonumber\\
& =&\P \bigl[ \xi^{ 0, 1 } \in A \bigr] \Biggl( \sum
_{ l = 0 }^{ \ld(N) - 1 } \P \bigl[ B \cap \{
\tilde{L}_N = l \} \bigr] \Biggr)\nonumber\\
& =& \P \bigl[ \xi^{ 0, 1 } \in
A \bigr] \cdot \P [ B ]\nonumber %
\end{eqnarray}
for all
$ A \in\mathcal{B}(\R) $,
$
B \in
\tilde{\mathcal{F}}_{\tilde{L}_N}^N
$,
$ k \in\N$
and all
$
N \in\{ 2^1, 2^2, 2^3, \ldots\}
$.
Next we note that
$
\eta_N \dvtx
\Omega\rightarrow[0,\infty)
$
is
$
\tilde{ \mathcal{F} }_{ \tilde{L}_N }^N
/
\mathcal{B}( \R)
$-measurable
for every
$ N \in\{ 2^1, 2^2, 2^3, \ldots\} $.
Indeed, observe that
%
%eA.43 #&#
%eA.44 #&#
\begin{eqnarray}
\label{eqetaNmeasurable} %
 &&\{ \eta_N < c \} \cap \{
\tilde{L}_N = l \}\nonumber\\
&&\qquad = \biggl\{ \max \biggl\{\bigl | \xi^{L_N,k}\bigr | \in
\mathbb{R} \dvtx k \in \biggl\{ 1, 2, \ldots, \frac{ N }{ 2^{ L_N } } \biggr\} \biggr
\} < c \biggr\} \cap \{ \tilde{L}_N = l \}
\nonumber
\\[-8pt]
\\[-8pt]
\nonumber
&&\qquad = \underbrace{ \biggl\{ \max \biggl\{ \bigl| \xi^{ \ld( N ) - l , k }\bigr | \in\mathbb{R}
\dvtx k \in \biggl\{ 1, 2, \ldots, \frac{ N }{ 2^{ \ld( N ) - l } } \biggr\} \biggr\} < c
\biggr\} }_{
\in
\tilde{ \mathcal{F} }^N_l
} \cap \underbrace{ \{ \tilde{L}_N = l \}
}_{
\in
\tilde{ \mathcal{F} }^N_l
} \\
&&\qquad\in \tilde{ \mathcal{F} }^N_l\nonumber
\end{eqnarray}
for all
$ c \in\mathbb{R} $,
$
l \in\{ 0, 1, \ldots, \ld( N ) - 1 \}
$
and all
$
N \in\{ 2^1, 2^2, 2^3, \ldots\}
$.
In the next step we observe that
\eqref{eqabc},
\eqref{eqetaNmeasurable},
the fact that
$
L_N \dvtx
\Omega\rightarrow
\{ 1, 2, \ldots, \ld(N) \}
$
is measurable with respect to
$
\tilde{ \mathcal{F} }_{ \tilde{L}_N }^N
$
for all
$
N \in\{ 2^1, 2^2, 2^3, \ldots\}
$
and the inequality
$
\P [
|
| \xi^{0,1} | - x
|
\leq\eps
]
\leq
\P [
| \xi^{0,1} |
\leq2 \eps
]
\leq
2 \eps\sigmab^{-1}
$
for all
$ x \in\R$
and all
$ \eps\in(0,\infty) $
show
%
%eA.45 #&#
%eA.46 #&#
%eA.47 #&#
%eA.48 #&#
%eA.49 #&#
\begin{eqnarray}
\label{eqdiffthetaeta} %
&& \P \bigl[ \llvert \theta_N -
\eta_N \rrvert \leq 4^{
( - 2^{ ( L_N - 1 ) }  )
} 2^{
{ ( L_N - 1 ) }/{ 4 }
}
T^{ - { 1 }/{ 4 } } N | \tilde{ \mathcal{F} }_{ \tilde{L}_N }^N \bigr]
\nonumber\\
& &\qquad\leq \P \biggl[ \exists k \in \biggl\{ 1, 2, \ldots, \frac{ N }{ 2^{ (L_N - 1) }
}
\biggr\} \dvtx \nonumber\\
&&\hspace*{23pt}\qquad{} \bigl| \bigl| \xi^{ L_N - 1, k } \bigr| - \eta_N\bigr |\leq
4^{
( - 2^{ ( L_N - 1 ) }
)
} 2^{
{ (L_N-1) }/{ 4 }
} T^{ - {1}/{4} } N | \tilde{ \mathcal{F}
}_{ \tilde{L}_N
}^N \biggr]
\nonumber
\\[-8pt]
\\[-8pt]
\nonumber
&&\qquad \leq \sum_{ k = 1 }^{{N}/{(2^{(L_N-1)})}} \P \bigl[ \bigl
\llvert \bigl|\xi^{ L_N - 1 , k }\bigr|-\eta_N \bigr\rrvert
\leq4^{ (-2^{(L_N-1)} )} 2^{
{ (L_N-1) }/{ 4 }
} T^{ - {1}/{4} } N | \tilde{
\mathcal{F}}_{\tilde{L}_N}^N \bigr]
\\
&&\qquad = \frac{N}{2^{(L_N-1)}} \cdot \P \bigl[ \bigl\llvert \bigl|\xi^{0,1}\bigr|-
\eta_N \bigr\rrvert \leq4^{ (-2^{(L_N-1)} )} 2^{
{ (L_N-1) }/{ 4 }
}
T^{ - {1}/{4} } N | \tilde{ \mathcal{F} }_{ \tilde{L}_N }^N \bigr]
\nonumber\\
&&\qquad \leq \frac{ N }{ 2^{ ( L_N - 1 ) }
} \cdot 2 \cdot 4^{
( - 2^{ ( L_N - 1 ) }
)
} 2^{
{ (L_N-1) }/{ 4 }
}
T^{ - {1}/{4} } N\cdot\sigmab^{-1} \leq \frac{
2 N^2
}{
\sigmab
T^{ { 1 }/{ 4 } }
2^{  ( 2^{ L_N }  ) }
}\nonumber
\end{eqnarray}
$ \mathbb{P} $-almost
surely
for all
$
N \in
\{ 2^1, 2^2, 2^3, \ldots\}
$.
Now we
apply
inequality~\eqref{eqdiffthetaeta}
to obtain
%
%eA.50 #&#
%eA.51 #&#
%eA.52 #&#
%eA.53 #&#
\begin{eqnarray}
\label{eqfirstestiA4} %
&& \mathbb{ P } \bigl[ A_{ N }^{ (4) }
\cap \bigl( A_{ N }^{ (2) } \bigr)^{ c } \cap \bigl(
A_{ N }^{ (1) } \bigr)^{ c } \bigr]
\nonumber\hspace*{-15pt}\\
&&\qquad \leq \mathbb{ P } \bigl[ \bigl\{ \llvert \eta_N -
\theta_N \rrvert \leq 4^{
( - 2^{ ( L_N - 1 ) }
)
} \eta_N \bigr\}
\cap \bigl\{ \eta_N< 2^{
{ (L_N-1) }/{ 4 }
} T^{ - {1}/{4} } N \bigr\} \cap
\bigl( A_{ N }^{ (1) } \bigr)^{ c } \bigr]
\nonumber\hspace*{-15pt}\\
&&\qquad \leq \mathbb{ P } \bigl[ \bigl\{ \llvert \eta_N -
\theta_N \rrvert \leq 4^{
( - 2^{ ( L_N - 1 ) }  )
} 2^{
{ (L_N-1) }/{ 4 }
}
T^{ - {1}/{4} } N \bigr\} \cap \bigl( A_{ N }^{ (1) }
\bigr)^{ c } \bigr]\hspace*{-15pt}
\\
&&\qquad = \mathbb{ E } \bigl[ \1_{
(
A_{ N }^{ (1) }
)^{ c }
} \cdot \mathbb{ P } \bigl[ \llvert
\eta_N-\theta_N\rrvert \leq4^{ (-2^{(L_N-1)} )}
2^{
{ (L_N-1) }/{ 4 }
} T^{ - {1}/{4} } N | \tilde{\mathcal{F}}_{\tilde{L}_N}^N
\bigr] \bigr]\nonumber\hspace*{-15pt}\\
&&\qquad \leq \mathbb{ E } \biggl[ \1_{
(
A_{ N }^{ (1) }
)^{ c }
} \cdot
\frac{
2 N^2
}{
\sigmab
T^{ { 1 }/{ 4 } }
2^{  ( 2^{ L_N }  ) }
} \biggr] \nonumber\hspace*{-15pt}
\end{eqnarray}
for all
$
N \in\{ 2^1, 2^2, 2^3, \ldots\}
$.
Next we observe that
%
%eA.54 #&#
\begin{eqnarray}
\label{eqtwotoLN} 2^{ L_N } &\geq& 2^{
\lfloor
2
\operatorname{ld} (
\sigmab^2
T^{
{1}/{2}
}
\operatorname{ln}(N)
)
\rfloor
} \geq 2^{
(
2
\operatorname{ld} (
\sigmab^2
T^{
{1}/{2}
}
\operatorname{ln}(N)
)
-1
)
}
\nonumber
\\[-8pt]
\\[-8pt]
\nonumber
&=& \tfrac{ 1 }{ 2 } \cdot 2^{
\operatorname{ld} (
(
\sigmab^2
T^{
{1}/{2}
}
\operatorname{ln}(N)
)^2
)
} = \tfrac12 \sigmab^4 T
\bigl(\ln(N) \bigr)^2
\end{eqnarray}
on
$
(
A_N^{ (1) }
)^{ c }
$
for all $N\in\{2^1,2^2,2^3,\ldots\}$.
Inserting~\eqref{eqtwotoLN} into~\eqref{eqfirstestiA4} results in
%
%eA.55 #&#
%eA.56 #&#
\begin{eqnarray}
\label{eqsecondestiA4} %
 \mathbb{ P } \bigl[ A_{ N }^{ (4) }
\cap \bigl( A_{ N }^{ (2) } \bigr)^{ c } \cap \bigl(
A_{ N }^{ (1) } \bigr)^{ c } \bigr] &\leq& \mathbb{ E }
\biggl[ \1_{
(
A_{ N }^{ (1) }
)^{c}
} \cdot \frac{ 4 N^2
}{
\sigmab
T^{ { 1 }/{ 4 } }
2^{
(
(1/2)\sigmab^4 T
( \ln(N)  )^2
)
}
} \biggr]
\nonumber\\
& \leq& \frac{ 2 N^2 }{
\sigmab T^{ { 1 }/{ 4 }
}
} \exp \biggl( - \frac{ \ln( 2 ) }{ 2 }
\sigmab^4 T \bigl(\ln(N) \bigr)^2 \biggr) \\
&=& 2 \bar{
\sigma}^{-1} T^{ - { 1 }/{ 4 } } N^{
(
2 -
\ln(2)\sigmab^4 T \ln(N) / 2
)
} \nonumber %
\end{eqnarray}
for all $N\in\{2^1,2^2,2^3,\ldots\}$.
Combining~\eqref{eqsecondestiA4},
Lemmas~\ref{lAN1}
and~\ref{lAN2} then
shows
%
%eA.57 #&#
%eA.58 #&#
%eA.59 #&#
\begin{eqnarray}
\sum_{ n = 1 }^{ \infty} \mathbb{ P }
\bigl[ A_{ 2^n }^{ (4) } \bigr] &=& \sum
_{ n = 1 }^{ \infty} \mathbb{ P } \bigl[ A_{ 2^n }^{ (4) }
\cap \bigl( A_{ 2^n }^{ (2) } \bigr)^{ c } \cap \bigl(
A_{ 2^n }^{ (1) } \bigr)^{ c } \bigr] \nonumber\\
&&{}+ \sum
_{ n = 1 }^{ \infty} \mathbb{ P } \bigl[ A_{ 2^n }^{ (4) }
\cap \bigl( \bigl( A_{ 2^n }^{ (2) } \bigr)^{ c } \cap
\bigl( A_{ 2^n }^{ (1) } \bigr)^{ c }
\bigr)^{ c } \bigr]
\nonumber\\
& \leq& \sum_{ n = 1 }^{ \infty} \mathbb{ P }
\bigl[ A_{ 2^n }^{ (4) } \cap \bigl( A_{ 2^n }^{ (2) }
\bigr)^{ c } \cap \bigl( A_{ 2^n }^{ (1) }
\bigr)^{ c } \bigr] + \sum_{ n = 1 }^{ \infty}
\mathbb{ P } \bigl[ A_{ 2^n }^{ (2) } \cup A_{ 2^n }^{ (1) }
\bigr]
\\
&\leq& \sum_{ N = 1 }^{ \infty} 4 \bar{
\sigma}^{-1} T^{ -{1}/{4} } N^{ (2 -
\ln(2)\sigmab^4 T \ln(N)/2
)} \nonumber\\
&&{}+ \sum
_{ n = 1 }^{ \infty} \mathbb{ P } \bigl[ A_{ 2^n }^{ (2) }
\bigr] + \sum_{ n = 1 }^{ \infty} \mathbb{ P }
\bigl[ A_{ 2^n }^{ (1) } \bigr] <\infty. \nonumber%
\end{eqnarray}
This completes the proof
of Lemma~\ref{lAN4}.
\end{pf}

We now present
the proof of
Lemma~\ref{lem_n1}.
It makes use of
Lemmas~\ref{lAN1}--\ref{lAN4}
above.
\begin{pf*}{Proof of
Lemma~\ref{lem_n1}}
Combining the subadditivity of
the probability measure
$ \mathbb{P} $ and
Lemmas~\ref{lAN1},~\ref{lAN2},~\ref{lAN3} and~\ref{lAN4}
shows
%
%eA.60 #&#
\begin{eqnarray}
&&\sum_{ n = 1 }^{ \infty} \mathbb{P}
\bigl[ A_{ 2^n }^{ (1) } \cup A_{ 2^n }^{ (2) }
\cup A_{ 2^n }^{ (3) } \cup A_{ 2^n }^{ (4) }
\bigr]
\nonumber
\\[-8pt]
\\[-8pt]
\nonumber
&&\qquad\leq \sum_{ n = 1 }^{ \infty} \mathbb{P}
\bigl[ A_{ 2^n }^{ (1) } \bigr] + \sum
_{ n = 1 }^{ \infty} \mathbb{P} \bigl[ A_{ 2^n }^{ (2) }
\bigr] + \sum_{ n = 1 }^{ \infty} \mathbb{P} \bigl[
A_{ 2^n }^{ (3) } \bigr] + \sum_{ n = 1 }^{ \infty}
\mathbb{P} \bigl[ A_{ 2^n }^{ (4) } \bigr] < \infty. %
\end{eqnarray}
The lemma of Borel--Cantelli
(e.g., Theorem~2.7 in~\cite{klenke2008})
therefore implies
%
%eA.61 #&#
\begin{equation}
\mathbb{P} \Bigl[ \limsup_{ n \rightarrow\infty} \bigl( A^{ (1) }_{ 2^n }
\cup A^{ (2) }_{ 2^n } \cup A^{ (3) }_{ 2^n }
\cup A^{ (4) }_{ 2^n } \bigr) \Bigr] = 0 .
\end{equation}
Hence, we obtain
%
%eA.62 #&#
%eA.63 #&#
%eA.64 #&#
\begin{eqnarray}
\mathbb{P} [ N_1 < \infty ] &=& \mathbb{P} \Biggl[
\Biggl\{ \omega\in\Omega\dvtx \exists n \in \bigl\{ N_0,
2^1 N_0, 2^2 N_0, \ldots
\bigr\} \dvtx\nonumber\\
&&\hspace*{18pt}{} \forall m \in \bigl\{ n, 2^1n, 2^2n, \ldots
\bigr\} \dvtx \omega \notin \bigcup_{ i = 1 }^4
A^{ (i) }_{ m } \Biggr\} \Biggr]
\nonumber\\
& =& \mathbb{P} \bigl[ \bigl\{ \omega\in\Omega\dvtx \exists n \in\mathbb{N}
 \dvtx\forall m \in\{ n, n+1, \ldots\} \dvtx\\
 &&\hspace*{50pt}{} \omega\notin A^{ (1) }_{ 2^m }
\cup A^{ (2) }_{ 2^m } \cup A^{ (3) }_{ 2^m }
\cup A^{ (4) }_{ 2^m } \bigr\} \bigr]
\nonumber\\
& =& \mathbb{P} \Bigl[ \liminf_{ n \rightarrow\infty} \bigl( A^{ (1) }_{ 2^n }
\cup A^{ (2) }_{ 2^n } \cup A^{ (3) }_{ 2^n }
\cup A^{ (4) }_{ 2^n } \bigr)^{ c } \Bigr] = 1 .\nonumber
\end{eqnarray}
This completes the proof
of Lemma~\ref{lem_n1}.
\end{pf*}
\end{appendix}

\section*{Acknowledgments}

We would
like to express our sincere
gratitude to
Weinan~E,
Klaus Ritter,
Andrew M.~Stuart,
Jan van Neerven
and
Konstantinos Zygalakis
for their
very helpful advice.

%

%
% imsref loaded by akundreckaite, 2012-12-06 14:08\dvtx 36
%
% imsref loaded by akundreckaite, 2012-12-10 09:07:28

%suskaldyti doi

\printaddresses

\end{document}